\documentclass[preprint,3p,sort&compress]{elsarticle}

\usepackage{afterpage}        %% Do stuff when starting a new page
\usepackage{amssymb}          %% Math symbols
\usepackage{booktabs}         %% Nice table (mid lines)
\usepackage{commath}          %% Math symbols for abs and norm
\usepackage[overload]{empheq} %% System of equations with align and others
\usepackage[inline]{enumitem} %% Inline enumerate
\usepackage[pdftex]{hyperref} %% Hyperlinks
\usepackage{subcaption}       %% Subfloat
\usepackage{xspace}           %% Smart spaces

\newlist{myenum*}{enumerate*}{1}
\setlist[myenum*]{label=\itshape\roman*\upshape)}

\newcommand{\ie}{\textit{i.e.}\xspace}
\newcommand{\eg}{\textit{e.g.}\xspace}
\renewcommand{\vec}[1]{\mathbf{#1}}
\newcommand{\imag}{\jmath}
\newcommand{\dd}{\mathrm{d}}
\newcommand{\deriv}[2]{\frac{\partial#1}{\partial#2}}
\newcommand{\derivv}[2]{\frac{\partial^2#1}{\partial#2}}
\newcommand{\tc}[2]{\uppercase{#1}\textsuperscript{#2}}
\newcommand{\Tab}[1]{Table~\ref{#1}}
\newcommand{\Fig}[1]{Figure~\ref{#1}}
\newcommand{\Figs}[1]{Figures~\ref{#1}}
\newcommand{\Code}[1]{\texttt{#1}}
\renewcommand{\epsilon}{\varepsilon}
\renewcommand{\url}[1]{\href{https://#1}{\Code{#1}}}

\DeclareMathOperator{\Grad}{\vec{grad}}
\DeclareMathOperator{\Div}{div}
\DeclareMathOperator{\Trsm}{\mathcal{S}}
\DeclareMathOperator{\Ident}{\mathcal{I}}
\DeclareMathOperator{\OpA}{\mathcal{A}}
\DeclareMathOperator{\Error}{\mathcal{E}}
\DeclareMathOperator{\Op}{Op}
\DeclareMathOperator{\DtN}{DtN}
\DeclareMathOperator{\Eig}{Eig}

\newcommand{\Sij}[2]{\Trsm_{ij}^{\text{#2, #1}}}
\newcommand{\Sji}[2]{\Trsm_{ji}^{\text{#2, #1}}}
\newcommand{\Lij}[2]{\lambda_{ij}^{\text{#2, #1}}}
\newcommand{\Lji}[2]{\lambda_{ji}^{\text{#2, #1}}}

\begin{document}
\begin{frontmatter}
  \title{Transmission operators for the non-overlapping Schwarz method
    for solving Helmholtz problems in rectangular cavities}
  \author[1]{Nicolas Marsic\corref{cor1}}\ead{marsic@temf.tu-darmstadt.de}
  \author[2]{Christophe Geuzaine}\ead{cgeuzaine@uliege.be}
  \author[1]{Herbert De Gersem}\ead{degersem@temf.tu-darmstadt.de}
  \address[1]{Technische Universit\"at Darmstadt,
    Institute for Accelerator Science and Electromagnetic Fields (TEMF),
    D-64289 Darmstadt, Germany}
  \address[2]{Universit\'e de Li\`ege, Institut Montefiore B28,
    B-4000 Li\`ege, Belgium}
  \cortext[cor1]{Corresponding author}

  \begin{abstract}
    In this paper we discuss different transmission operators
    for the non-overlapping Schwarz method
    which are suited for solving the time-harmonic Helmholtz equation
    in cavities
    (\ie closed domains which do not feature an outgoing wave condition).
    Such problems are heavily impacted by back-propagating waves
    which are often neglected
    when devising optimized transmission operators for the Schwarz method.
    This work explores new operators taking into account
    those back-propagating waves and compares them
    with well-established operators neglecting these contributions.
    Notably, this paper focuses on the case of rectangular cavities,
    as the optimal (non-local) transmission operator can be easily determined.
    Nonetheless, deviations from this ideal geometry are considered as well.
    In particular, computations of the acoustic noise
    in a three-dimensional model of the helium vessel of a beamline cryostat
    with optimized Schwarz schemes are discussed.
    Those computations show a reduction of $46\%$ in the iteration count,
    when comparing an operator optimized for cavities with
    those optimized for unbounded problems.
  \end{abstract}

\end{frontmatter}

\section{Introduction}
\label{sec:intro}
It is well known that large-scale time-harmonic Helmholtz problems
are hard to solve because of
\begin{myenum*}
\item the pollution effect~\cite{Ihlenburg1995} and
\item the indefiniteness of the discretized operator~\cite{Ernst2012}.
\end{myenum*}
While the pollution effect can be alleviated
by using higher order discretization schemes~\cite{Ihlenburg1997},
the indefiniteness is an intrinsic property of time-harmonic wave problems,
at least with standard variational formulations~\cite{Moiola2014, Diwan2019},
and significantly limits the performance of classical iterative solvers,
such as the generalized minimal residual method (GMRES) for instance.
Of course, as an alternative to iterative algorithms,
direct solvers can be used.
However, because of the fill-in effect,
whose minimization is known to be a NP-complete problem~\cite{Yannakakis1981},
the amount of memory needed to treat large-scale systems
can become prohibitively high (see for instance~\cite{Marsic2018b}).

As an alternative to direct and (unpreconditioned) iterative methods
for solving large-scale, high-frequency time-harmonic Helmholtz problems,
domain decomposition (DD) algorithms, and
optimized Schwarz (OS) techniques~\cite{Despres1990, Boubendir2007,
  Gander2002, Boubendir2012} in particular,
have attracted a lot of attention during the last decades.
The key idea thereof is:
\begin{myenum*}
\item to decompose the computational domain
  into (possibly overlapping) subdomains, creating thus new subproblems,
\item to solve each subproblem \emph{independently},
\item to exchange data at the interfaces between the subdomains
  via an appropriate \emph{transmission operator} and
\item to repeat this ``solve and exchange'' procedure
  until convergence of the solution.
\end{myenum*}
Since all subproblems are solved independently,
domain decomposition methods are parallel by nature\footnote{It is also possible
  to solve the subproblems sequentially
  and to exchange data after each single solve.
  This family of DD methods are often referred to as sweeping algorithms,
  and offer some advantages, notably in terms of iteration count,
  which will not be further discussed in this work.
  More details can be found for instance in~\cite{Vion2014a, Gander2019}.}
and are thus very well suited for the treatment of large-scale problems.
Furthermore, as the subproblems are of reduced size, direct solvers can be used.
Let also note that DD methods are rarely used as stand-alone solvers,
but most of the time as a \emph{preconditioner} for a Krylov subspace method
such as GMRES.
The design of such preconditioners for time-harmonic Helmholtz problems
remains an active and challenging topic~\cite{Gander2019}.

The convergence rate of an OS scheme strongly depends on its
transmission operator.
It is well known that the optimal operator is the
Dirichlet-to-Neumann ($\DtN$) map
at the interface between two subdomains~\cite{Dolean2015a}
(\ie the operator relating the trace of the unknown field
to its normal derivative at a given interface).
However, the $\DtN$ map is rarely employed,
as it is a \emph{non-local} operator which leads to
a numerically expensive scheme.
Instead, in practice, \emph{local approximations} of the $\DtN$ map are used,
which lead to many different computational schemes~\cite{Despres1990,
  Gander2002, Boubendir2007, Boubendir2012}.
To the best of our knowledge,
those OS techniques share a common drawback:
they ignore the impact of \emph{back-propagating waves}.
While this assumption is legitimate in many cases
(antenna arrays~\cite{Peng2011},
medical imaging reconstruction~\cite{Tournier2017}
or photonic waveguides~\cite{Marsic2016a} just to cite a few),
it becomes questionable when the geometry allows resonances
(even if the source does not oscillate exactly at a resonance frequency),
as found for instance in lasers~\cite{Saleh2007},
accelerator cavities~\cite{Ko2006}
or quantum electrodynamic devices~\cite{Marsic2018b}.

The objective of this work is to develop new transmission conditions
taking into account the effect of back-propagating waves, and
to compare them with well-established operators neglecting these contributions.
To this end, we will study a rectangular cavity,
determine the $\DtN$ map and localize it by following different strategies.
We will then apply the resulting new transmission operators
to more general geometries.
This paper is organized as follows.
In sections~\ref{sec:problem} and~\ref{sec:S:dirichlet} the
model problem with Dirichlet boundary conditions
and the associated $\DtN$ map are presented
for both overlapping and non-overlapping decompositions.
New transmission operators are afterwards presented in section~\ref{sec:tc}
and generalized (\ie multiple subdomains and Neumann boundary conditions)
in section~\ref{sec:generalization}.
This is followed by section~\ref{sec:open}
showing a comparison with the classical $\DtN$ map
related to unbounded problems and the use of transmission operators
optimized for unbounded problems as an approximation
of the \emph{cavity} $\DtN$ map.
The new transmission operators are then validated and compared
with numerical experiments involving the reference rectangular cavity
in section~\ref{sec:validation}.
The case of geometries deviating from this last model problem
is further discussed in section~\ref{sec:sensitivity}
and an engineering problem involving a model of the helium vessel
of a beamline cryostat is analyzed in section~\ref{sec:cryo}.
Finally, conclusions and final remarks are drawn
in section~\ref{sec:conclusion}.

\section{Model problem and Schwarz domain decomposition method}
\label{sec:problem}
\begin{figure}[ht]
  \centering
  \begin{subfigure}[b]{0.49\textwidth}
    \centering
    \includegraphics{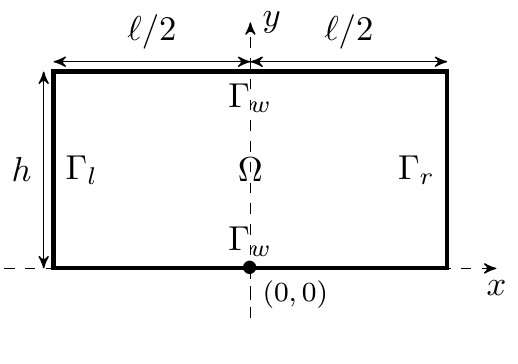}
    \caption{Domain $\Omega$\dots}
    \label{fig:domain:omega}
  \end{subfigure}
  \begin{subfigure}[b]{0.49\textwidth}
    \includegraphics{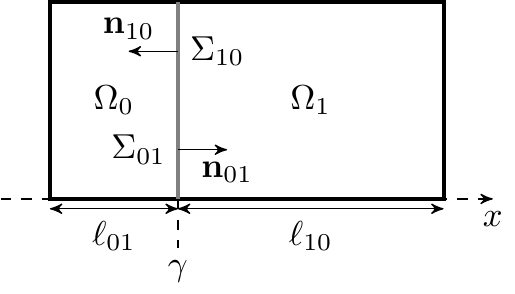}
    \caption{\dots and its decomposition into
      $\Omega_1$ and $\Omega_2$.}
    \label{fig:domain:dd}
  \end{subfigure}
  \caption{Reference computational domain.}
  \label{fig:domain}
\end{figure}

Let $\Omega$ be the two-dimensional domain $[-\ell/2, +\ell/2]\times[0, h]$
depicted in \Fig{fig:domain:omega},
and let $\Gamma=\Gamma_l\cup\Gamma_w\cup\Gamma_r$ be its boundary.
This domain is separated into two non-overlapping subdomains
$\Omega_0=[-\ell/2,\gamma(t)]\times[0, h]$ and
$\Omega_1=[\gamma(t),+\ell/2]\times[0, h]$,
where $\gamma(t) = t\ell - \ell/2$ with $t\in[0, 1]$ a fixed parameter
controlling the position of the interface shared by the two subdomains,
as shown in \Fig{fig:domain:dd}.
In addition,
the resulting subdomains have a length of $\ell_{01}$ (for $\Omega_0$)
and $\ell_{10}$ (for $\Omega_1$) respectively
and $\ell = \ell_{01}+\ell_{10}$.
This splitting has introduced a new artificial boundary on each subdomain:
we denote by $\Sigma_{01}$ the artificial boundary of $\Omega_0$ and
by $\Sigma_{10}$ the artificial boundary of $\Omega_1$.
Furthermore, $\vec{n}_{ij}$ denotes the outwardly oriented unit vector
normal to $\Sigma_{ij}$.

Let us solve the following Helmholtz problem on $\Omega$:
\begin{subequations}
  \label{eq:helmholtz}
  \begin{align}[left = \empheqlbrace]
    \Div{\Grad{p}} + k^2p & = g & \text{on}~\Omega,\label{eq:helmholtz:omega}\\
    p                     & = 0 & \text{on}~\Gamma,\label{eq:helmholtz:gamma}
  \end{align}
\end{subequations}
where $p(x, y)$ is the unknown function, $g(x, y)$ is a known source term
and $k\in\mathbb{R}$ is the fixed wavenumber of the Helmholtz problem.
Because of its boundary condition, it is obvious that~\eqref{eq:helmholtz}
models a \emph{cavity problem} exhibiting both forward-
\emph{and back-propagating waves}.
It is important to stress that for this problem to be well defined,
we must assume that $k^2$ is not an eigenvalue of~\eqref{eq:helmholtz}.

Let us now set up the following optimized non-overlapping Schwarz iterative scheme,
indexed by $n$, to solve the cavity Helmholtz problem~\eqref{eq:helmholtz}:
\begin{subequations}
  \label{eq:ddm}
  \begin{align}[left = \empheqlbrace]
    \Div{\Grad{p_0^{n+1}}} + k^2p_0^{n+1}
    & = g
    & \text{on}~\Omega_0, \\
    +\vec{n}_{01}\cdot\Grad{p_0^{n+1}} + \Trsm_{01}{(p_0^{n+1})}
    & = +\vec{n}_{01}\cdot\Grad{p_1^{n}} + \Trsm_{01}{(p_1^{n})}
    & \text{on}~\Sigma_{01}, \\
    p_0^{n+1}
    & = 0
    & \text{on}~\Gamma,  \label{eq:ddm:bc0} \\
    \nonumber\\
    \Div{\Grad{p_1^{n+1}}} + k^2p_1^{n+1}
    & = g
    & \text{on}~\Omega_1, \\
    -\vec{n}_{01}\cdot\Grad{p_1^{n+1}} + \Trsm_{10}{(p_1^{n+1})}
    & = -\vec{n}_{01}\cdot\Grad{p_0^{n}} + \Trsm_{10}{(p_0^{n})}
    & \text{on}~\Sigma_{10}, \\
    p_1^{n+1}
    & = 0
    & \text{on}~\Gamma, \label{eq:ddm:bc1}
  \end{align}
\end{subequations}
where $\Trsm_{ij}$ is the transmission operator
of the optimized Schwarz algorithm at $\Sigma_{ij}$
and $p_i^n(x, y)$ is the solution at iteration $n$ on domain $\Omega_i$.
Let us stress that, since the subdomains do not overlap,
the following holds true: $\vec{n}_{01} = -\vec{n}_{10}$.
Once the Schwarz algorithm has converged,
the solution $p(x, y)$ of the original problem~\eqref{eq:helmholtz}
is recovered by concatenating the solutions $p_0(x, y)$ and $p_1(x, y)$.

In practice, let us note that the above fixed-point scheme is usually recast
into the linear system~\cite{Dolean2015a}:
\begin{equation}
  \label{eq:ddm:lin}
  (\Ident-\OpA)\vec{d} = \vec{b},
\end{equation}
where one application of the operator $\OpA$ amounts to one iteration of the
fixed-point method with \emph{homogeneous} Dirichlet boundary conditions,
where $\Ident$ is the identity operator,
where $\vec{d}$ concatenates all $\vec{n}\cdot\Grad{p}+\Trsm(p)$
at the interface between the subdomains
and where the right hand side vector $\vec{b}$ results from the
\emph{non-homogeneous} Dirichlet boundary conditions.
This linear system can then be solved with a matrix-free Krylov subspace method
such as GMRES.

\section{Optimal transmission operator for the rectangular cavity problem with
  homogeneous Dirichlet boundary conditions}
\label{sec:S:dirichlet}
In this section, we will first determine the optimal transmission operators
at $\Sigma_{01}$ and $\Sigma_{10}$ of the Schwarz scheme~\eqref{eq:ddm}
involving the \emph{non-overlapping} subdomains in \Fig{fig:domain:dd}.
While this work focuses on non-overlapping decompositions,
the impact of an overlap on the optimal transmission operator
is also discussed in this section.

\subsection{Non-overlapping case}
In order to further simplify the problem at hand,
let us now assume that the source term $g$ is zero.
Obviously, by not imposing a source in our problem
the solution $p(x,y)$ is trivially $p=0$ since $k$ is not an eigenvalue.
This however does not jeopardize the generality of
the results derived in this section.

Let us start by taking
the sine Fourier series of $p_i^n(x,y)$ along the $y$-axis:
\begin{equation}
  \label{eq:fourier}
  p_i^n(x,y) = \sum_{s\in\mathbb{S}}\widehat{p}_i^n(x, s)\sin(sy),
\end{equation}
where the functions $\widehat{p}_i^n(x, s)$ are the Fourier coefficients
and where $s$ is the Fourier variable, whose values are restricted to the set
\begin{equation}
  \label{eq:fourier:set}
  \mathbb{S} =
  \Big\{s\in\mathbb{R}~\Big|~s=m\frac{\pi}{h}, \forall m\in\mathbb{N}_0\Big\}.
\end{equation}
Indeed, by restricting $s$ to the set $\mathbb{S}$, the boundary conditions
\begin{align*}[left = \empheqlbrace]
  p^n_i(x, 0) & = 0 & \forall x\in\left[-\frac{\ell}{2}, +\frac{\ell}{2}\right],
  \\
  p^n_i(x, h) & = 0 & \forall x\in\left[-\frac{\ell}{2}, +\frac{\ell}{2}\right],
\end{align*}
are automatically satisfied.
Then, by exploiting decomposition~\eqref{eq:fourier},
the partial differential equation~\eqref{eq:ddm}
becomes the following ordinary differential equation:
\begin{subequations}
  \label{eq:ddm:fourier}
  \begin{align}[left = \empheqlbrace]
    \label{eq:ddm:fourier:ode:0}
    \derivv{\widehat{p}_0^{n+1}}{x^2} + (k^2-s^2)\widehat{p}_0^{n+1}
    & = 0
    & \forall x\in\left[-\frac{\ell}{2}, \gamma\right]
      \text{and}~\forall s\in\mathbb{S}, \\
    \label{eq:ddm:fourier:trsm:0}
    +\deriv{\widehat{p}_0^{n+1}}{x}+\lambda_{01}\,\widehat{p}_0^{n+1}
    & = +\deriv{\widehat{p}_1^{n}}{x}+\lambda_{01}\,\widehat{p}_1^{n}
    & \text{on}~x=\gamma~\text{and}~\forall s\in\mathbb{S}, \\
    \label{eq:ddm:fourier:bnd:0}
    \widehat{p}_0^{n+1}
    & = 0
    & \text{on}~x=-\frac{\ell}{2}~\text{and}~\forall s\in\mathbb{S}, \\
    \nonumber \\
    \label{eq:ddm:fourier:ode:1}
    \derivv{\widehat{p}_1^{n+1}}{x^2} + (k^2-s^2)\widehat{p}_1^{n+1}
    & = 0
    & \forall x\in\left[\gamma, +\frac{\ell}{2}\right]
      \text{and}~\forall s\in\mathbb{S}, \\
    \label{eq:ddm:fourier:trsm:1}
    -\deriv{\widehat{p}_1^{n+1}}{x}+\lambda_{10}\,\widehat{p}_1^{n+1}
    & = -\deriv{\widehat{p}_0^{n}}{x}+\lambda_{10}\,\widehat{p}_0^{n}
    & \text{on}~x=\gamma~\text{and}~\forall s\in\mathbb{S}, \\
    \label{eq:ddm:fourier:bnd:1}
    \widehat{p}_1^{n+1}
    & = 0
    & \text{on}~x=+\frac{\ell}{2}~\text{and}~\forall s\in\mathbb{S},
  \end{align}
\end{subequations}
where $\lambda_{ij}$ is the Fourier symbol of $\Trsm_{ij}$.
Furthermore, and for simplicity, let us define $P_i^n(s)$ as:
\begin{equation}
  \label{eq:ddm:fourier:bnd:x0}
  P_i^n(s) = \widehat{p}_i^n(0, s).
\end{equation}

In order to find the best symbol $\lambda_{ij}$,
we need to determine the convergence radius
of the iterative scheme~\eqref{eq:ddm:fourier}.
This objective can be achieved
by following the strategy discussed in~\cite{Gander2002},
which can be summarized as follows.
\begin{enumerate}
\item Determine the solutions of~\eqref{eq:ddm:fourier:ode:0}
  and~\eqref{eq:ddm:fourier:ode:1}
  with the boundary conditions~\eqref{eq:ddm:fourier:bnd:0} and
  \eqref{eq:ddm:fourier:bnd:1} and the definition~\eqref{eq:ddm:fourier:bnd:x0}.

\item Compute $\displaystyle\deriv{\widehat{p}_i^n(x, s)}{x}$ at $x=\gamma$
  from the solutions $\widehat{p}_i^n(x,s)$ found in the previous step.

\item The convergence radius is obtained
  by simplifying the transmission conditions~\eqref{eq:ddm:fourier:trsm:0}
  and~\eqref{eq:ddm:fourier:trsm:1} with the expressions found above.
\end{enumerate}
By following this approach,
it can be shown (see~\ref{sec:dtn}) that the convergence radius
$\rho$ of~\eqref{eq:ddm:fourier} satisfies\footnote{In what follows,
  we distinguish the cavity and unbounded contexts with
  the superscripts ``c'' and ``u'' respectively.}:
\begin{equation}
  \label{eq:fourier:rho}
  \rho^2(s)
  =
  \frac{\lambda_{01}(s)-\lambda_{01}^{\text{c, dtn}}(s)}
       {\lambda_{01}(s)+\lambda_{10}^{\text{c, dtn}}(s)}
  \frac{\lambda_{10}(s)-\lambda_{10}^{\text{c, dtn}}(s)}
       {\lambda_{10}(s)+\lambda_{01}^{\text{c, dtn}}(s)}
  =
  \frac{n_{01}(s)\,n_{10}(s)}{d_{01}(s)\,d_{10}(s)},
\end{equation}
where
\begin{subequations}
  \label{eq:close:lambda:dtn}
  \begin{align}[left = {\Lij{dtn}{c}(s)=\empheqlbrace}]
    & \sqrt{k^2-s^2}
      \cot\mathopen{}\left[\ell_{ji}\sqrt{k^2-s^2}\mathclose{}\right]
    & \text{if}~s^2<k^2,\\
    & 1/\ell_{ji}
    & \text{if}~s^2=k^2,\\
    & \sqrt{s^2-k^2}
      \coth\mathopen{}\left[\ell_{ji}\sqrt{s^2-k^2}\mathclose{}\right]
    & \text{if}~s^2>k^2.
  \end{align}
\end{subequations}
The best transmission operator $\Sij{dtn}{c}$,
that is the Dirichlet-to-Neumann map, is thus
\begin{equation}
  \label{eq:close:S:dtn}
  \Sij{dtn}{c} = \Op\mathopen{}\Bigg\{\Lij{dtn}{c}\mathclose{}\Bigg\}.
\end{equation}

\subsection{Overlapping case}
\begin{figure}[ht]
  \centering
  \includegraphics{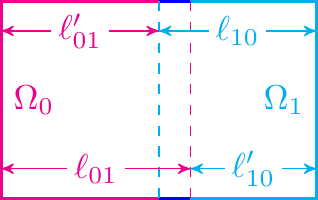}
  \caption{Overlapping partitioning of the computational domain.}
  \label{fig:domain:overlap}
\end{figure}

Let us now assume a partitioning of the domain in \Fig{fig:domain:omega} into
two overlapping rectangles, as shown in \Fig{fig:domain:overlap}.
As suggested by this figure, we define $\ell_{01}$ (resp. $\ell_{10}$)
as the length of $\Omega_0$ (resp. $\Omega_1$) including the overlap
and $\ell_{01}^\prime$ (resp. $\ell_{10}^\prime$)
as the length of $\Omega_0$ (resp. $\Omega_1$) without the overlap.
By following the same strategy as in the non-overlapping case,
but taking into account that $\Sigma_{01}$ and $\Sigma_{10}$
have now \emph{different} locations,
the convergence radius $\rho^\text{overlap}$
of the overlapping variant of~\eqref{eq:ddm:fourier} reads
(see~\ref{sec:dtn:overlap}):
\begin{subequations}
  \label{eq:fourier:rho:overlap}
  \begin{align}[left = {(\rho^\text{overlap})^2 = \empheqlbrace}]
    & \left(
        \frac{\lambda_{01}-\alpha\cot(\alpha\ell_{10}^\prime)}
             {\lambda_{01}+\alpha\cot(\alpha\ell_{01})}
        \frac{\lambda_{10}-\alpha\cot(\alpha\ell_{01}^\prime)}
             {\lambda_{10}+\alpha\cot(\alpha\ell_{10})}
      \right)
      \left(
        \frac{\sin(\alpha\ell_{10}^\prime)}{\sin(\alpha\ell_{10})}
        \frac{\sin(\alpha\ell_{01}^\prime)}{\sin(\alpha\ell_{01})}
      \right)
    & \text{if}~s^2<k^2,\label{eq:fourier:rho:overlap:prop}\\
    & \left(
        \frac{\lambda_{01}-(\ell_{10}^\prime)^{-1}}
             {\lambda_{01}+(\ell_{01})^{-1}}
        \frac{\lambda_{10}-(\ell_{01}^\prime)^{-1}}
             {\lambda_{10}+(\ell_{10})^{-1}}
      \right)
      \left(
        \frac{\ell_{10}^\prime}{\ell_{10}}
        \frac{\ell_{01}^\prime}{\ell_{01}}
      \right)
    & \text{if}~s^2=k^2,\\
    & \left(
        \frac{\lambda_{01}-\alpha\coth(\alpha\ell_{10}^\prime)}
             {\lambda_{01}+\alpha\coth(\alpha\ell_{01})}
        \frac{\lambda_{10}-\alpha\coth(\alpha\ell_{01}^\prime)}
             {\lambda_{10}+\alpha\coth(\alpha\ell_{10})}
      \right)
      \left(
        \frac{\sinh(\alpha\ell_{10}^\prime)}{\sinh(\alpha\ell_{10})}
        \frac{\sinh(\alpha\ell_{01}^\prime)}{\sinh(\alpha\ell_{01})}
      \right)
    & \text{if}~s^2>k^2,\label{eq:fourier:rho:overlap:evan}
  \end{align}
\end{subequations}
where
\begin{subequations}
  \label{eq:alpha}
  \begin{align}[left = {\alpha(s) = \empheqlbrace}]
    & -\imag\sqrt{k^2-s^2} & \text{if}~s^2<k^2,\\
    & 0                    & \text{if}~s^2=k^2,\\
    & \sqrt{s^2-k^2}       & \text{if}~s^2>k^2.
  \end{align}
\end{subequations}
A few conclusions can be drawn from the above expressions.
Firstly, it is clear from~\eqref{eq:fourier:rho:overlap} that
the optimal transmission operator writes
\begin{subequations}
  \label{eq:close:lambda:dtn:overlap}
  \begin{align}[left = {\Lij{dtn, overlap}{c}(s)=\empheqlbrace}]
    & \sqrt{k^2-s^2}
      \cot\mathopen{}\left[\ell_{ji}^\prime\sqrt{k^2-s^2}\mathclose{}\right]
    & \text{if}~s^2<k^2,\\
    & 1/\ell_{ji}^\prime
    & \text{if}~s^2=k^2,\\
    & \sqrt{s^2-k^2}
      \coth\mathopen{}\left[\ell_{ji}^\prime\sqrt{s^2-k^2}\mathclose{}\right]
    & \text{if}~s^2>k^2,
  \end{align}
\end{subequations}
that is $\Lij{dtn, overlap}{c}(s)$ is equal to $\Lij{dtn}{c}(s)$,
up to the substitution $\ell_{ji}\rightarrow\ell_{ji}^\prime$.
Secondly, it is easy to notice that when there is no overlap,
it holds that $\ell_{ji}=\ell_{ji}^\prime$
and the non-overlapping case is recovered.
Thirdly and most importantly,
it is obvious that an overlap does \emph{not necessarily} improve
the convergence radius $\rho^\text{overlap}$
\emph{unlike} for the unbounded case~\cite{Dolean2015a},
since when $s^2<k^2$ (\ie propagating modes)
the term involving the $\sin$ functions in~\eqref{eq:fourier:rho:overlap:prop}
exhibits poles and oscillates with respect to the size of the overlap.
Nonetheless,
the term with the $\sinh$ functions in~\eqref{eq:fourier:rho:overlap:evan}
introduces a damping proportional to the overlap,
as in the unbounded case, but only when $s^2>k^2$ (\ie evanescent modes).

\section{Some local transmission operators for the cavity problem}
\label{sec:tc}
In this section,
we discuss some local transmission operators
based on the Fourier symbol~\eqref{eq:close:lambda:dtn}.

\subsection{Zeroth-order transmission condition}
A zeroth-order transmission condition (\tc{oo0}{c}) can easily be constructed by
approximating the symbol of the $\DtN$ map
with the constant term of its Taylor series expansion around $s=0$.
For the considered cavity setting, we obtain:
\begin{equation}
  \label{eq:close:lambda:oo0}
  \Lij{dtn}{c}(s) \approx \Lij{oo0}{c}(s) = k\cot(k\ell_{ji})
\end{equation}
and the \tc{oo0}{c} transmission condition reads
\begin{equation}
  \label{eq:close:S:oo0}
  \Sij{oo0}{c} = k\cot(k\ell_{ji}).
\end{equation}

As such, this operator exhibits a rather poor behavior.
Indeed,
the \emph{denominator} of the convergence radius~\eqref{eq:fourier:rho}
involves terms of the form
\begin{equation}
  \label{eq:dij}
  d_{ij}(s) = \lambda_{ij}(s) + \Lji{dtn}{c}(s),
\end{equation}
which can \emph{change their sign multiple times} when $s^2<k^2$,
since $\Lji{dtn}{c}$ is nothing but a cotangent and
$\Lij{oo0}{c}$ is \emph{constant}.
Consequently, it is possible that $d_{ij} \approx 0$ for some $s\in\mathbb{S}$,
leading to a very large convergence radius.
In the worst case scenario, one can also have $d_{ij} = 0$
for some $s\in\mathbb{S}$ and the problem becomes ill-posed.
Regularization procedures for preventing this behavior
are further discussed in sections~\ref{sec:tc:reg} and~\ref{sec:tc:mixed}.

\subsection{Truncated Mittag-Leffler expansion based transmission condition}
In order to improve the performance of the above
\tc{oo0}{c} transmission condition for cavity problem,
we need a condition whose symbol is a better approximation of $\Lij{dtn}{c}$.
To this end, an option is to exploit
the Mittag-Leffler~\cite{Ahlfors1979} expansion of $\cot(z)$
according to its poles, leading to the following
partial fraction decomposition~\cite{Ahlfors1979}:
\begin{equation}
  \label{eq:closed:ml:cot}
  z\cot(az) =
  \frac{1}{a}\left(
    1 + 2\sum_{n=1}^\infty\frac{z^2}{z^2-\left(\frac{n\pi}{a}\right)^2}
  \right),
\end{equation}
that can be exploited to expand the symbol of the $\DtN$ map as
\begin{equation}
  \label{eq:closed:ml:opt}
  \Lij{dtn}{c}(s) =
  \frac{1}{\ell_{ji}}\left(
    1 + 2\sum_{n=1}^\infty\frac{1-\frac{s^2}{k^2}}
                             {1-\left(\frac{n\pi}{k\ell_{ji}}\right)^2
                               -\frac{s^2}{k^2}}
  \right).
\end{equation}
This symbol can hence be localized
by \emph{truncating} the series up to the $N$\textsuperscript{th} term,
enabling us to form a $N$-term truncated Mittag-Leffler expansion based
(\tc{ml}{c}) transmission condition:
\begin{equation}
  \label{eq:closed:ml:S}
  \Sij{ml($N$)}{c} =
  \frac{1}{\ell_{ji}}\left\{
    1 + 2\sum_{n=1}^N\Bigg[1+\frac{\Div_{\Sigma_{ij}}\Grad_{\Sigma_{ij}}}{k^2}\Bigg]
                     \Bigg[1-\Bigg(\frac{n\pi}{k\ell_{ji}}\Bigg)^2
                            +\frac{\Div_{\Sigma_{ij}}\Grad_{\Sigma_{ij}}}{k^2}\Bigg]^{-1}
  \right\}.
\end{equation}

As increasing the number of terms in this expansion
makes the transmission condition arbitrarily precise,
the poor convergence rate of the previous \tc{oo0}{c} operator can be alleviated.
This comes however at the cost of $N$ auxiliary computations to account for the
$N$ inverse operations appearing in~\eqref{eq:closed:ml:S}~\cite{Boubendir2012}.
For this reason,
it is desirable to devise an approximation of $\Lij{dtn}{c}(s)$
with a limited number of auxiliary terms.
Let us also mention that while the above operator has been developed
for a one-dimensional interface, it can be used as well
with a two-dimensional one, as shown in section~\ref{sec:cryo} for instance.

\subsection{Pad\'e approximant based transmission condition}
A Pad\'e rational approximation exhibits usually a good convergence rate
with respect to its order $[M/N]$,
where $M$ (resp. $N$) denotes the order of the numerator (resp. denominator).
One can construct it by exploiting
the continued fraction expansion of the function to approach~\cite{Baker1996}.
By taking the reciprocal of the continued fraction expansion
of the tangent function~\cite{Oldham2009}, we have:
\begin{equation}
  \label{eq:closed:pade:cot:cfrac}
  z\cot(z)
  = 1 - \cfrac{z^2}{3-\cfrac{z^2}{7-\cfrac{z^2}{9-\dots}}}
  = b_0 + \cfrac{a_1}{b_1+\cfrac{a_2}{b_2+\cfrac{a_3}{b_3+\dots}}},
\end{equation}
where
\begin{subequations}
  \begin{align*}[left = \empheqlbrace]
    b_n & = 2n+1 & \forall{} n \geq 0,\\
    a_n & = -z^2 & \forall{} n \geq 1.
  \end{align*}
\end{subequations}
The Pad\'e approximant
can then be determined from the following recurrence formula~\cite{Baker1996}:
\begin{equation}
  \label{eq:closed:pade:euler}
  \frac{A_m}{B_l} = \frac{b_mA_{m-1}+a_mA_{m-2}}{b_lB_{l-1}+a_lB_{l-2}}
  \qquad\forall{}m\geq{}2~\text{and}~\forall{}l\geq{}2,
\end{equation}
with
\begin{equation*}
  \begin{array}{c@{\quad\text{and}\quad}c}
    \left\{
      \begin{array}{r@{~}c@{~}l}
        A_0 & = & b_0,\\
        B_0 & = & 1,
      \end{array}
    \right.
    &
    \left\{
      \begin{array}{r@{~}c@{~}l}
        A_1 & = & b_1b_0 + a_1,\\
        B_1 & = & b_1.
      \end{array}
    \right.
  \end{array}
\end{equation*}
That is, we have for the $z\cot(z)$ function:
\begin{equation}
  \label{eq:closed:pade:euler:cot}
  \begin{array}{c@{\quad\text{and}\quad}c@{\quad\text{and}\quad}c}
    \left\{
      \begin{array}{r@{~}c@{~}l}
        A_0 & = & 1,\\
        B_0 & = & 1,
      \end{array}
    \right.
    &
    \left\{
      \begin{array}{r@{~}c@{~}l}
        A_1 & = & 3-z^2,\\
        B_1 & = & 3,
      \end{array}
    \right.
    &
    \left\{
      \begin{array}{r@{~}c@{~}l}
        A_m & = & (2m+1)A_{m-1} - z^2A_{m-2},\\
        B_l & = & (2l+1)B_{l-1} - z^2B_{l-2}.
      \end{array}
    \right.
  \end{array}
\end{equation}
Starting from this recurrence formula and choosing $l=m$,
we can devise a $N$-term decomposition of the form
\begin{equation}
  \label{eq:pade:cot}
  z\cot(z)
  \approx \frac{A_N}{B_N}
  = \tilde{C}_0 + \sum_{i=0}^{N}\frac{\tilde{A}_i}{z^2-\tilde{B}_i}.
\end{equation}
However, compared with the unbounded case
where the coefficients of the Pad\'e approximant of the $\DtN$ map
are known analytically and exploited to construct the \tc{pade}{u}
operator~\cite{Boubendir2012},
no closed form formulae were found for
the coefficients  $\tilde{C}_0$, $\tilde{A}_i$ and $\tilde{B}_i$
of~\eqref{eq:pade:cot}.
Nevertheless, those can be computed numerically by
\begin{enumerate}
\item performing a polynomial long division of $A_N/B_N$,
  that is $A_N = \tilde{C}_0B_N + R$,
\item computing the poles $\tilde{B}_i$ of $R/B_N$ and
\item determining the residues $\tilde{A}_i$ of $R/B_N$.
\end{enumerate}
The numerically demanding part of this approach is the calculation
of the poles of $R/B_N$, \ie the zeros of $B_N$,
which requires arbitrary precision arithmetic
as the coefficients of the monomials appearing in $B_N$ can be vary large.

In this work, this is achieved with the \Code{MPSolve} library%
\footnote{See \url{github.com/robol/MPSolve}.}~\cite{Bini2014}.
Within that framework, it takes less than 5 minutes\footnote{This
  computation was carried out with a dual-core laptop-class Intel i7-7500U CPU.}
to compute the $\tilde{C}_0$, $\tilde{A}_i$ and $\tilde{B}_i$ coefficients
in the very large case of $N=1024$ for instance.
Of course, these coefficients can be pre-computed and tabulated for
various values of $N$
and the actual transmission condition recovered with the change of variable
$z = \ell_{ji}\sqrt{k^2-s^2}$ (see paragraph below).
For illustration purposes, the Pad\'e coefficients are
presented in \Tab{tab:closed:pade} for $N\in[1, 4]$.
\begin{table}[ht]
  \centering
  \begin{tabular}[t]{cccc@{\quad}c}
    \toprule
    $N$ & $\tilde{C}_0$      & $\tilde{A}_i$      & $\tilde{B}_i$\\
    \midrule
    $1$ & $6.00\times{}10^0$ & $7.50\times{}10^1$ & $1.50\times{}10^1$\\
    \midrule
    $2$ & $1.50\times{}10^1$ & $2.05\times{}10^1$ & $9.94\times{}10^0$\\
        &                    & $1.13\times{}10^3$ & $9.51\times{}10^1$\\
    \midrule
    $3$ & $2.80\times{}10^1$ & $1.97\times{}10^1$ & $9.87\times{}10^0$\\
        &                    & $1.09\times{}10^2$ & $4.20\times{}10^1$\\
        &                    & $7.31\times{}10^3$ & $3.26\times{}10^2$\\
    \bottomrule
  \end{tabular}$\quad$
  \begin{tabular}[t]{cccc@{\quad}c}
    \toprule
    $N$ & $\tilde{C}_0$      & $\tilde{A}_i$      & $\tilde{B}_i$\\
    \midrule
    $4$ & $4.50\times{}10^1$ & $1.97\times{}10^1$ & $9.87\times{}10^0$\\
        &                    & $8.03\times{}10^1$ & $3.96\times{}10^1$\\
        &                    & $4.03\times{}10^2$ & $1.06\times{}10^2$\\
        &                    & $3.02\times{}10^4$ & $8.35\times{}10^2$\\
    \\
    \\
    \\[-2pt]
    \bottomrule
  \end{tabular}
  \caption{Pad\'e coefficients $\tilde{C}_0$, $\tilde{A}_i$ and $\tilde{B}_i$
    for $N\in[1,4]$.}
  \label{tab:closed:pade}
\end{table}

Capitalizing on the above development,
we can now devise a new approximation (\tc{pade}{c}) of $\Lij{dtn}{c}$
of the form
\begin{equation}
  \label{eq:closed:pade:symbol}
  \Lij{dtn}{c}(s) \approx \Lij{pade($N$)}{c}(s)
  = \frac{1}{\ell_{ji}}
    \left(
      \tilde{C}_0 +
      \sum_{n=1}^N\frac{\tilde{A}_n}{\Big(k\ell_{ji}\Big)^2
                                     \Big(1-\frac{s^2}{k^2}\Big)-\tilde{B}_n}
    \right),
\end{equation}
by exploiting the change of variable $z = \ell_{ji}\sqrt{k^2-s^2}$.
The operator associated with this symbol then reads:
\begin{equation}
  \label{eq:closed:pade:S}
  \Sij{pade($N$)}{c} =
  \frac{1}{\ell_{ji}}
    \left\{
      \tilde{C}_0 +
      \sum_{n=1}^N\tilde{A}_n\Bigg[
        \Bigg(k\ell_{ji}\Bigg)^2\Bigg(1+\frac{\Div_{\Sigma_{ij}}\Grad_{\Sigma_{ij}}}{k^2}\Bigg)-\tilde{B}_n
      \Bigg]^{-1}
    \right\}.
\end{equation}
As in the truncated Mittag-Leffler expansion case,
the above operator has been developed for a one-dimensional interface,
but can be used as well with a two-dimensional one,
as shown in section~\ref{sec:cryo} for instance.

\begin{figure}[t]
  \centering
  \includegraphics{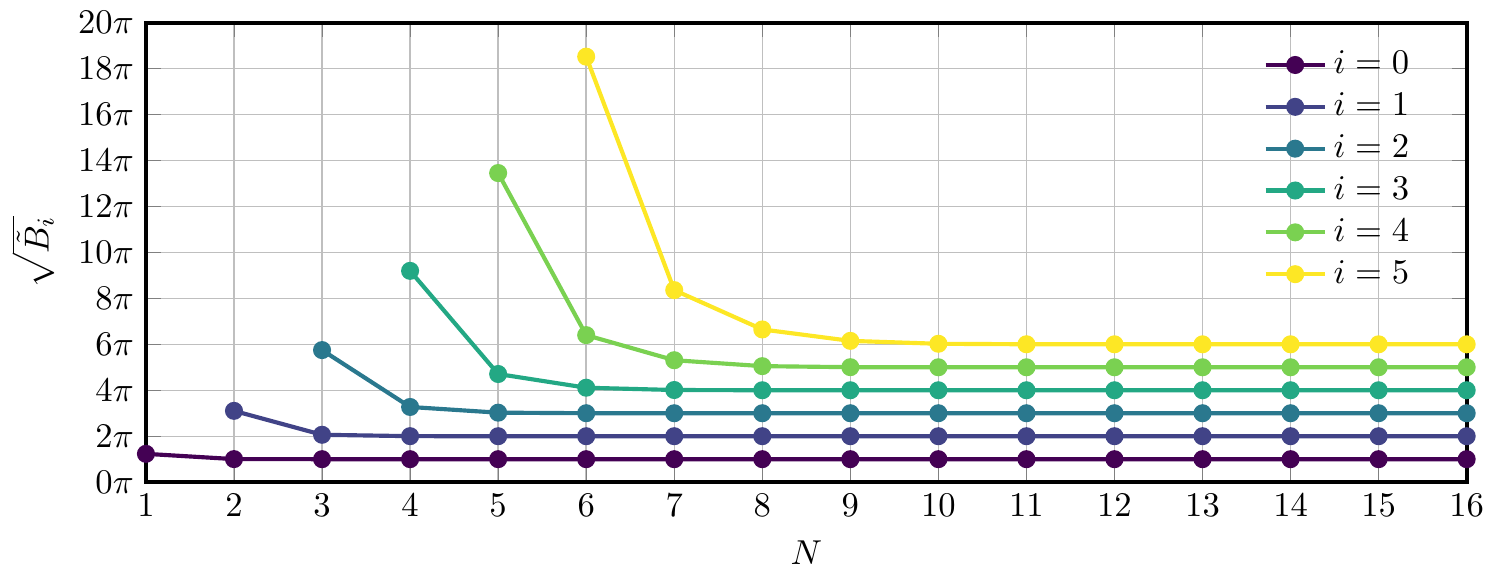}
  \caption{Convergence of the square-root of the $i$\textsuperscript{th} pole
    of the Pad\'e approximant with $N$.}
  \label{fig:pade:ml:poles}
\end{figure}

Before concluding this subsection,
let us mention that the truncated Mittag-Leffler expansion of $z\cot(z)$
is related to its Pad\'e approximant in the following way:
the square-root of the $i$\textsuperscript{th} pole
appearing in~\eqref{eq:pade:cot} is converging towards $(i+1)\pi$,
\ie the $(i+1)$\textsuperscript{th} pole of $z\cot(z)$
appearing in~\eqref{eq:closed:ml:cot},
for a fixed value of $i$ and as $N$ increases.
Formally, we have that
\begin{equation}
  \label{eq:pade:ml:poles}
  \lim_{N\to\infty}\sqrt{\tilde{B}_i} = (i+1)\pi\quad\forall i\geq0,
\end{equation}
as show in \Fig{fig:pade:ml:poles} for the first six poles.

\begin{figure}[b]
  \centering
  \includegraphics{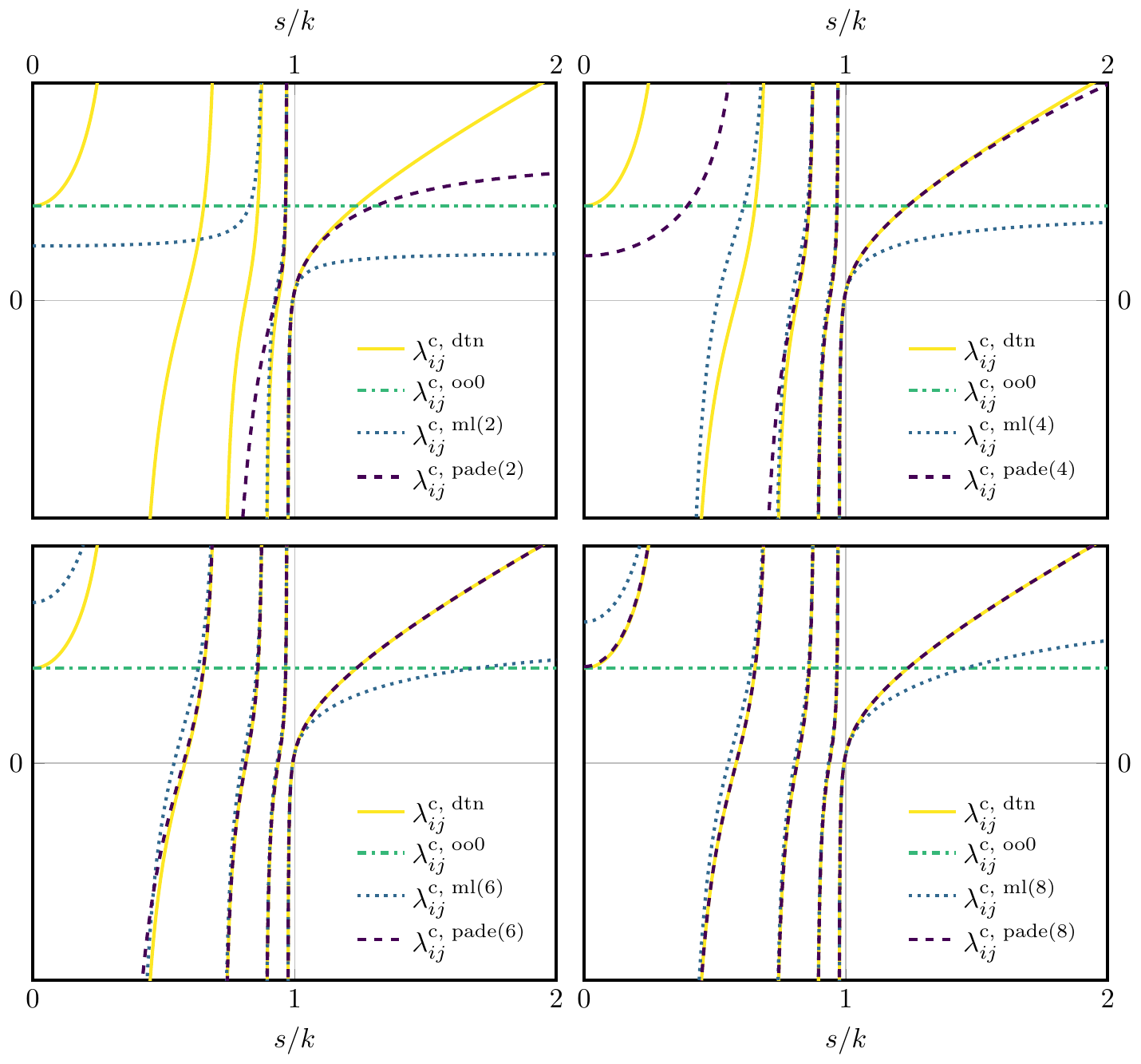}
  \caption{Approximation quality of the different transmission operators
    in a one-dimensional setting with $\ell/\lambda_w=4.3$
    and two subdomains of equal size
    (arbitrary values for the y-axis,
    only the zero is shown to highlight the sign changes).}
  \label{fig:approx:1d}
\end{figure}
\afterpage{\clearpage}

\subsection{Quality of the transmission operators in the one-dimensional case}
As the $\DtN$ and the three transmission operators discussed above
are purely real,
they reduce to a simple real function in the one-dimensional case,
allowing a simple graphical comparison,
as suggested in \Fig{fig:approx:1d}.
For clarity reasons, this figure is restricted
to a relatively low frequency problem with a ratio
$\ell/\lambda_w=4.3$, where $\lambda_w = 2\pi/k$ is the wavelength,
and two subdomains of equal size.
Nonetheless, the discussion below remains general.

To start with, it is clear that \tc{oo0}{c} is a quite poor approximation
of the $\DtN$, as it obviously cannot capture its oscillations and its poles.
It is also apparent that changing the expansion point of the Taylor series
(here $s=0$, as a recall) will not improve the situation.

On the other hand, both \tc{ml}{c} and \tc{pade}{c} are able to capture
the oscillation and poles of the $\DtN$,
at least for a sufficiently high value of $N$.
Those approximations exhibit however major differences,
which are summarized in what follows.
\begin{itemize}
\item For the evanescent modes (\ie when $s/k>1$), the convergence of
  $\Lij{ml($N$)}{c}$ towards $\Lij{dtn}{c}$ as $N$ increases is quite slow,
  at least when compared with $\Lij{pade($N$)}{c}$.
  Let us mention that in the unbounded case, $\Lij{pade($N$)}{u}$
  is known to perform well when $s/k>1$.
  This aspect will be further discussed in section~\ref{sec:open}.

\item For the non-evanescent modes (\ie when $s/k<1$),
  each pole introduced by $\Lij{ml($N$)}{c}$ coincides,
  by construction, with a pole of $\Lij{dtn}{c}$.
  Thus for low values of $N$
  where all the poles of $\Lij{dtn}{c}$ are not yet present,
  $\Lij{ml($N$)}{c}$ leads to better approximations of $\Lij{dtn}{c}$
  than $\Lij{pade($N$)}{c}$,
  at least when restricting $s/k$ to the range spanned
  by the first and last poles of $\Lij{ml($N$)}{c}$.
  Conversely, for higher values of $N$,
  the approximation given by $\Lij{pade($N$)}{c}$
  is clearly the best one.
\end{itemize}
To conclude this subsection,
let us also note that the approximation is better as $s/k$ tends to $1$
for both $\Lij{ml($N$)}{c}$ and $\Lij{pade($N$)}{c}$.

\subsection{Regularization with a constant imaginary part}
\label{sec:tc:reg}
As already mentioned earlier, the zeroth-order transmission condition optimized
for the cavity problem can lead to an ill-posed problem
when one of the $d_{ij}(s)$ terms in the denominator
of the convergence radius~\eqref{eq:fourier:rho}
equals (or is sufficiently close to) zero.
This problem is however not peculiar to the \tc{oo0}{c} condition
and affects the \tc{ml}{c} and \tc{pade}{c} ones as well,
requiring therefore a regularization procedure.

The most straightforward and simple strategy to regularize
the \tc{oo0}{c}, \tc{ml}{c} and \tc{pade}{c} conditions
is to exploit the fact that those operators are purely real-valued.
Consequently, by adding a purely imaginary part,
the $d_{ij}(s)$ terms can be pushed away from zero
and the convergence radius $\rho(s)$ can be guaranteed to be well defined.
Formally, the regularized \tc{oo0}{c}, \tc{ml}{c} and \tc{pade}{c} operators
read\footnote{We use the superscript ``r($\chi$)'' to denote the regularization
  with a constant imaginary part proportional to $\chi$.}
\begin{equation}
  \label{eq:closed:reg}
  \Sij{oo0/ml($N$)/pade($N$), r($\chi$)}{c} =
  \Sij{oo0/ml($N$)/pade($N$)}{c} + \imag\chi k,
\end{equation}
where $\chi\in\mathbb{R}$ is the regularization parameter.
Numerical experiments showing the impact of this regularization approach
are further discussed in sections~\ref{sec:validation}
and~\ref{sec:sensitivity}.

\subsection{Regularization by mixing operators
  optimized for cavity and unbounded problems}
\label{sec:tc:mixed}
The above regularization is in some sense suboptimal
as it acts on all $s\in\mathbb{S}$,
while regularization is required only in the $s^2<k^2$ case.
A more selective approach can be achieved by exploiting the
\tc{pade}{u} operator.
Indeed, assuming a sufficiently high value of $N$
and an appropriate rotation of the branch cut~\cite{Milinazzo1997}
of \tc{pade}{u},
the latter exhibits the following properties~\cite{Milinazzo1997}:
\begin{enumerate}
\item it is approximately imaginary when $s^2<k^2$,
\item it is approximately real when $s^2>k^2$ and
\item it is a good approximation of the $\DtN$ map when $s^2>k^2$
  (see discussion in section~\ref{sec:open}).
\end{enumerate}
Therefore, a regularized operator can be constructed by combining either
the \tc{oo0}{c}, the \tc{ml}{c} or the \tc{pade}{c} operator
with the \tc{pade}{u} one in a convex way.
These new operators will be further referred to as
\emph{mixed} operator and write\footnote{We use the superscript
  ``m($\epsilon, M$)'' to denote mixed operators involving
  a $M$-term \tc{pade}{u} operator with a weight of $(1-\epsilon)$.}
\begin{equation}
  \label{eq:closed:mixed}
  \Sij{oo0/ml($N$)/pade($N$)}{m($\epsilon, M$)} =
  \epsilon\Sij{oo0/ml($N$)/pade($N$)}{c} + (1-\epsilon)\Sij{pade($M$)}{u},
\end{equation}
where $\epsilon\in]0,1[$ denotes the regularization parameter
of the mixed formulation.

\subsection{Estimate for the minimum number of auxiliary unknowns}
\label{sec:tc:nmin}
Given the \tc{ml}{c} and \tc{pade}{c} transmission conditions,
one question naturally arises: how many auxiliary terms should be selected?
In order to answer this question, one could opt for the following criterion:
the number of auxiliary terms $N$ should at least be equal to the number
of poles of $\Lij{dtn}{c}$,
as it seems legitimate to assume that the behavior of
$\Lij{dtn}{c}$ is driven by its poles in the range $s^2<k^2$.

Given the Mittag-Leffler expansion of $\Lij{dtn}{c}$,
it is clear that the poles must satisfy
\begin{equation}
  \label{eq:closed:nmin:poles}
  s^2 = k^2-\left(\frac{n\pi}{\ell_{ij}}\right)^2\quad\forall{}n>0,
\end{equation}
which implies that
\begin{equation}
  \label{eq:closed:nmin:pos}
  k^2-\left(\frac{n\pi}{\ell_{ij}}\right)^2 > 0
\end{equation}
and thus
\begin{equation}
  \label{eq:closed:nmin:range}
  0 < n < 2\frac{\ell_{ij}}{\lambda_w},
\end{equation}
since $k$, $n$ and $\ell_{ij}$ are positive by construction.

Therefore, according to the pole criterion stated above,
the minimum number of terms $N_\text{min}^\text{poles}$ for localizing
$\Lij{dtn}{c}$ is
\begin{equation}
  \label{eq:closed:nmin}
  N_\text{min}^\text{pole} = \left\lceil2\frac{\ell_{ij}}{\lambda_w}\right\rceil
\end{equation}
and depends on the size of the subdomains and on the wavelength.
As discussed further in section~\ref{sec:validation:naux},
this criterion seems however pessimistic,
as lower values of $N_\text{min}^\text{pole}$
provide already acceptable results.

\section{Generalizations}
\label{sec:generalization}
The transmission operators detailed in the pervious section are restricted
to the very particular case of a cavity with \emph{Dirichlet} boundary conditions
divided into \emph{two} subdomains.
In order to generalize this setting,
we discuss in section the cases of many subdomains
in a one-dimensional partitioning and of Neumann boundary conditions.

\subsection{One-dimensional partitioning with more than two subdomains}
\begin{figure}[t]
  \centering
  \begin{subfigure}[b]{0.49\textwidth}
    \centering
    \includegraphics{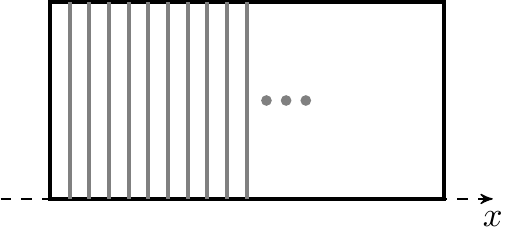}
    \caption{Partitioning into $D$ subdomains.}
    \label{fig:domain:many}
  \end{subfigure}
  \begin{subfigure}[b]{0.49\textwidth}
    \centering
    \includegraphics{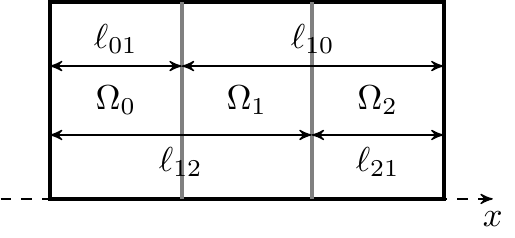}
    \caption{Definition of $\ell_{ji}$  (here for $D=3$ subdomains).}
    \label{fig:domain:lij}
  \end{subfigure}
  \caption{One-dimensional partitioning of the rectangular cavity.}
\end{figure}

Let us start with the one-dimensional partitioning
of the computational domain into $D$~subdomains,
as shown in \Fig{fig:domain:many}.
In such a one-dimensional domain decomposition,
the physical meaning of the $\ell_{ji}$ coefficient appearing in
$\Lij{dtn}{c}(s)$ must be clarified.
In the two subdomains case,
\emph{the $\ell_{ji}$ coefficients represents the distance between
$\Sigma_{ij}$ and the reflecting wall located in the $\vec{n}_{ij}$ direction}.
This interpretation can be directly applied to the $D$~subdomains case
to define the $\ell_{ji}$ coefficients,
as shown in \Fig{fig:domain:lij} for the $D=3$ case.

\subsection{Neumann boundary conditions}
So far, we considered only the situation
where the reflecting wall associated with $\Trsm_{ij}$
is implemented with an homogeneous Dirichlet boundary condition,
\ie a soft-wall condition.
In the case of an homogeneous Neumann boundary condition,
\ie a hard-wall condition,
the following $\DtN$ map is obtained:
\begin{subequations}
  \label{eq:close:lambda:dtn:neumann}
  \begin{align}[left = {\Lij{dtn, neumann}{c}(s) = \empheqlbrace}]
    & -\sqrt{k^2-s^2}
      \tan\mathopen{}\left[\ell_{ji}\sqrt{k^2-s^2}\mathclose{}\right]
    & \text{if}~s^2<k^2,\\
    & 0
    & \text{if}~s^2=k^2,\\
    & \sqrt{s^2-k^2}
      \tanh\mathopen{}\left[\ell_{ji}\sqrt{s^2-k^2}\mathclose{}\right]
    & \text{if}~s^2>k^2,
  \end{align}
\end{subequations}
by following the same strategy as in section~\ref{sec:S:dirichlet}.
It is worth stressing the similarities
between~\eqref{eq:close:lambda:dtn} and~\eqref{eq:close:lambda:dtn:neumann}.
This $\DtN$ map can then be localized using the previously presented approaches
and a \tc{oo0}{c}, \tc{ml}{c} or \tc{pade}{c} transmission condition
can be devised.
In this regard, let us note that the Mittag-Leffler expansion of $\tan(z)$
reads~\cite{Titchmarsh1976}
\begin{equation}
  \label{eq:closed:ml:tan}
  z\tan(az) =
  -\frac{2}{a}
  \sum_{n=0}^\infty\frac{z^2}{z^2-\left[\frac{(n+\frac{1}{2})\pi}{a}\right]^2}
\end{equation}
and its continued fraction expansion has the following form~\cite{Oldham2009}:
\begin{equation}
  \label{eq:closed:pade:tan:cfrac}
  z\tan(z)
  = \cfrac{z^2}{1-\cfrac{z^2}{3-\cfrac{z^2}{5-\dots}}}
  = b_0 + \cfrac{a_1}{b_1+\cfrac{a_2}{b_2+\cfrac{a_3}{b_3+\dots}}}.
\end{equation}
These results are given here for the sake of completeness
and will not be further discussed in this work.

\section{Operators optimized for unbounded problems without obstacles
  in a cavity context}
\label{sec:open}
In this section we derive some estimates on the performance
of operators optimized for unbounded problems without obstacles
when used in a cavity problem.
In particular,
\begin{myenum*}
\item we first compare the $\DtN$ operator
  related to an unbounded problem without obstacles $\Sij{dtn}{u}$
  with its cavity counter part $\Sij{dtn}{c}$, then
\item we discuss the convergence radius of the OS scheme
  when using $\Sij{dtn}{u}$ as a transmission operator
  for the rectangular cavity problem~\eqref{eq:helmholtz}, as well as
\item the particular case of the optimized order 0 operator
  (\tc{oo0}{u})~\cite{Despres1990}.
\end{myenum*}

\subsection{Dirichlet-to-Neumann operators}
Let us consider the following Helmholtz problem without obstacles:
\begin{subequations}
  \label{eq:helmholtz:open}
  \begin{align}[left = \empheqlbrace]
    \Div{\Grad{p}} + k^2p & = g \quad\text{on}~\mathbb{R}^2,\\
    \lim_{r\to\infty}\sqrt{r}\left(\deriv{p}{r}-\imag{}kp\right) & = 0,
    \label{eq:helmholtz:open:rad}
  \end{align}
\end{subequations}
where $r^2 = x^2+y^2$.
In this case, it can be shown that the optimal transmission operator
$\Sij{dtn}{u}$
for solving this problem with an OS scheme is~\cite{Boubendir2012}:
\begin{equation}
  \label{eq:open:S:dtn}
  \Sij{dtn}{u} = \Op\mathopen{}\Big(\Lij{dtn}{u}\mathclose{}\Big) =
  -\imag{}k\sqrt{1+\frac{\Div_{\Sigma_{ij}}\Grad_{\Sigma_{ij}}}{k^2}},
\end{equation}
where
\begin{equation}
  \label{eq:open:lambda:dtn}
  \Lij{dtn}{u} = -\imag{}k\sqrt{1-\frac{s^2}{k^2}}.
\end{equation}

By comparing~\eqref{eq:open:lambda:dtn} and~\eqref{eq:close:lambda:dtn},
it is easy to realize that
\begin{subequations}
  \label{eq:open:close}
  \begin{align}[left = {\Lij{dtn}{c}(s) - \Lij{dtn}{u}(s)=\empheqlbrace}]
    & \sqrt{k^2-s^2}
      \cot\mathopen{}\left[\ell_{ji}\sqrt{k^2-s^2}\mathclose{}\right]
      +\imag{}\sqrt{k^2-s^2}
    & \text{if}~s^2<k^2, \label{eq:open:close:s<k}\\
    & 1/\ell_{ji}
    & \text{if}~s^2=k^2,\\
    & \sqrt{s^2-k^2}
      \coth\mathopen{}\left[\ell_{ji}\sqrt{s^2-k^2}\mathclose{}\right]
      -\sqrt{s^2-k^2}
    & \text{if}~s^2>k^2.
  \end{align}
\end{subequations}
Interestingly, by exploiting the definition of the hyperbolic
cotangent~\cite{Oldham2009}, the case $s^2>k^2$ can be further simplified into
\begin{align}
  \Lij{dtn}{c}(s)-\Lij{dtn}{u}(s)
  & = \sqrt{s^2-k^2}
    \coth\mathopen{}\left[\ell_{ji}\sqrt{s^2-k^2}\mathclose{}\right]
    -\sqrt{s^2-k^2}\nonumber\\
  & = \sqrt{s^2-k^2}
      \left(\frac{\exp\mathopen{}\left(2\ell_{ji}\sqrt{s^2-k^2}\mathclose{}\right)+1}
                 {\exp\mathopen{}\left(2\ell_{ji}\sqrt{s^2-k^2}\mathclose{}\right)-1}
            - 1\right)\nonumber\\
  & = \frac{2}{\exp\mathopen{}\left(2\ell_{ji}\sqrt{s^2-k^2}\mathclose{}\right)-1}
      \sqrt{s^2-k^2}
  & \text{if}~s^2>k^2, \label{eq:open:close:s>k}
\end{align}
which yields:
\begin{equation}
  \label{eq:open:close:lim}
  \lim_{s\to\infty} \Lij{dtn}{c}(s) - \Lij{dtn}{u}(s) = 0.
\end{equation}
In other words, for the case $s^2>k^2$,
the symbol $\Lij{dtn}{u}(s)$ is converging towards $\Lij{dtn}{c}(s)$
as $s$ grows.
Furthermore,
as the difference between both symbols is decreasing exponentially,
\emph{$\Lij{dtn}{u}(s)$ is an excellent approximation of $\Lij{dtn}{c}(s)$
  when $s^2>k^2$.}

\begin{figure}[t]
  \centering
  \includegraphics{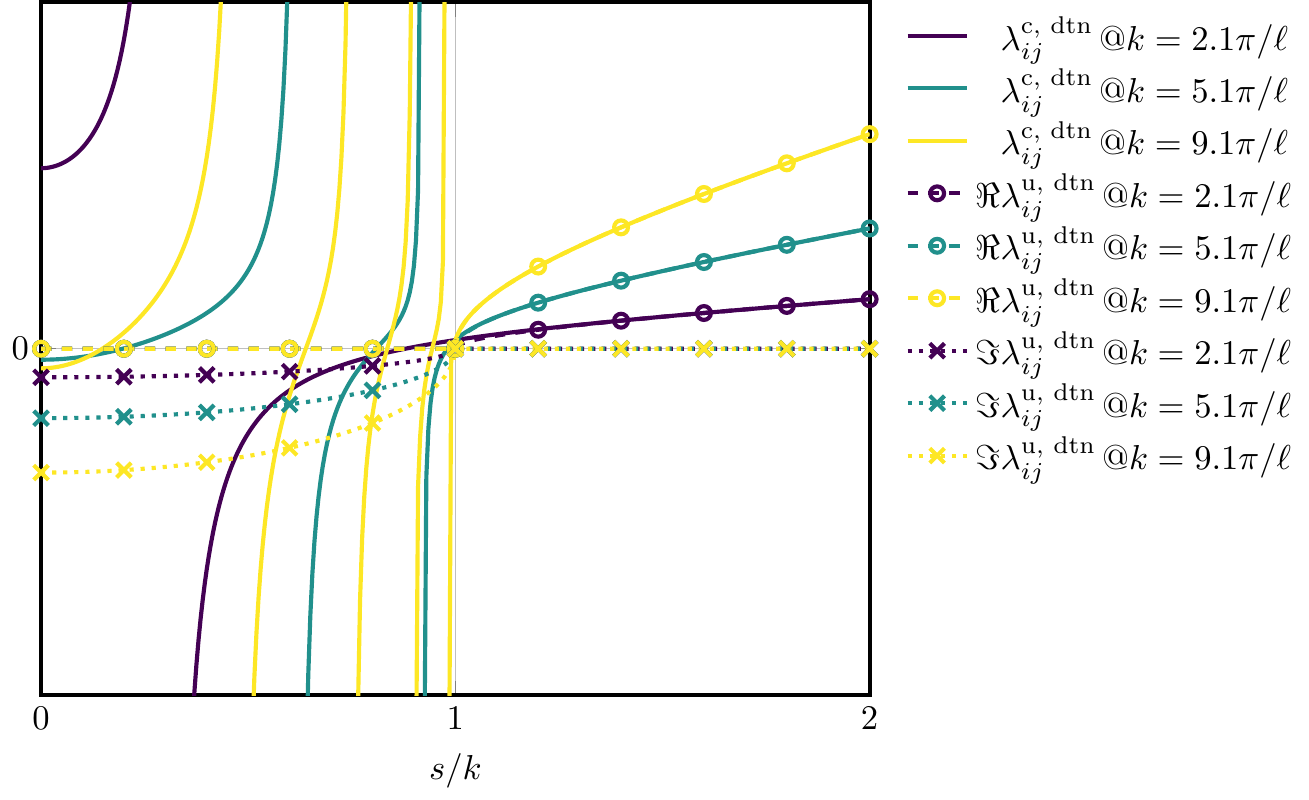}
  \caption{Graphs of
    $\Lij{dtn}{u}(s/k)$ and $\Lij{dtn}{c}(s/k)$ for different values of $k$
    (arbitrary values for the y-axis,
    only the zero is shown to highlight the sign changes).}
  \label{fig:lambda:k}
\end{figure}

Regarding the case $s^2<k^2$, as the codomains of
$\Lij{dtn}{u}(s)$ (which is purely imaginary)
and $\Lij{dtn}{c}(s)$ (which is purely real)
do not match, the expression in~\eqref{eq:open:close:s<k}
cannot be further simplified.
For illustration purposes, the graphs of
$\Lij{dtn}{u}$ and $\Lij{dtn}{c}$
are depicted in \Fig{fig:lambda:k} for different values of $k$
($\Re$ and $\Im$ are respectively denoting
the real and imaginary part functions).

\subsection{Best convergence radius}
We already know from the previous section that
$\Lij{dtn}{u}(s)$ is a good approximation of $\Lij{dtn}{c}(s)$ when $s^2>k^2$,
that is for \emph{evanescent modes}.
For this reason, local approximations of $\Lij{dtn}{u}$
are legitimate, yet suboptimal, candidates for approximating $\Lij{dtn}{c}$.

In terms of convergence radius, as defined in~\eqref{eq:fourier:rho},
it is easy to show that
\begin{subequations}
  \begin{align}[left = {\rho^2(s)=\empheqlbrace}]
    & \frac{\imag-\cot\mathopen{}\left[\ell_{ji}\sqrt{k^2-s^2}\mathclose{}\right]}
           {\imag+\cot\mathopen{}\left[\ell_{ij}\sqrt{k^2-s^2}\mathclose{}\right]}
      \frac{\imag-\cot\mathopen{}\left[\ell_{ij}\sqrt{k^2-s^2}\mathclose{}\right]}
           {\imag+\cot\mathopen{}\left[\ell_{ji}\sqrt{k^2-s^2}\mathclose{}\right]}
    & \text{if}~s^2<k^2,\\
    & 1
    & \text{if}~s^2=k^2,\\
    & \frac{1-\coth\mathopen{}\left[\ell_{ji}\sqrt{s^2-k^2}\mathclose{}\right]}
           {1+\coth\mathopen{}\left[\ell_{ij}\sqrt{s^2-k^2}\mathclose{}\right]}
      \frac{1-\coth\mathopen{}\left[\ell_{ij}\sqrt{s^2-k^2}\mathclose{}\right]}
           {1+\coth\mathopen{}\left[\ell_{ji}\sqrt{s^2-k^2}\mathclose{}\right]}
    & \text{if}~s^2>k^2,
  \end{align}
\end{subequations}
and
\begin{subequations}
  \label{eq:open:rho}
  \begin{align}[left = {\abs{\rho(s)}=\empheqlbrace}]
    & 1
    & \text{if}~s^2\leq{}k^2,\\
    & \exp\mathopen{}\left(-\ell\sqrt{s^2-k^2}\mathclose{}\right)
    & \text{if}~s^2>k^2,
  \end{align}
\end{subequations}
when  $\Lij{dtn}{c}$ is approximated with $\Lij{dtn}{u}$.
For this reason, the transmission operators that are good approximations of
$\Lij{dtn}{u}$,
such that the optimized order 2~(\tc{oo2}{u})~\cite{Gander2002}
or the $N$-term Pad\'e-localized~(\tc{pade}{u})~\cite{Boubendir2012} operators,
should exhibit a convergence radius close to~\eqref{eq:open:rho}.
In other words, those local operators should exhibit a slow convergence
for the \emph{propagating modes} ($s^2<k^2$)
and a fast convergence for the \emph{evanescent} ones ($s^2>k^2$).

\subsection{Particular case of the optimized order 0 operator}
Before concluding this subsection, it is worth mentioning that in the case
of the \tc{oo0}{u} operator,
we have that $\Lij{oo0}{u}=-\imag{}k$ and therefore
\begin{equation}
  \label{eq:open:rho:oo0}
  \abs{\rho(s)} = 1.
\end{equation}

\section{Numerical validation and comparison between the different operators}
\label{sec:validation}
In this section we analyze the performance
of the different transmission conditions developed in section~\ref{sec:tc}
and compare them with the operators of section~\ref{sec:open}.
To this end, we consider the rectangular cavity
shown in \Fig{fig:domain:omega}
and solve the time-harmonic Helmholtz equation over $\Omega$
with the following boundary conditions imposed on $\Gamma$:
\begin{subequations}
  \label{eq:num:bc}
  \begin{align}[left = \empheqlbrace]
    p & = 0 & \text{on}~\Gamma_w\cup\Gamma_r,\\
    p & = \displaystyle
        \sum_{m=1}^K\sin\left(m\frac{\pi}{h}y\right) & \text{on}~\Gamma_l,
  \end{align}
\end{subequations}
where $K$ designates the number of modes used to excite the cavity.
This Helmholtz problem is then decomposed into $D$ subdomains of equal size
and solved with an optimized Schwarz scheme
combined with a GMRES algorithm \emph{without restart}.
Let us mention as well that the software implementation relies on
the \Code{GmshDDM} and \Code{GmshFEM}~\cite{Royer2021} frameworks\footnote{See
  \url{git.rwth-aachen.de/marsic/closeddm},
  \url{gitlab.onelab.info/gmsh/ddm} and
  \url{gitlab.onelab.info/gmsh/fem}.}
and exploits a finite element (FE) discretization of the subproblems.
In the case of the
$\Sij{pade($N$)}{c}$ and $\Sij{ml($N$)}{c}$ operators,
the FE variational formulations involve auxiliary unknowns
for the treatment of the inverse operation,
as proposed in~\cite{Boubendir2012}.
However, let us note that since the new transmission operators are not symmetric,
\ie $\Sij{oo0/ml($N$)/pade($N$)}{c}\neq\Sji{oo0/ml($N$)/pade($N$)}{c}$,
the amount of auxiliary fields must be doubled.

As the novel transmission conditions involve operators
oscillating rapidly in a wide range,
the finite precision arithmetic aspects of the linear solver
must be treated with care.
In this regard, the orthogonalization step of the GMRES solver is critical
and the \emph{modified} version of the Gram-Schmidt algorithm~\cite{Saad2003, Saad1986}
is required for convergence.
Concerning the software implementation,
the linear solvers of the \Code{PETSc}~\cite{Balay1997} (GMRES)
and \Code{MUMPS}~\cite{Amestoy2001} (LU factorization) libraries are used.

Unless stated otherwise, the cavity has an aspect ratio of $\ell/h = 2$
and a length-to-wavelength ratio of
$\ell/\lambda_w = \frac{157.085}{2\pi} \approx 25.001$.
This configuration allows $25$ non-evanescent modes and
is excited with the $K=50$ first modes.
The geometry is discretized with a mesh consisting
of $8$ triangular elements per wavelength
and the subproblems are discretized with an FE method of order $4$.
The stopping criterion of the GMRES solver is set
to a relative tolerance decrease of $\norm{\vec{r}_i}/\norm{\vec{r}_0}=10^{-6}$,
where $\vec{r}_i$ is the residual vector at iteration~$i$ and
$\vec{r}_0$ is the residual vector of the first guess,
which was chosen equal to zero.

\subsection{Two subdomains case}
\begin{figure}[ht]
  \centering
  \includegraphics{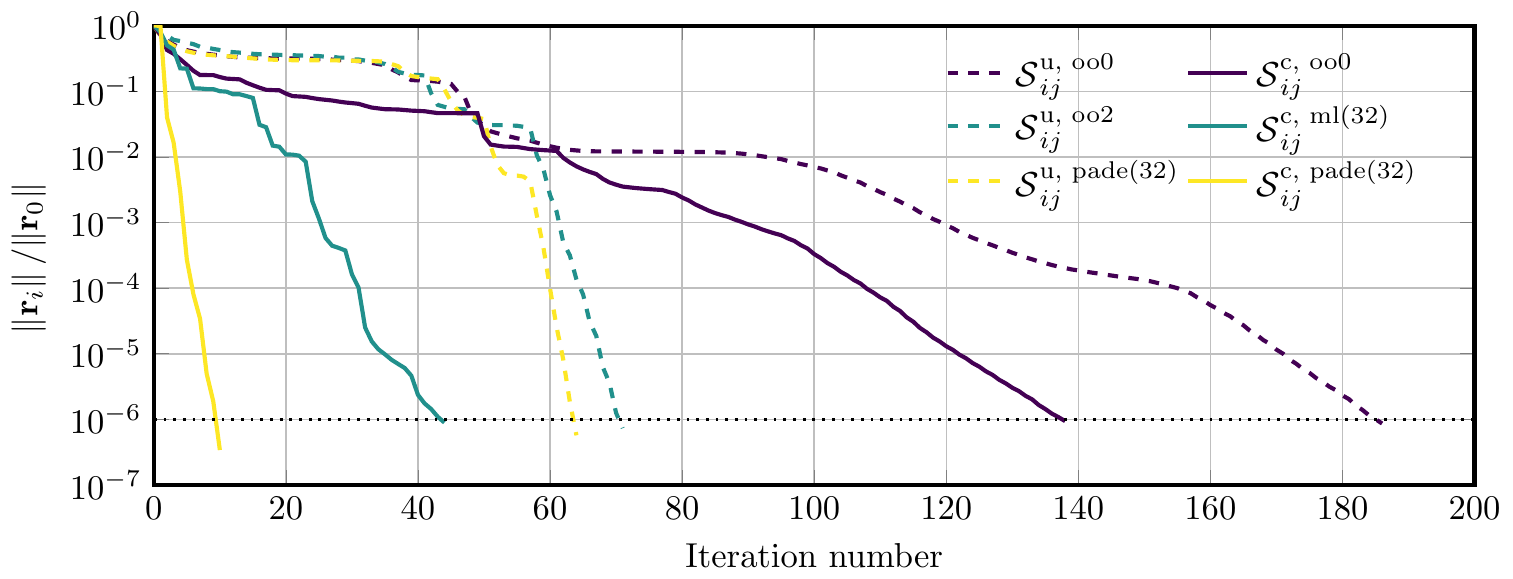}
  \caption{Convergence history of the GMRES solver
    -- rectangular cavity with $D=2$ subdomains.}
  \label{fig:val:d2}
\end{figure}

Let us start by studying the performance of the transmission conditions
developed in section~\ref{sec:tc} when applied to
a rectangular cavity divided into two subdomains.
The convergence history of the GMRES solver is shown in \Fig{fig:val:d2}.
From these data, it is clear that all transmission conditions converge.
In addition,
it appears clearly that $\Sij{ml($32$)}{c}$ and $\Sij{pade($32$)}{c}$ outperform
the transmission conditions devised for unbounded problems.
Let us mention that, while $\Sij{oo0}{c}$ converges with less iterations
than $\Sij{oo0}{c}$ in this example, the converse may happen
as shown in section~\ref{sec:validation:reg}.
Furthermore, the best operator in this numerical experiment is
$\Sij{pade($32$)}{c}$ that converges without any noticeable plateau.

As the closed form solution $p_\text{CF}$ of this canonical problem is known,
\ie
\begin{equation}
  \label{eq:pAna:rect}
  p_\text{CF} = \sum_{m=1}^K\frac{\sin\left[k_x(m)(\ell-x)\right]}
                               {\sin\left[k_x(m)\ell\right]}
                \sin\left[k_y(m)\right]
\end{equation}
with
\begin{equation}
  \label{eq:kxy}
  k_y(m) = m\frac{\pi}{h}\quad\text{and}\quad k_x(m) = \sqrt{k^2-k_y^2(m)},
\end{equation}
let us determine the accuracy of the above simulations
by computing the relative $L_2$ error $\Error$
between the FE solution $p_\text{FE}$ and $p_\text{CF}$:
\begin{equation}
  \label{eq:error:l2}
  \Error = \sqrt{
    \frac{\displaystyle\int_\Omega\abs{p_\text{CF} - p_\text{FE}}^2\dd\Omega}
         {\displaystyle\int_\Omega\abs{p_\text{CF}}^2\dd\Omega}
  },
\end{equation}
where the above integrals are evaluated using a quadrature rule with twice
the amount of integration points than the number used for the FE computations.
The $L_2$ errors associated with the different OS schemes of \Fig{fig:val:d2}
are gathered in \Tab{tab:val:d2:error}
together with the error associated with the \Code{MUMPS} direct solver.
These data show clearly that all transmission conditions and the direct solver
lead to errors of the same order of magnitude.
\begin{table}[ht]
  \centering
  \begin{tabular}{ccccccc}
    \toprule
    $\Sij{oo0}{u}$ & $\Sij{oo2}{u}$ & $\Sij{pade($32$)}{u}$
    & $\Sij{oo0}{c}$ & $\Sij{ml($32$)}{c}$ & $\Sij{pade($32$)}{c}$
    & \Code{MUMPS} \\
    \midrule
      $1.76\times{}10^{-4}$ & $1.76\times{}10^{-4}$ & $1.78\times{}10^{-4}$
    & $1.77\times{}10^{-4}$ & $1.76\times{}10^{-4}$ & $1.77\times{}10^{-4}$
    & $1.78\times{}10^{-4}$ \\
    \bottomrule
  \end{tabular}
  \caption{Relative $L_2$ errors $\Error$
    -- rectangular cavity with $D=2$ subdomains.}
  \label{tab:val:d2:error}
\end{table}

Let us now focus on the wall-clock time required
to complete the above computations\footnote{Those
  calculations were carried out
  with an eight-core desktop-class Intel Xeon E5-2630 CPU
  and parallelized with two processes with 4 threads each.},
as reported in \Tab{tab:val:d2:time}.
In order to analyze these values,
let us stress that they heavily depend on the actual software implementation
of the FE and OS tools.
Nonetheless, some general remarks can be drawn.
\begin{enumerate}
\item The \tc{oo0}{u}, \tc{oo2}{u} and \tc{oo0}{c} operators are
  the computationally cheapest to apply,
  as they do not involve auxiliary unknowns.
\item The \tc{ml}{c} and \tc{pade}{c} operators
  are computationally more expensive than \tc{pade}{u},
  since they require two sets of auxiliary unknowns as they are not symmetric,
  \ie $\Sij{ml($N$)/pade($N$)}{c} \neq \Sji{ml($N$)/pade($N$)}{c}$.
\item The regularization procedure involving \tc{pade}{u}
  further increases the computational cost,
  as additional auxiliary unknowns are introduced.
\end{enumerate}
In the current software implementation,
the subproblems are solved with the direct solver \Code{MUMPS}
and the resulting LU factorization is reused in the subsequent
OS iterations.
Therefore, the first iteration is more time consuming than the other ones
and the data in \Tab{tab:val:d2:time} are thus split in different subquantities.
By defining  $T_i$ as the share of the total wall clock time $T_\text{tot}$
dedicated to the $i$\textsuperscript{th} iteration,
\ie $T_\text{tot}=\sum_{i=1}^IT_i$
with $I$ the total number of GMRES iterations required for convergence,
we report in \Tab{tab:val:d2:time} the following values:
$T_\text{tot}$, $T_1$, $\overline{T_{2,I}}$, where $\overline{T_{i,j}}$
is the mean value of the sequence $[T_i, \dots, T_I]$, and $I$.
\begin{table}[ht]
  \centering
  \begin{tabular}{cllllllc}
    \toprule
    Quantity of interest
    & $\Sij{oo0}{u}$ & $\Sij{oo2}{u}$      & $\Sij{pade($32$)}{u}$
    & $\Sij{oo0}{c}$ & $\Sij{ml($32$)}{c}$ & $\Sij{pade($32$)}{c}$
    & Unit \\
    \midrule
    $T_\text{tot}$
    & 113  & 48   & 48    & 90   & 41    & 17    & s \\
    $T_1$
    & 9.52 & 9.37 & 10.14 & 9.44 & 10.90 & 10.99 & s \\
    $\overline{T_{2,I}}$
    & 0.56 & 0.56 & 0.60  & 0.58 & 0.69  & 0.69  & s \\
    $I$
    & 186  & 71   & 64    & 138  & 44    & 10    & - \\
    \bottomrule
  \end{tabular}
  \caption{Wall-clock times -- rectangular cavity with $D=2$ subdomains
    ($T_i$ is the share of the total wall-clock time, in seconds, taken
    by the $i$\textsuperscript{th} iteration,
    $I$ is the total number of GMRES iterations and
    $\overline{T_{i,j}}$ is the mean value of the sequence $[T_i, \dots, T_j]$).}
  \label{tab:val:d2:time}
\end{table}

From the data gathered in \Tab{tab:val:d2:time},
it is clear that $\Sij{pade($32$)}{c}$ leads to both
the minimal amount of iterations and the fastest computation with respect to
the wall clock time.
Nonetheless, it is evident that the cost of $T_1$ and $\overline{T_{2,I}}$
is higher for $\Sij{pade($32$)}{c}$ and $\Sij{ml($32$)}{c}$
than for the other transmission operators,
in accordance with the remarks drawn above.
The novel operators will therefore lead to the fastest computations
only when the reduction of the iteration count is sufficiently high,
when compared with the transmission conditions optimized for unbounded problems
(see related discussion in sections~\ref{sec:cryo}).

\subsection{Spectrum of the iteration operators}
\label{sec:validation:spectrum}
In this section, let us briefly discuss the spectra
of the discretized iteration operator $\mathcal{F} = \Ident-\OpA$,
which we will refer to as $F$ from now on\footnote{We use
  a calligraphic (resp. non-calligraphic) typeface
  for continuous (resp. discrete) operators.}.
As the explicit construction of $F$ is computationally heavy,
only the case of $D=2$ subdomains is considered here.
It is also important to stress that $F$
is not a normal matrix\footnote{That is $FF^H \neq F^HF$,
  where $F^H$ is the conjugate transpose of $F$.},
as it can be directly seen from the formal expression of $F$
in the case of two subdomains
(see~\cite[section~2.3.1]{Dolean2015a} for instance).
As a consequence, the behavior of the GMRES
cannot be predicted from the spectrum $\Eig(F)$
of the system matrix $F$~\cite{Greenbaum1996}.
Nonetheless, eigenvalues that are clustered near $1$ are
a good indicator that those modes are well captured
by a given transmission operator.

\begin{figure}[ht]
  \centering
  \includegraphics{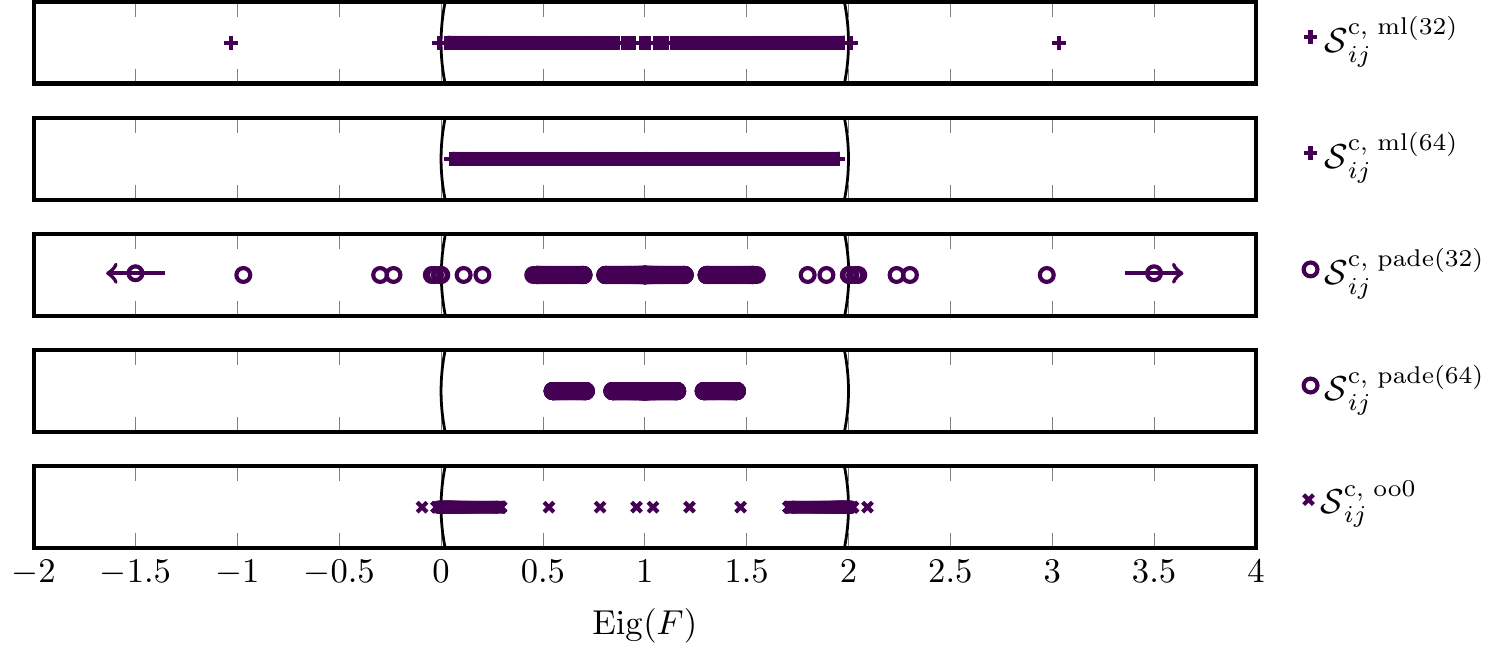}
  \caption{Spectrum of the discrete iteration operator $F$
    for the different transmission operators optimized for rectangular cavities
    -- rectangular cavity with $D=2$ subdomains
    (note that
    \textit{i)}~the spectrum is real,
    \textit{ii)}~the black arcs crossing $0$ and $2$ represent
                 a portion of the unit circle centered around $1$ and
    \textit{iii)}~the arrows indicate that a few eigenvalues
                  are outside the range shown here).}
  \label{fig:spect:closed}
\end{figure}

The spectra in \Fig{fig:spect:closed} are associated with the novel
transmission operators optimized for the rectangular cavity.
Let us stress that the eigenvalues are all real in the present example.
However, the spectrum of $F$ can be complex,
even when its coefficients are all real, as it is non-normal.
With the $\Sij{ml($N$)}{c}$ and $\Sij{pade($N$)}{c}$ operators,
we can directly see from the displayed data
that for sufficiently large values of $N$ (\ie here $N=64$)
the spectrum of $F$ lies within the unit circle centered around~$1$.
This shows that those operators are good approximations of the $\DtN$ map
of the rectangular cavity.
On the other hand, only a few eigenvalues are located around~$1$
with the $\Sij{oo0}{c}$ operator.

\begin{figure}[ht]
  \centering
  \includegraphics{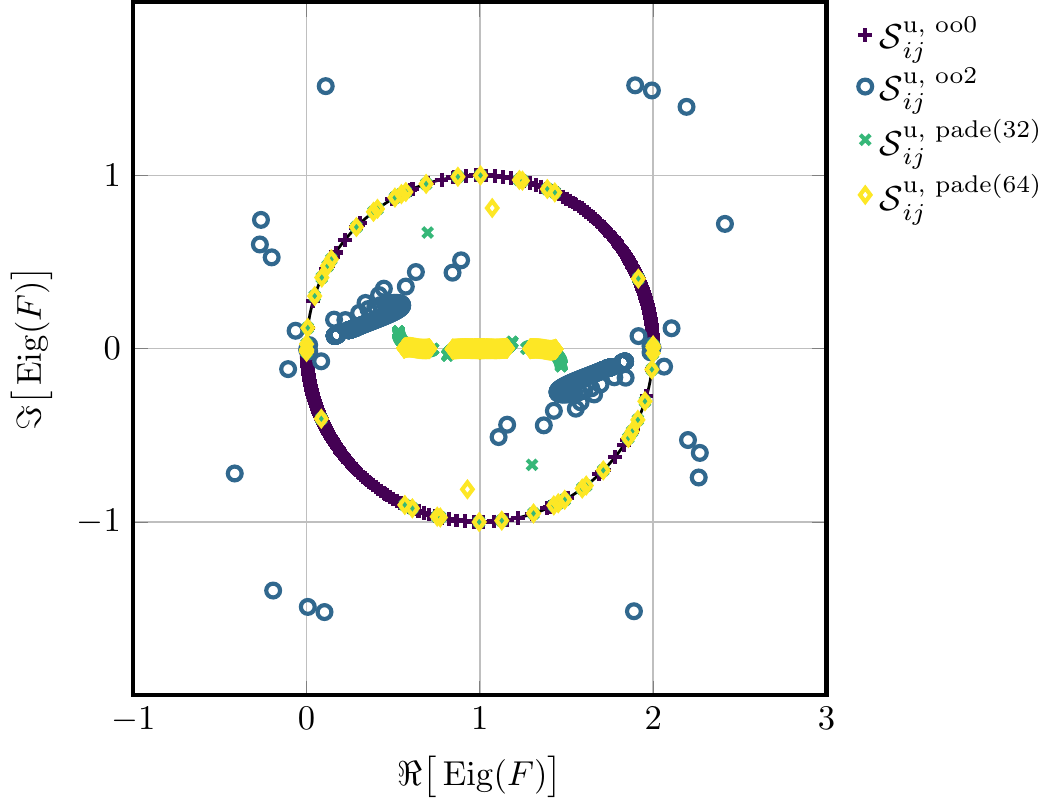}
  \caption{Spectrum of the discrete iteration operator $F$
    for the different transmission operators optimized for unbounded problems
    -- rectangular cavity with $D=2$ subdomains.}
  \label{fig:spect:open}
\end{figure}

When using a transmission operator optimized for unbounded problem
in a cavity setting,
the spectra in \Fig{fig:spect:open} are obtained.
As expected from the analysis of section~\ref{sec:open},
the $\Sij{oo0}{u}$ operator leads to eigenvalues
lying on the unit circle $\mathcal{C}$ centered around~$1$.
In addition, the behavior of the \tc{pade}{u} operators
has been well anticipated by our previous discussions as well.
It is indeed clear that the eigenvalues of $F$ fall into two categories:
\begin{enumerate}
\item those associated with evanescent modes,
  which form a cluster \emph{similar to the \tc{pade}{c} one} and
\item those associated with non-evanescent modes,
  that lie on $\mathcal{C}$.
\end{enumerate}
To conclude this subsection, let us also stress that despite a clustering
significantly better for the \tc{pade}{u} operators than for \tc{oo2}{u},
both operators exhibit rather similar GMRES convergence curves,
as shown in \Fig{fig:val:d2}.
This shows the difficulty in predicting the convergence of GMRES
from $\Eig(F)$, since $F$ is non-normal.

\subsection{Regularized and mixed transmission conditions}
\label{sec:validation:reg}
\begin{figure}[ht]
  \centering
  \begin{subfigure}[b]{0.49\textwidth}
    \centering
    \includegraphics{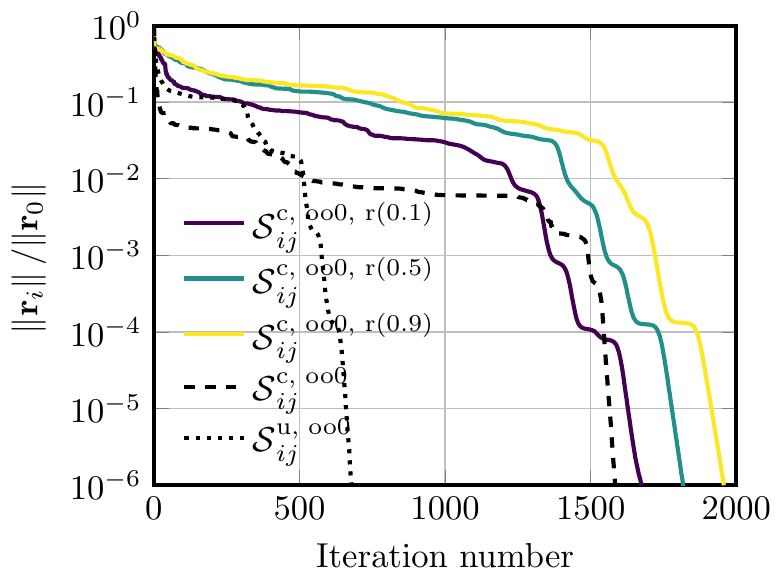}
    \caption{Regularization with a constant imaginary part.}
    \label{fig:val:reg:oo0}
  \end{subfigure}
  \begin{subfigure}[b]{0.49\textwidth}
    \centering
    \includegraphics{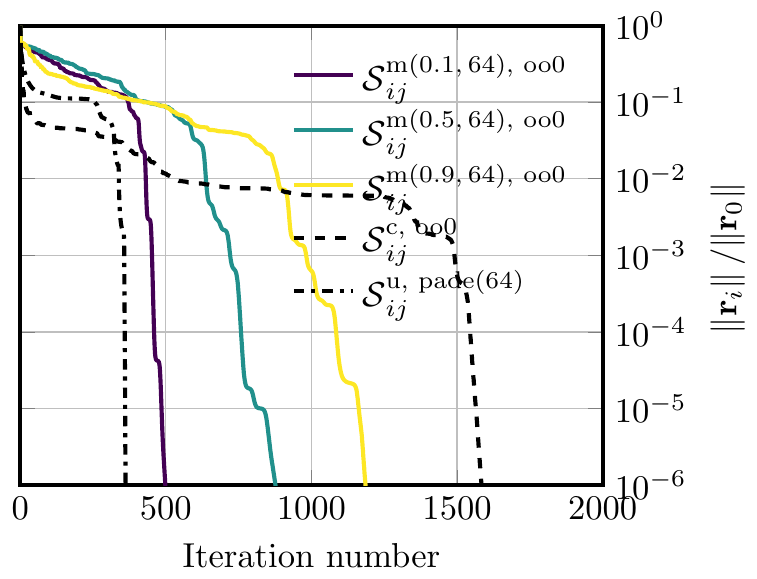}
    \caption{Regularization with a mixed operator.}
    \label{fig:val:mix:oo0}
  \end{subfigure}
  \caption{Regularized and mixed \tc{oo0}{c}
    -- rectangular cavity with $D=8$ subdomains.}
\end{figure}

As discussed in sections~\ref{sec:tc:reg} and~\ref{sec:tc:mixed},
a regularization term or a mixed transmission condition can be used
to prevent the convergence radius from becoming very large
(or, in the worst case, ill-defined).
The performance of those transmission conditions is shown
in \Figs{fig:val:reg:oo0} and~\ref{fig:val:mix:oo0}
for $\Sij{oo0, r($\chi$)}{c}$ and $\Sij{oo0}{m($\epsilon, M$)}$
respectively with $D=8$ subdomains.
It is clear from the figures that regularizing the \tc{oo0}{c}
operator by adding an imaginary part
comes at the cost of an increased number of iterations.
On the other hand, mixing \tc{oo0}{c} with \tc{pade}{u}
leads to an improvement in terms of iteration count.
Nonetheless, in all cases,
the regularized operators do not lead to any improvement
with respect to the \tc{pade}{u} conditions.
Let us also note that the regularized operators tend to the original \tc{oo0}{c}
as $\chi$ becomes small (resp. $\epsilon$ becomes large).

\begin{figure}[ht]
  \centering
  \begin{subfigure}[b]{0.49\textwidth}
    \centering
    \includegraphics{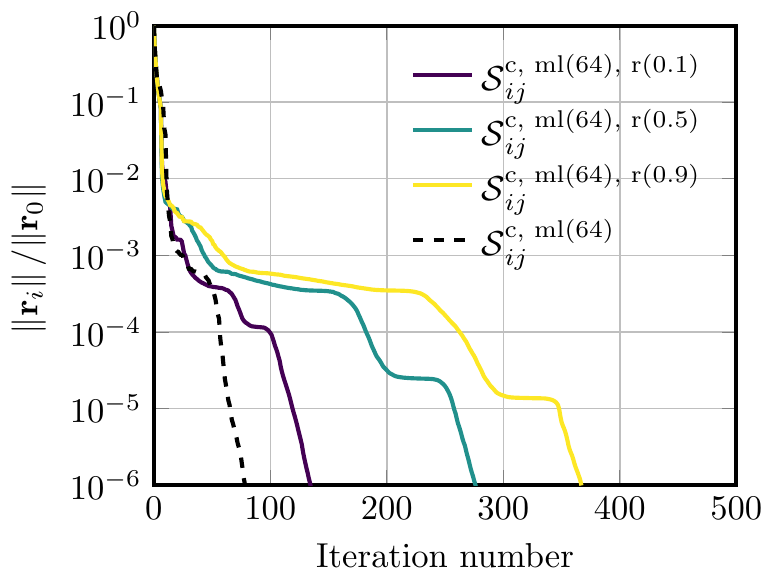}
    \caption{Regularization with a constant imaginary part.}
    \label{fig:val:reg:ml}
  \end{subfigure}
  \begin{subfigure}[b]{0.49\textwidth}
    \centering
    \includegraphics{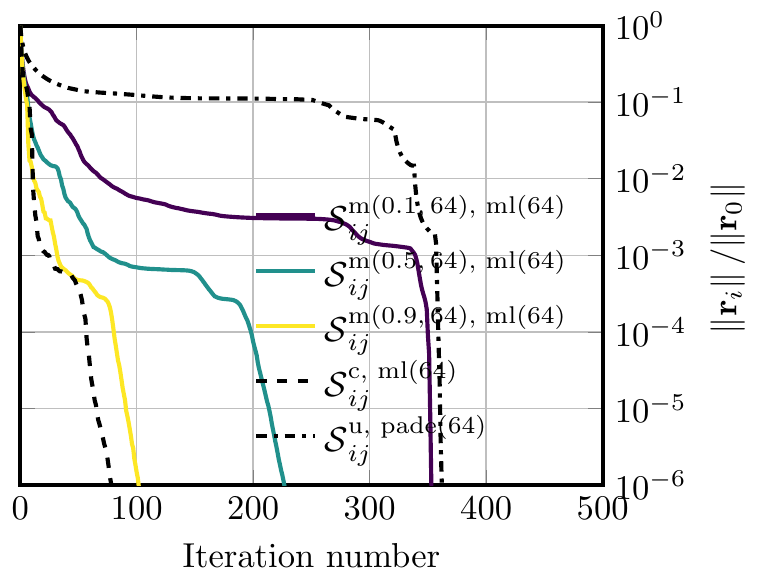}
    \caption{Regularization with a mixed operator.}
    \label{fig:val:mix:ml}
  \end{subfigure}
  \caption{Regularized and mixed \tc{ml}{c}
    -- rectangular cavity with $D=8$ subdomains.}
    \label{fig:val:regmix:ml}
\end{figure}

A similar numerical experiment can also be carried out for the regularization of
the \tc{ml}{c} condition.
From the GMRES convergence histories depicted in \Fig{fig:val:regmix:ml},
it can be noticed that while the performance of \tc{ml}{c}
remains better than \tc{pade}{u},
the regularization increases the number of iterations required
to reach convergence.
Nonetheless, this increase declines as the regularization becomes lighter,
\ie when $\chi$ (resp. $\epsilon$) becomes small (resp. large).

\begin{figure}[ht]
  \centering
  \begin{subfigure}[b]{0.49\textwidth}
    \centering
    \includegraphics{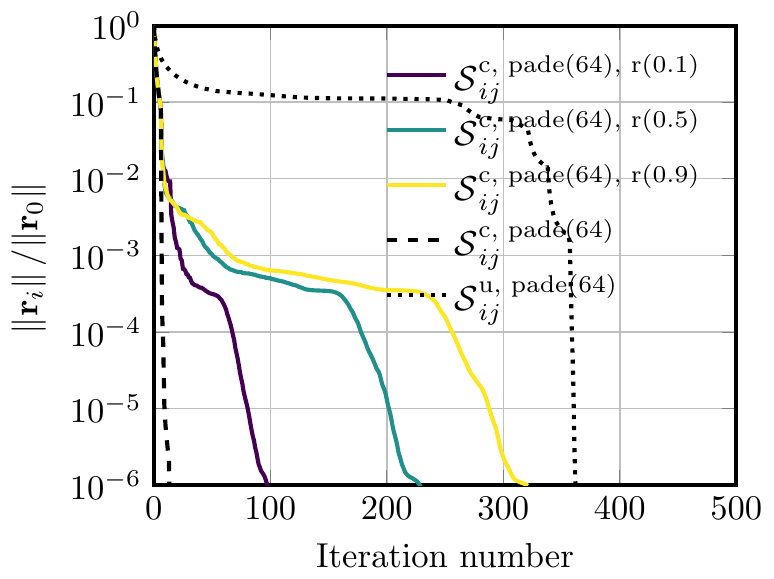}
    \caption{Regularization with a constant imaginary part.}
    \label{fig:val:reg:pade}
  \end{subfigure}
  \begin{subfigure}[b]{0.49\textwidth}
    \centering
    \includegraphics{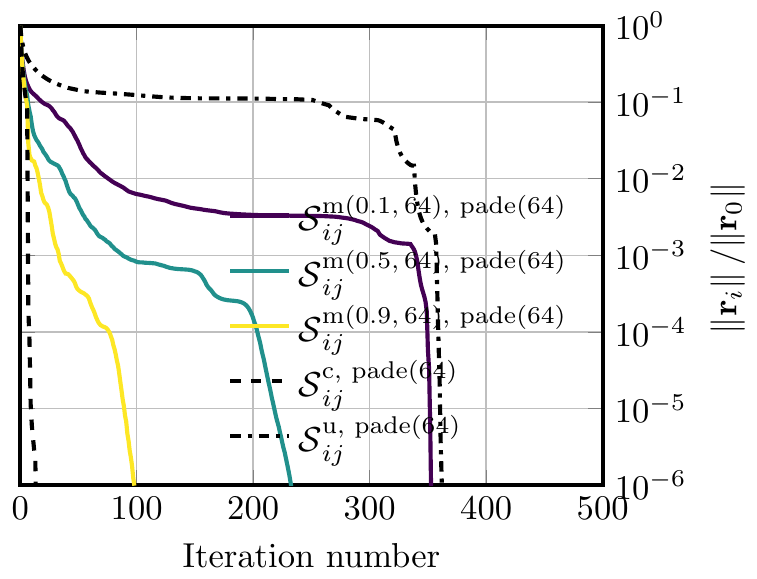}
    \caption{Regularization with a mixed operator.}
    \label{fig:val:mix:pade}
  \end{subfigure}
  \caption{Regularized and mixed \tc{pade}{c}
    -- rectangular cavity with $D=8$ subdomains.}
    \label{fig:val:regmix:pade}
\end{figure}

An analogous numerical experiment is performed once more
for the regularization of the \tc{pade}{c} condition.
This last scenario shows a behavior similar to \tc{ml}{c},
as it can be directly observed in \Fig{fig:val:regmix:pade}.

\clearpage
\subsection{Impact of the number of auxiliary unknowns}
\label{sec:validation:naux}
\begin{figure}[ht]
  \centering
  \begin{subfigure}{1\textwidth}
    \centering
    \includegraphics{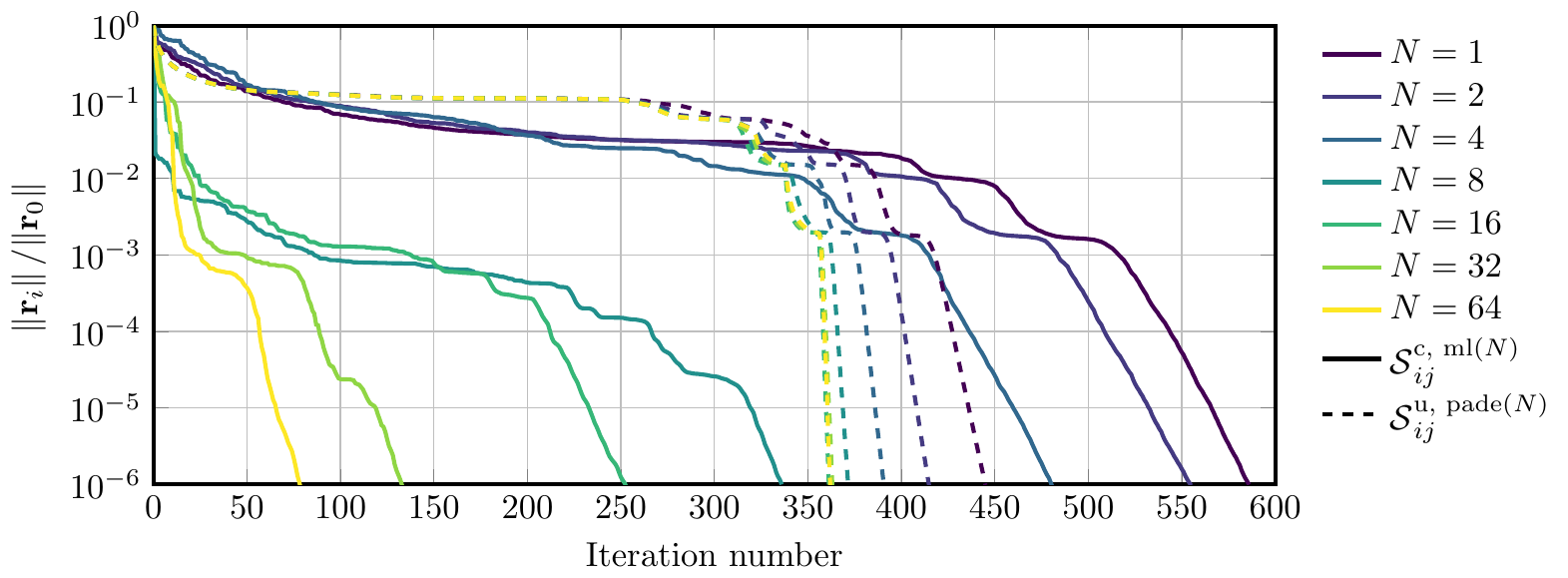}
    \caption{Mittag-Leffler expansion.}
  \end{subfigure}
  \begin{subfigure}{1\textwidth}
    \centering
    \includegraphics{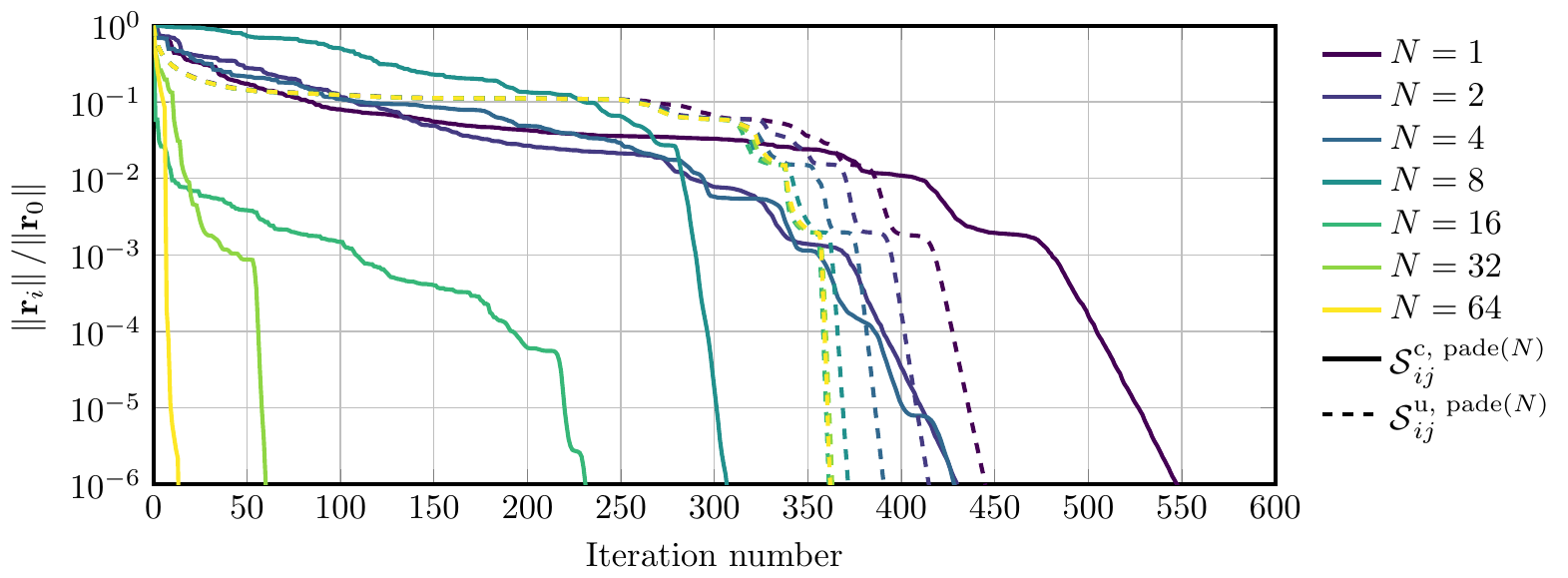}
    \caption{Pad\'e expansion.}
  \end{subfigure}
  \caption{Convergence history of the GMRES solver
    for different number of auxiliary unknown in the case of $D=8$ subdomains.}
  \label{fig:val:naux}
\end{figure}

Let us now investigate the impact of the number of auxiliary unknowns $N$
appearing in
$\Sij{pade($N$)}{u}$, $\Sij{ml($N$)}{c}$ and $\Sij{pade($N$)}{c}$
on the number of iteration of the GMRES solver.
The numerical results are shown in \Fig{fig:val:naux}
for a numerical experiment involving $D=8$ subdomains of equal size.
It is clear from the data that every operator converges
with values of $N$ as low as $N=1$.
Nonetheless, for the novel operators to outperform the convergence
of the \tc{pade}{u} operator,
a minimum value of $N_\text{min}^\text{data}=8$ is required in this example.
It is also evident in this example that $\Sij{pade($N$)}{c}$
always performs better than $\Sij{ml($N$)}{c}$ for a fixed value of $N$.
Last but not least, it is evident from \Fig{fig:val:naux} that
$N$ has only a mild effect on the performance of $\Sij{pade($N$)}{u}$.

It is also interesting to
compare the value of $N_\text{min}^\text{data}$ found experimentally
and the value predicted with the pole criterion discussed
in section~\ref{sec:tc:nmin}.
In this numerical experiment,
the largest $\ell_{ij}$ is $\ell_{ij}^\text{max}=\frac{7}{8}\ell$ and thus
$N_\text{min}^\text{pole} = 44$
according to the pole criterion~\eqref{eq:closed:nmin}.
Consequently,
$N_\text{min}^\text{pole}$ is a pessimistic estimation in this case,
since $N_\text{min}^\text{data} < N_\text{min}^\text{pole}$.

\subsection{Increase in the number of rectangular subdomains}
\label{sec:val:dmany}
\begin{figure}[ht]
  \begin{subfigure}{1\textwidth}
    \centering
    \includegraphics{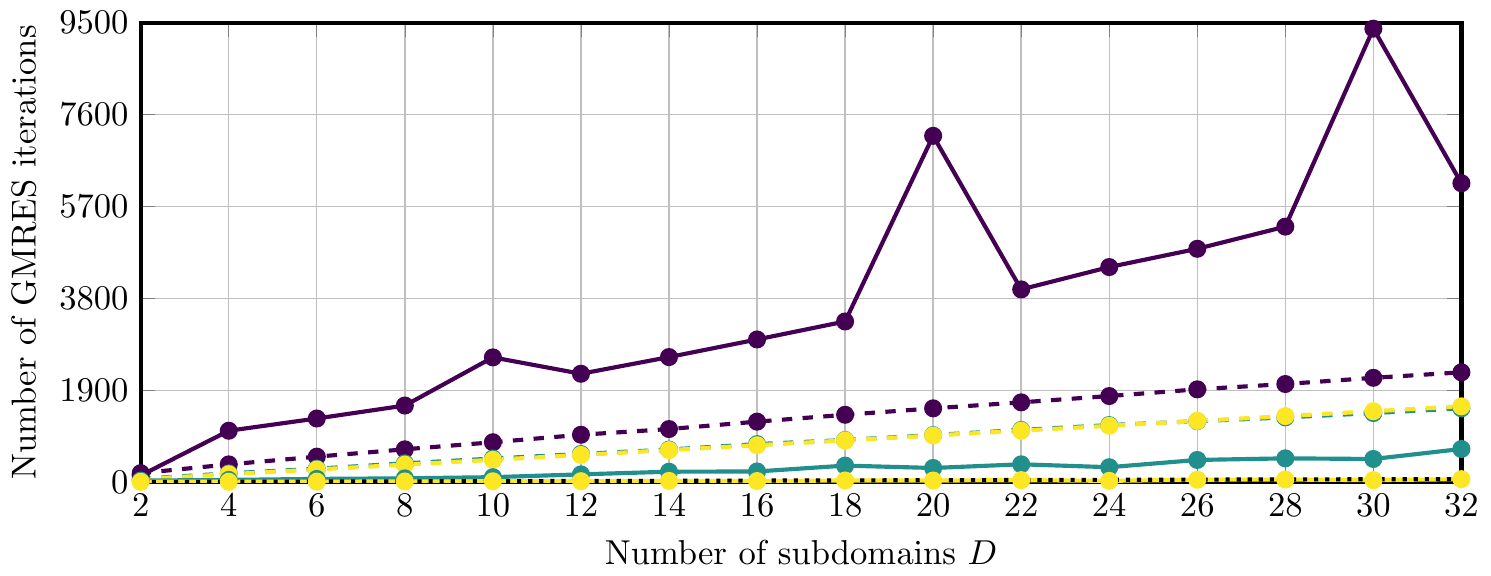}
    \caption{Zoom level 0: focus on \tc{oo0}{c}.}
  \end{subfigure}
  \begin{subfigure}{1\textwidth}
    \centering
    \includegraphics{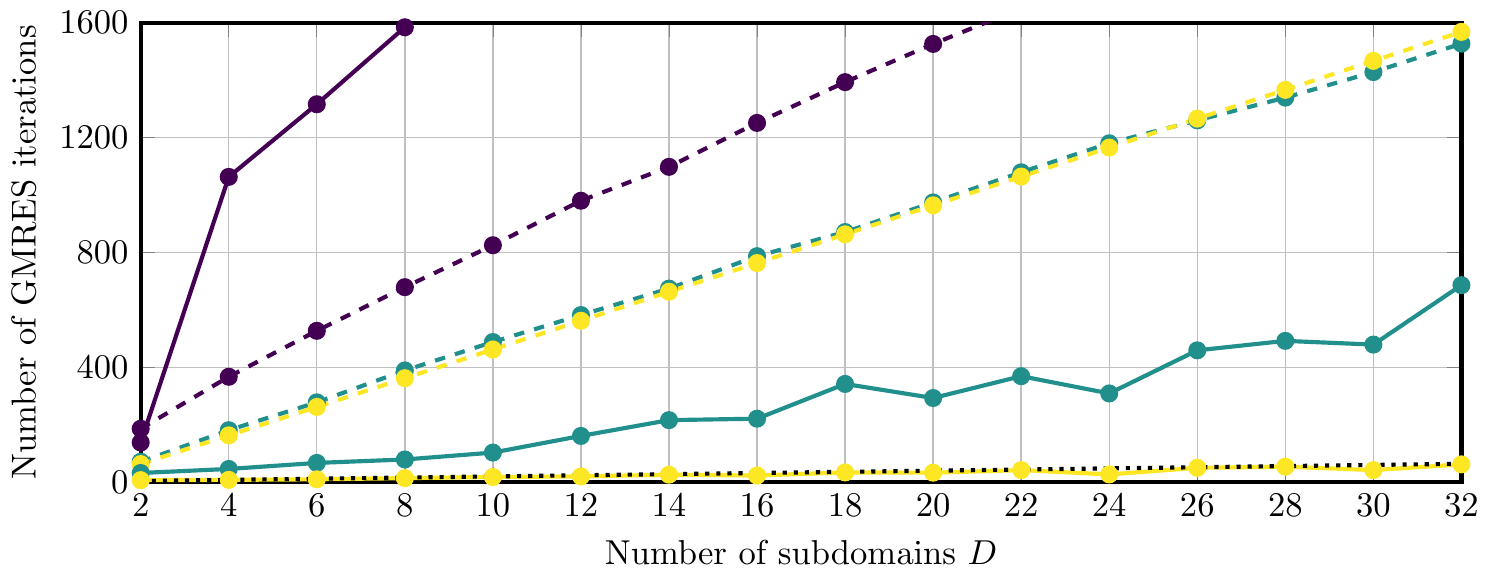}
    \caption{Zoom level 1: focus on \tc{oo0}{u}, \tc{oo2}{u},
      \tc{pade}{u} and \tc{ml}{c}.}
  \end{subfigure}
  \begin{subfigure}{1\textwidth}
    \centering
    \includegraphics{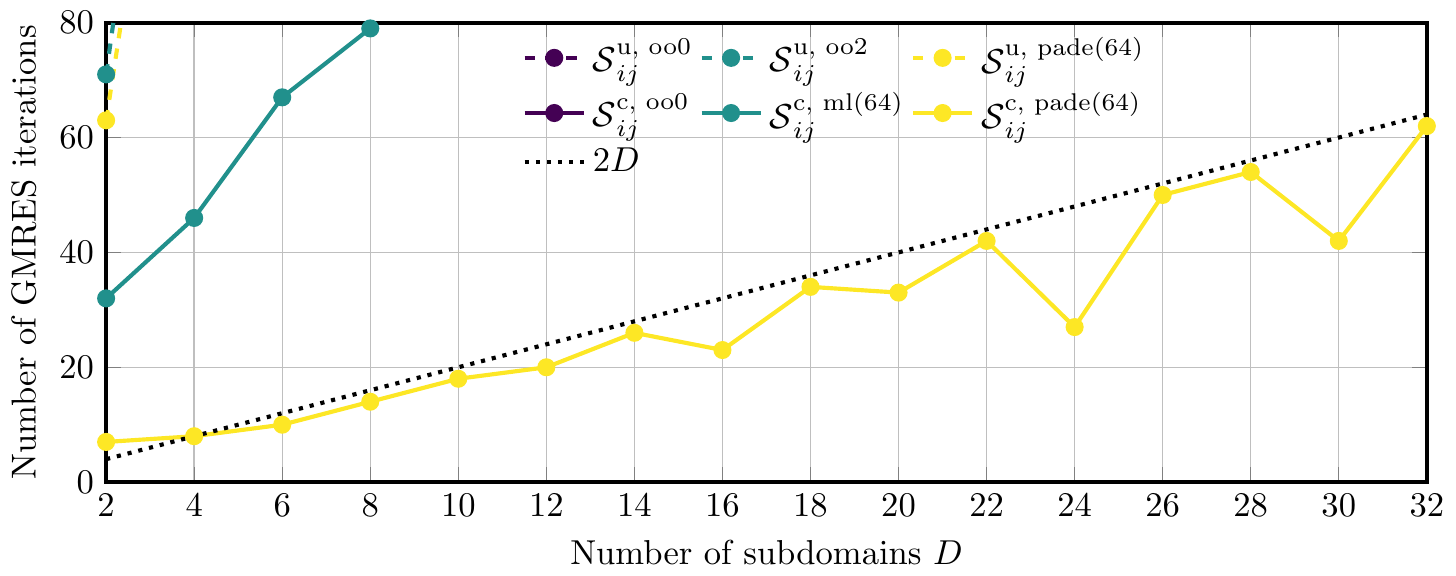}
    \caption{Zoom level 2: focus on \tc{pade}{c}}
  \end{subfigure}
  \caption{Number of GMRES iterations
    as a function of the number of subdomains
    (raw numbers available in \Tab{tab:val:dmany}).}
  \label{fig:val:dmany}
\end{figure}

In this section we focus on the impact of the number of subdomains $D$
onto the GMRES iteration count, as shown in \Fig{fig:val:dmany}.
From this plot, is clear that $\Sij{pade($64$)}{c}$ leads to an
increase in the number of GMRES iterations with the optimal slope of $2D$,
at least for the considered range of $D$.
Let us mention also that the unbounded transmission operators
exhibit a significantly larger slope,
motivating thus the use of transmission conditions specifically devised for
cavity problems.
Furthermore let us note that $\Sij{pade($64$)}{u}$
and $\Sij{oo2}{u}$ present the same slope,
as they are both excellent localization of the unbounded $\DtN$ map.
In addition, while the $\Sij{ml($64$)}{c}$ operator shows
a suboptimal scaling with respect to $D$,
it outperforms $\Sij{pade($64$)}{u}$ and $\Sij{oo2}{u}$.
Last but not least, it is evident from \Fig{fig:val:dmany}
that the scaling behavior of $\Sij{oo0}{c}$ is the worst of all.
\afterpage{\clearpage}

\subsection{Impact of the length-to-wavelength ratio}
\label{sec:val:k}
\begin{figure}[ht]
  \begin{subfigure}{1\textwidth}
    \centering
    \includegraphics{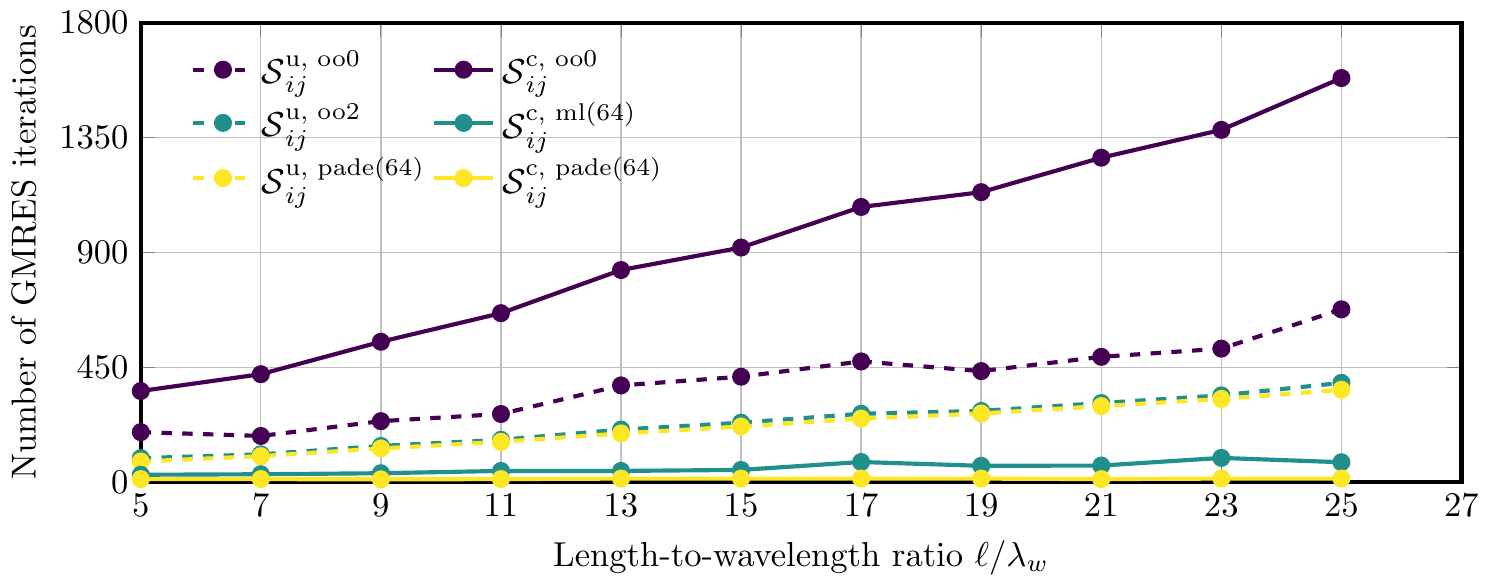}
    \caption{Zoom level 0: focus on \tc{oo0}{c}.}
  \end{subfigure}
  \begin{subfigure}{1\textwidth}
    \centering
    \includegraphics{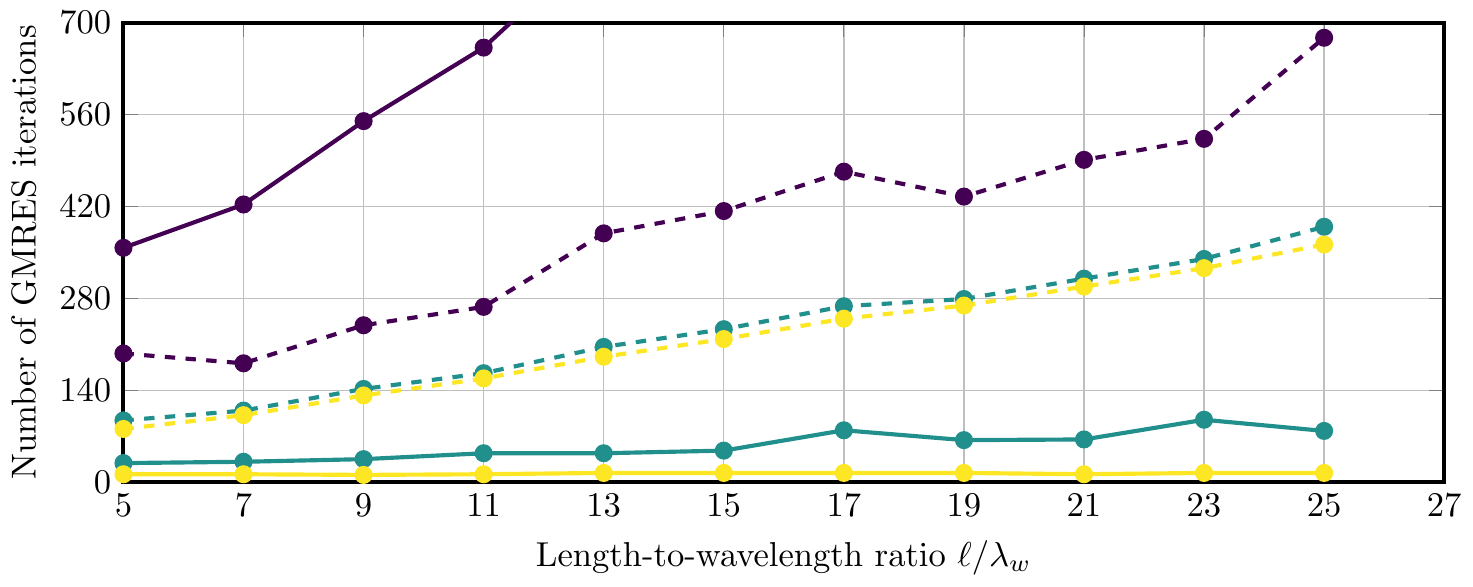}
    \caption{Zoom level 1: focus on \tc{oo0}{u}, \tc{oo2}{u} and \tc{pade}{u}.}
  \end{subfigure}
  \begin{subfigure}{1\textwidth}
    \centering
    \includegraphics{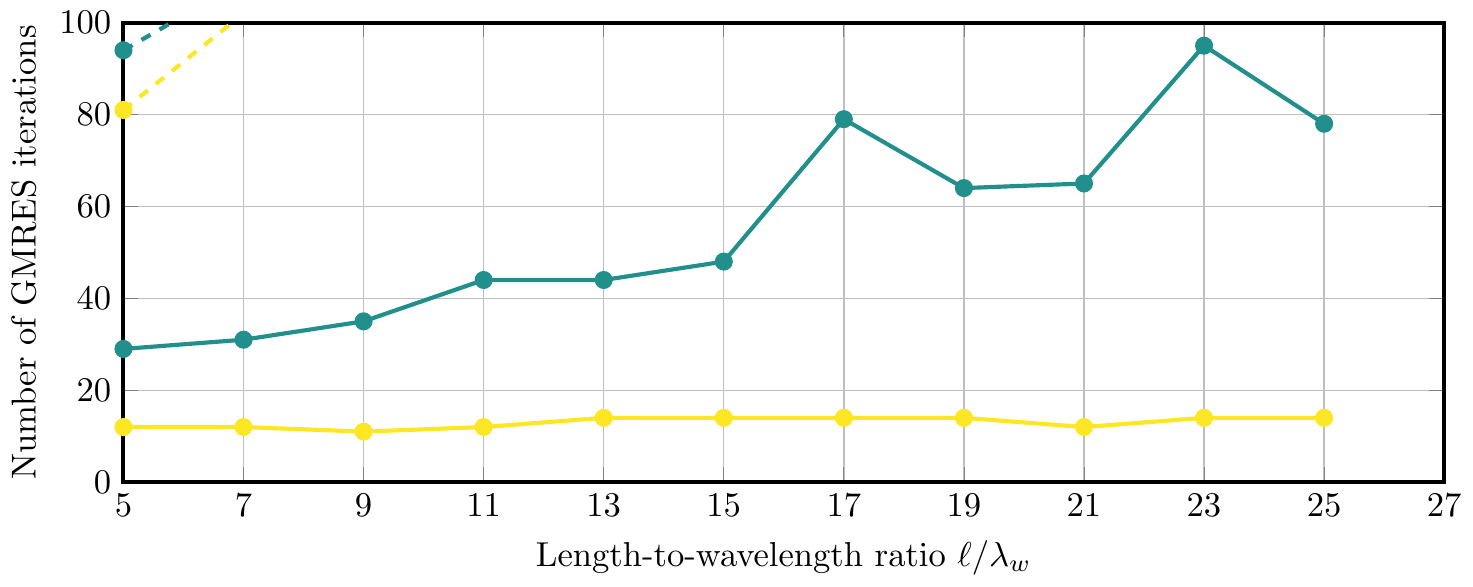}
    \caption{Zoom level 2: focus on \tc{ml}{c} and \tc{pade}{c}.}
  \end{subfigure}
  \caption{Number of GMRES iterations as a function
    of the length-to-wavelength ratio for a cavity with $D=8$ subdomains
    (raw numbers available in \Tab{tab:val:k}).}
  \label{fig:val:k}
\end{figure}

The dependence of the GMRES iteration count on the wavenumber
is a major performance indicator of an OS scheme,
and a numerical experiment analyzing this is therefore carried out.
In this investigation the values of $\ell/\lambda_w$ are chosen
close to an integer value, which corresponds to a cavity driven
at a frequency close to one of its resonance frequencies.
The computational domain is partitioned into $D=8$ subdomains
and the cavity is excited with double the amount of non-evanescent modes
(this number depending on  $\ell/\lambda_w$).
It is clear from the data shown in \Fig{fig:val:k}
that the iteration count increases rapidly with $\ell/\lambda_w$
with the unbounded transmission operators.
This growth in the iteration count is, nonetheless,
not as fast as for $\Sij{oo0}{c}$.
On the other hand,
the situation is significantly improved with $\Sij{ml($N$)}{c}$
and the iteration count becomes almost independent from
$\ell/\lambda_w$ with $\Sij{pade($N$)}{c}$,
at least in the considered $\ell/\lambda_w$ range.
Let us also note that similar results are obtained when $\ell/\lambda_w$
is selected away from a resonance.
\afterpage{\clearpage}

\subsection{Discussion}
From the data gathered in the previous numerical examples,
it is obvious that the $\Sij{oo0}{c}$ operator leads
to a poor scaling when $D>2$.
For this reason, this operator and its regularized variants
are not further considered in this work.

It is also clear that for rectangular cavities
the $\Sij{pade($N$)}{c}$ and $\Sij{ml($N$)}{c}$
operators converge with significantly less iterations than their
unbounded alternatives, \ie $\Sij{pade($N$)}{u}$, $\Sij{oo2}{u}$
and $\Sij{oo0}{u}$,
and exhibit a better scaling with increasing values of $D$.
Nonetheless, let us recall that the novel operators
\begin{myenum*}
\item are associated with a higher computational cost and
\item are more sensitive to numerical errors in finite precision arithmetic.
\end{myenum*}
This latter point motivates the use of the
\emph{modified} version of the Gram-Schmidt algorithm
in the GMRES orthogonalization step.

It is also worth noticing that the unbounded operator
$\Sij{pade($N$)}{u}$ and $\Sij{oo2}{u}$
lead to very similar convergence profiles,
which is not surprising as they both are excellent approximations
of the unbounded $\DtN$ map.
For this reason, only the $\Sij{pade($N$)}{u}$ operator
will be further considered.

Another interesting result concerns the impact of a change of the wavenumber
on the convergence profile of the above transmission operators.
For the considered range of length-to-wavelength ratios $\ell/\lambda_w$,
we observed that the $\Sij{pade($N$)}{c}$
behaves quasi-independently from $\ell/\lambda_w$,
while the other operators exhibit a (quasi) linear increase
of the iteration count with $\ell/\lambda_w$.
Nonetheless, the slope of this increase remains small for
$\Sij{ml($N$)}{c}$ when compared with
$\Sij{oo0}{u}$, $\Sij{oo2}{u}$ and $\Sij{pade($N$)}{u}$.

Let us finally note that similar experiments were carried out
on three-dimensional rectangular parallelepipedic cavities
with a square cross-section of $h\times{}h$ and a length of $\ell$.
Results analogous to the above two-dimensional test cases were obtained.

\newpage
\section{Sensitivity to geometrical parameters}
\label{sec:sensitivity}
The previous section considered only one geometry
(a rectangular cavity and its three-dimensional equivalent),
which refers to the specific configuration used to devise
the transmission operators discusses in section~\ref{sec:tc}.
Nonetheless, relevant simulations involve computational domains
that differ from this original setting.
We therefore discuss in this section two geometries
deviating from the canonical one.

\subsection{Trapezoidal geometry -- modified Gram-Schmidt variant}
\label{sec:sensitivity:trapeze:mgs}
\begin{figure}[ht]
  \centering
  \includegraphics{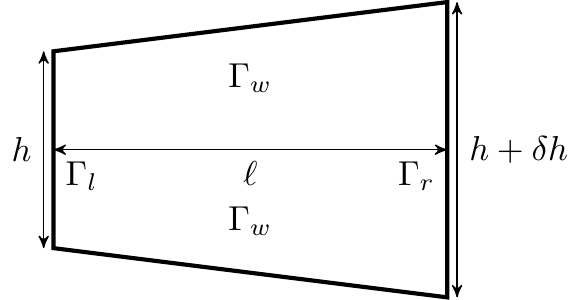}
  \caption{Isosceles trapezoid considered in the numerical experiment.}
  \label{fig:trapeze}
\end{figure}

For this first numerical experiment, let us consider an isosceles trapezoid
whose base is characterized with a parameter $\delta$,
as the one depicted in \Fig{fig:trapeze}.
This computational domain is partitioned into $D=16$ subdomains
and, as with the previous rectangular cavity,
the aspect ratio $\ell/h$ equals $2$
and the length-to-wavelength ratio $\ell/\lambda_w$ is approximately $25.001$.

\begin{figure}[ht]
  \centering
  \includegraphics{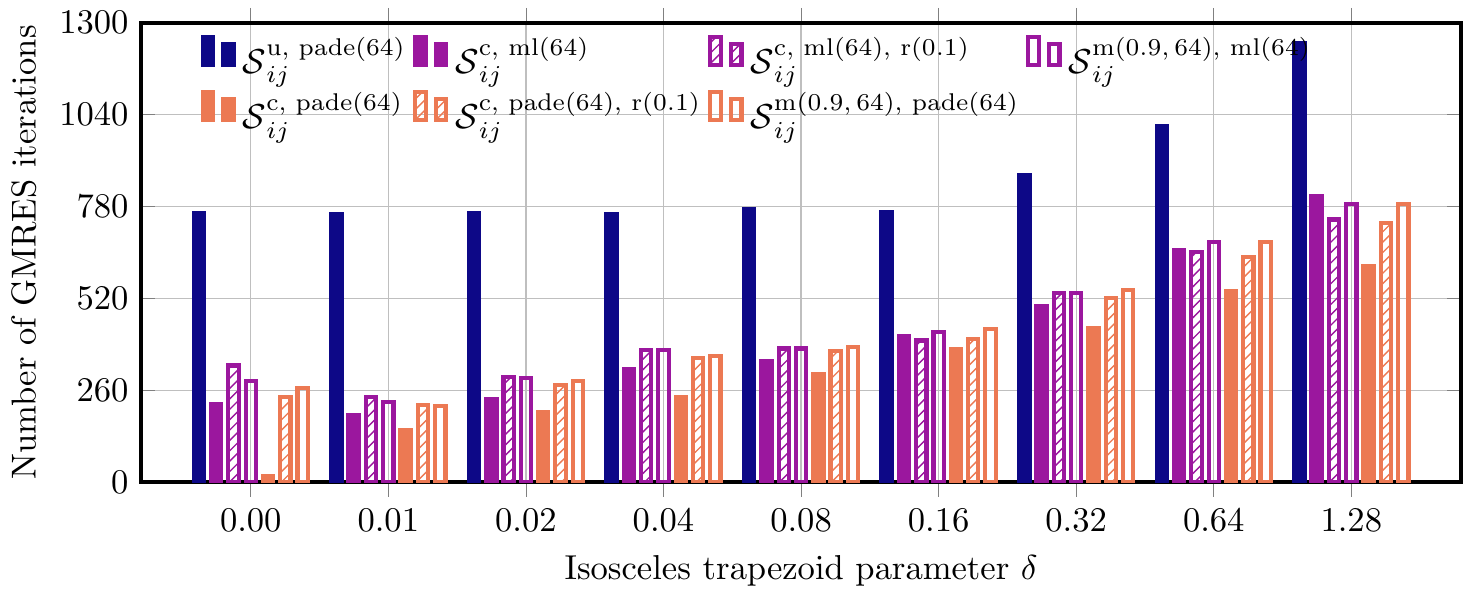}
  \caption{Iteration count of the GMRES solver
    when used on a trapezoidal geometry partitioned into $D=16$ subdomains
    (raw numbers available in \Tab{tab:trapeze:data}).}
  \label{fig:trapeze:data}
\end{figure}

The number of GMRES iterations associated with various operators
are shown in \Fig{fig:trapeze:data} for different values of $\delta$.
It is clear from the depicted results that
the novel rational operators outperform, in all considered cases,
the $\Sij{pade($64$)}{u}$ one in terms of iteration count,
\emph{even when no regularization is applied}.
Additionally, the $\Sij{pade($64$)}{c}$ is systematically the best choice
and its regularization always leads to a slight increase in the iteration count.
In the case of the \tc{ml}{c} variant,
regularization slightly improves the convergence
for the highest value of $\delta$.
Furthermore, there is an evident trend in the displayed data
showing an increase in the iteration count as $\delta$ increases.
However, as this increase also affects $\Sij{pade($64$)}{u}$,
it is hard to predict if and when this operator will overtake the novel ones,
at least in terms of iteration count.

\subsection{Trapezoidal geometry -- classical Gram-Schmidt variant}
In the introduction of section~\ref{sec:validation},
we mentioned that using the \emph{modified} version
of the Gram-Schmidt (GS) algorithm is critical for the convergence,
as the novel transmission conditions involve operators
oscillating rapidly in a wide range.
In order to assess the importance of this choice,
we carry out once more the previous numerical experiments using the
\emph{classical} version of the GS orthogonalization procedure~\cite{Saad2003}.

\begin{figure}[ht]
  \centering
  \includegraphics{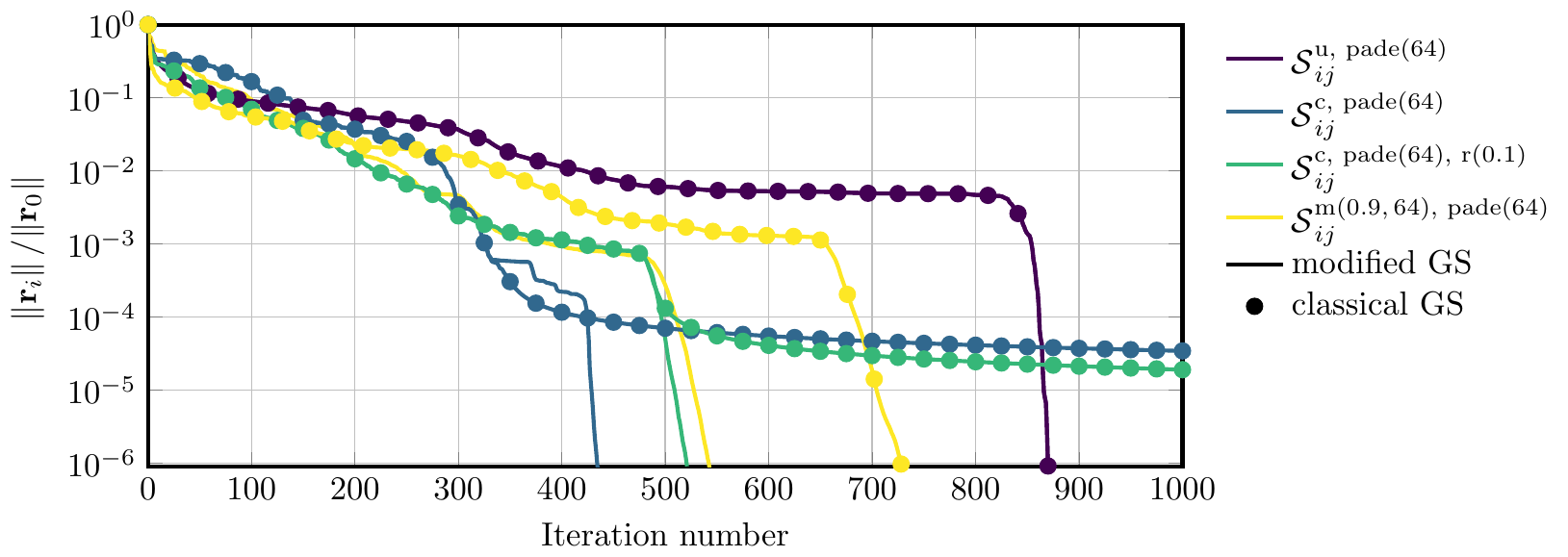}
  \caption{Convergence history of the GMRES solver
    with the classical and the modified Gram-Schmidt (GS) procedures,
    when used on a trapezoidal geometry partitioned into $D=16$ subdomains.}
  \label{fig:trapeze:gs}
\end{figure}

As it would be impractical to show the data for all possible combinations,
only the case with \tc{pade}{u} and \tc{pade}{c}
(with and without regularization) at $\delta=32$
is shown in \Fig{fig:trapeze:gs},
as it contains the different behavior we want to highlight.
First of all, it is clear from the displayed data that
the modified GS does not impact the behavior of the \tc{pade}{u} operator.
In addition, it is also evident that it prevents stagnation
and improves the overall convergence with \tc{pade}{c}.

Regarding the other combinations that are not shown in \Fig{fig:trapeze:gs},
we observed that:
\begin{myenum*}
\item a relative residual $\norm{\vec{r}_i}/\norm{\vec{r}_0}$
  smaller than $10^{-6}$ is not always reached
  when using the classical Gram-Schmidt procedure,
\item this relative residual is always reached with the modified GS,
\item the convergence history of the \tc{pade}{u} operator is identical
  for both GS variants and
\item the modified GS is never worse than the classical GS.
\end{myenum*}

\subsection{Rectangular cavity involving obstacles
  -- impact of the number of obstacles}
\label{sec:sensitivity:obstacles:N}
\begin{figure}[ht]
  \centering
  \includegraphics{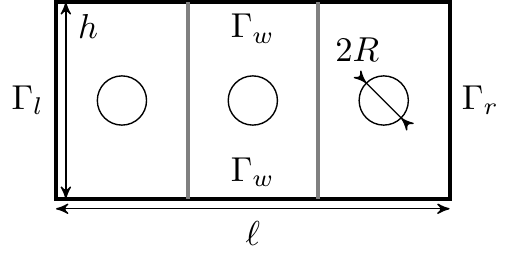}
  \caption{Rectangular cavity involving $O=3$ obstacles and $D=3$ subdomains.}
  \label{fig:obstacles}
\end{figure}

Let us now study again a rectangular cavity,
but let us introduce $O$ circular obstacles in the domain.
As in the previous case,
we have $\ell/h = 2$ and $\ell/\lambda_w\approx 25.001$.
Furthermore, we assume that the obstacles exhibit a hard-wall behavior
and that each subdomain includes
\emph{at most} one obstacle located in its center.
A sketch of a possible configuration is shown in \Fig{fig:obstacles},
with $O=3$ obstacles and $D=3$ subdomains.

\begin{figure}[ht]
  \centering
  \includegraphics{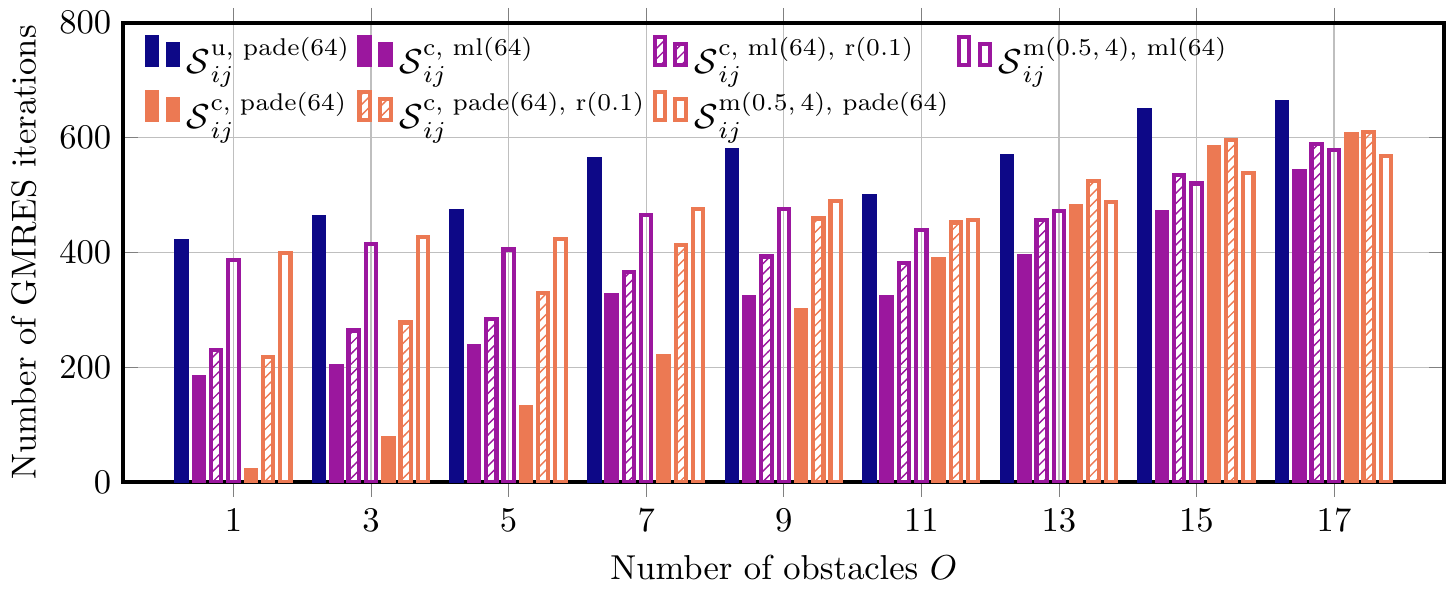}
  \caption{Iteration count of the GMRES solver
    when used on a rectangular cavity involving $O$ obstacles
    of radius $R=0.5\lambda_w$ and $D=17$ subdomains
    (raw numbers available in \Tab{tab:obstacles:N}).}
  \label{fig:obstacles:N}
\end{figure}

We carry out a first numerical experiment
consisting in determining the number of iterations
as a function of the number of obstacles $O$.
The obstacles are introduced by starting from the middle of the cavity
and then by adding them equally on both sides.
We consider a domain partitioning involving $D=17$ subdomains,
leading thus to a maximum number of obstacles of $O=17$,
and a radius of $R=0.5\lambda_w$ for the obstacles.
The results of this experiment are gathered
in \Fig{fig:obstacles:N}.
First of all, it is worth mentioning that when $O=1$,
the transmission condition $\Sij{pade($64$)}{c}$
exhibits an excellent performance compared with the other operators.
The operator $\Sij{ml($64$)}{c}$ also shows a good performance,
in comparison with the unbounded Pad\'e operator.
This lead with respect to $\Sij{pade($64$)}{u}$ decreases however
as $O$ increases.
Nonetheless, in all considered cases, the novel conditions are associated with
a lower iteration count than $\Sij{pade($64$)}{u}$,
although this gain is modest for higher values of $O$.
In this regard, the \tc{ml}{c} condition performs slightly
better for higher values of $O$,
while \tc{pade}{c} is better for the lower values.
Concerning the regularized variants,
they lead to higher iteration counts than their unregularized variants,
apart from the two highest values of $O$,
where regularization is slightly beneficial to the \tc{pade}{c} operator.
Before concluding this subsection,
let us stress that we used different values
for the number of terms $M$ and $N$ in the mixed operators
$\Sij{ml($N$)/pade($N$)}{m($\epsilon, M$)}$.

\subsection{Rectangular cavity involving obstacles
  -- impact of the size of the obstacles}
\label{sec:sensitivity:obstacles:R}
\begin{figure}[ht]
  \centering
  \includegraphics{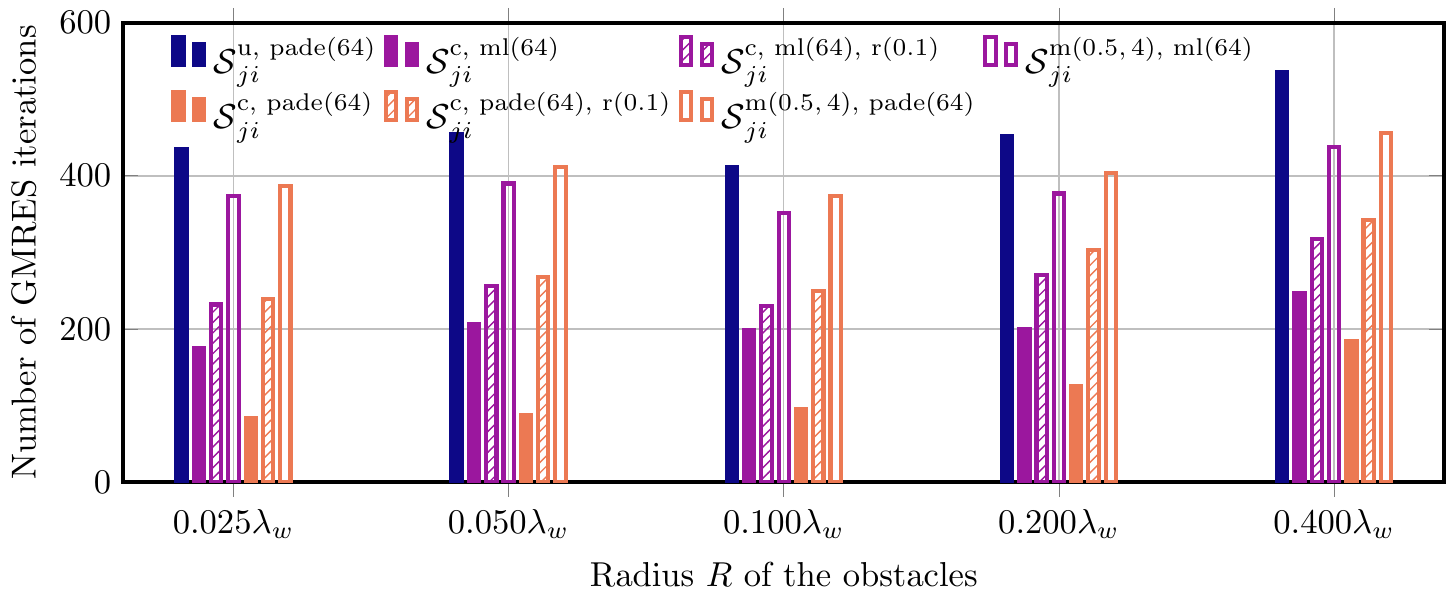}
  \caption{Iteration count of the GMRES solver
    when used on a rectangular cavity involving $O=7$ obstacles
    and $D=17$ subdomains
    (raw numbers available in \Tab{tab:obstacles:R}).}
  \label{fig:obstacles:R}
\end{figure}

In the previous subsection we assumed a rather large radius for the obstacles
compared with the wavelength
and focused only on the number of obstacles.
Let us now reverse the study
and let us determine the impact of the obstacle size
on the performance of the transmission operators.
To this end,
let us again consider a rectangular cavity with $D=17$ subdomains
and $O=7$ obstacles of radius $R$
and let us gather in \Fig{fig:obstacles:R}
the number of iterations required to solve this problem as $R$ varies.
Concerning the novel conditions,
the displayed data show a clear trend pointing toward
an increase of the iteration count as $R$ increases.
While this increase is systematic for $\Sij{pade($64$)}{c}$,
a plateau can be seen with $\Sij{ml($64$)}{c}$ for
the three intermediate values of $R$.
As in the previous numerical experiment,
let us note
\begin{myenum*}
\item that the novel conditions outperform $\Sij{pade($64$)}{u}$ in all cases
  and
\item that regularization leads to an increased iteration count.
\end{myenum*}

\subsection{Discussion}
The above numerical experiments clearly show that,
while the performance of the novel $\Sij{ml($N$)}{c}$ and $\Sij{pade($N$)}{c}$
operators deteriorates as the geometry of the problem
deviates from the reference rectangular case,
this deterioration does not exhibit sudden jumps as the geometry is modified.
In addition, the sensitivity of the novel operators
with respect to numerical errors stemming from finite precision arithmetic
is also clearly shown by our data,
motivating the use of the modified Gram-Schmidt algorithm
in the orthogonalization step of GMRES.

When compared to the \tc{pade}{u} operator,
\tc{pade}{c} and \tc{ml}{c} clearly outperform the former
for small deviations.
On the other hand, in all considered cases,
the reduction of the iteration count brought by
\tc{pade}{c} and \tc{ml}{c} can become modest (with respect to \tc{pade}{u})
for large deviations.
Nonetheless, we did not encounter situations where \tc{pade}{u}
leads to clear gain,
at least in terms of iteration count
and within the scope of the numerical experiment considered in this work.
This suggests a rather robust behavior of \tc{pade}{c} and \tc{ml}{c}.

\section{Engineering test case: acoustic noise in a three-dimensional
  model of the helium vessel of a beamline cryostat}
\label{sec:cryo}
Going back to the experiment with a rectangular cavity filled with obstacles,
it has been shown that the \tc{pade}{c} operator shows a very fast convergence
when the cavity \emph{exhibits only one obstacle},
as shown in \Fig{fig:obstacles:N}.
Additionally, the same behavior is observed for other configurations
exhibiting a unique circular obstacle.
A such characteristic can be very well suited when, for instance,
simulating the acoustic noise in the helium vessel of a beamline cryostat,
which basically consists of a rectangular cavity interrupted
by a circular obstacle (\eg see the cryostat discussed in~\cite{Haider2019}).
Such acoustic noise analyses can become critical, for instance, when designing
a cryogenic current comparator (CCC) with a large bore~\cite{Seidel2018}.
A CCC is one of the most sensitive instrument
for measuring very low electric currents with high accuracy
and can be used \eg in particle accelerators for the non-destructive monitoring
of slowly extracted charged particle beams
(current intensities below $1\mu{}A$)~\cite{Seidel2018}.

\begin{figure}[ht]
  \centering
  \includegraphics{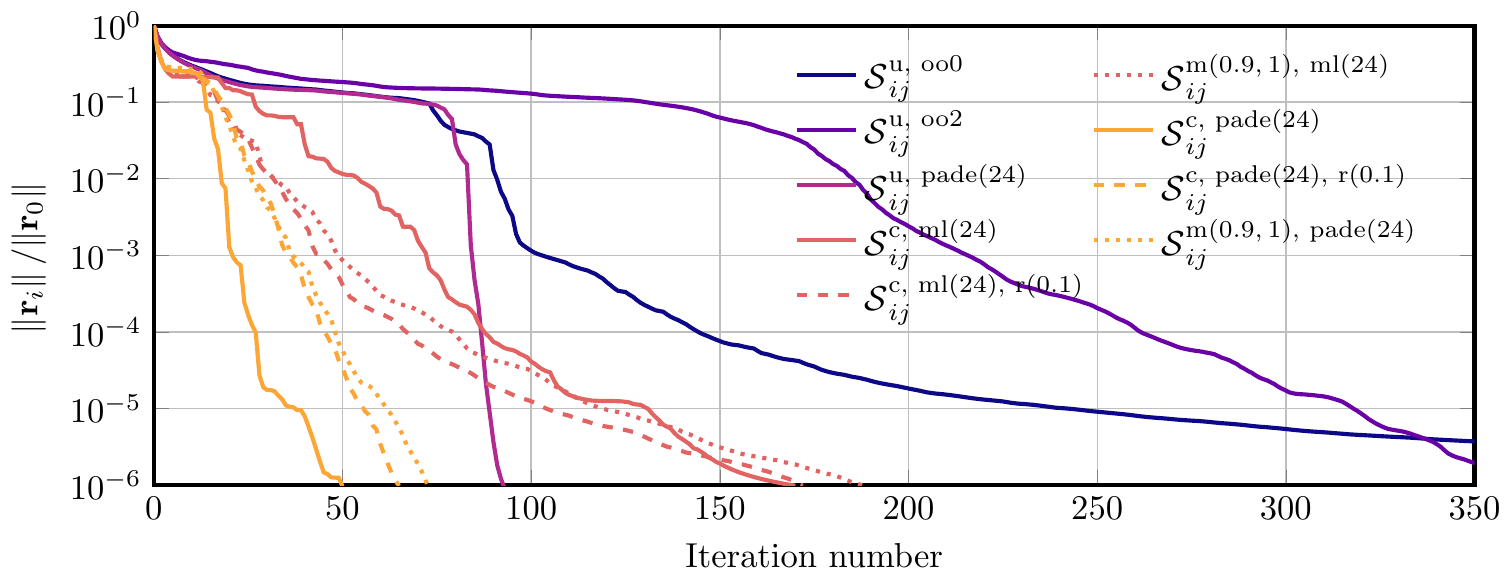}
  \caption{Convergence history of the GMRES solver
    -- helium vessel of a beamline cryostat with $D=8$ subdomains.}
  \label{fig:cryo}
\end{figure}

\begin{figure}[ht]
  \centering
  \includegraphics[width=7.5cm, angle=-90]{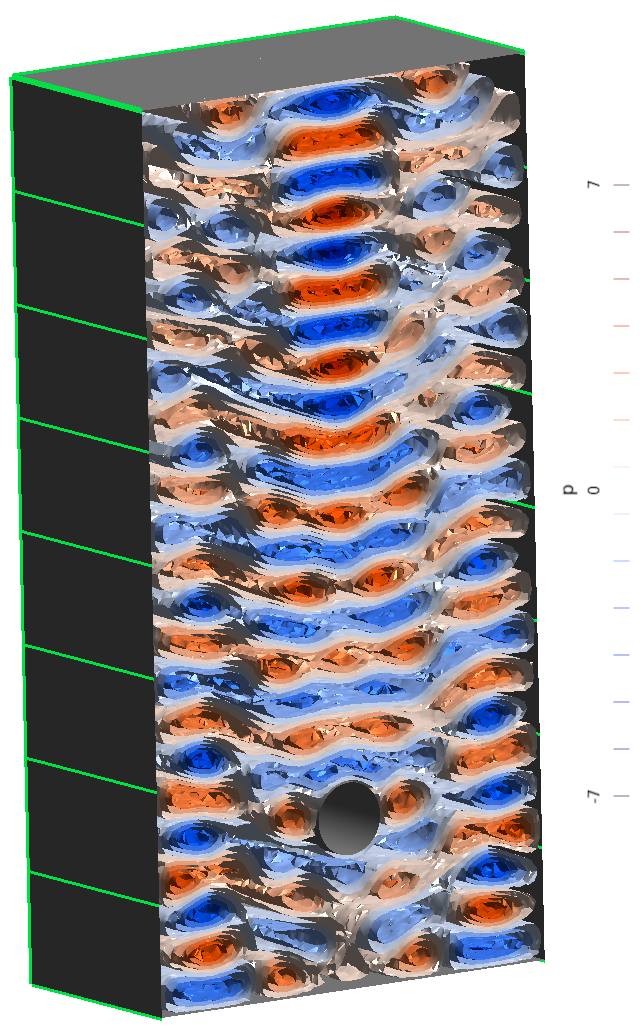}
  \caption{Field map of the simulated helium vessel
    (view rotated by $90^\circ$, subdomains are highlighted with green lines).}
  \label{fig:cryo:map}
\end{figure}

Within this context,
let us compare the behavior of the different transmission operators
on a helium vessel model consisting in a rectangular parallelepipedic
cavity interrupted with a cylindrical obstacle
and partitioned into $D=8$ subdomains.
The model is excited with its first spatial mode,
the length-to-wavelength ratio is chosen as $\ell/\lambda_w\simeq 12.5004$
and the cross-section is $h\times{}h$ with $h=\ell/2$.
This numerical experiment leads
to the convergence history displayed in \Fig{fig:cryo}
for each transmission operator.
It is clear from these data that the $\Sij{pade($24$)}{c}$
leads to the quickest convergence and requires only $50$ iterations to converge.
Compared with the $93$ iterations required by $\Sij{pade($24$)}{u}$,
a reduction of the iteration count of approximately $46\%$ is achieved.
\Fig{fig:cryo} shows as well the impact of a light regularization of
the \tc{pade}{c} and \tc{ml}{c} operators.
As expected, no performance gain is obtained,
as the unregularized counterparts already converge well.
For illustration purposes,
the computed field map is available in \Fig{fig:cryo:map}.

The above analysis would not be complete
without considering the wall clock time\footnote{The simulations
  of section~\ref{sec:cryo} were carried out on the \Code{NIC5} cluster
  hosted at the University of Li\`ege, Belgium.}.
As already discussed in section~\ref{sec:validation},
the novel transmission operators optimized for cavity problems
are associated with an increased computational cost,
when compared with their unbounded counterparts.
As a result, in the context of the considered acoustic study,
the \tc{oo0}{u} and \tc{oo2}{u} operators lead to the fastest computations
(\ie approximately $3.5$ and $3$ hours respectively),
despite their rather slow convergence.
In comparison with \tc{oo2}{u},
the rational operators see their wall clock time increased by a factor of
\begin{myenum*}
\item $2$ for \tc{pade}{u},
\item $2.7$ for \tc{ml}{c} and
\item $2.5$ for \tc{pade}{c}.
\end{myenum*}

Before concluding this section,
it is important to stress
that the current implementation does not treat the auxiliary unknowns
in the most efficient way.
For instance, it assembles all auxiliary unknowns into a single large system.
However, the extra unknowns stemming from the asymmetric nature of
$\Sij{ml($N$)/pade($N$)}{c}$ could be pulled out into smaller auxiliary systems,
which could improve the overall computational time.
Let us also mention that the OS iterative scheme can also be altered,
such that all auxiliary unknowns are decoupled from the subproblems,
as proposed in~\cite{Boubendir2018}.
Such improvements are left for a future work.

\section{Conclusion and final remarks}
\label{sec:conclusion}
In this work we presented new transmission operators
for the Schwarz method optimized for time-harmonic Helmholtz problems
in a rectangular cavity.
Those operators rely on different localizations of the Dirichlet-to-Neumann map
of this reference geometry.
Three different strategies were considered for devising local approximations
of the Fourier symbol of the $\DtN$, namely:
\begin{myenum*}
\item a zeroth-order Taylor approximation,
\item a truncated Mittag-Leffler partial fraction expansion and
\item a Pad\'e approximant.
\end{myenum*}
As these operators do not necessarily lead to a well defined scheme,
two regularization procedures where proposed to correct this problem.
Let us however mention that,
when combined with the modified variant of the Gram-Schmidt procedure
in the orthogonalization step of GMRES (without restart),
the proposed operators led to convergent iterative schemes
in all the numerical experiments considered in this work,
even in the unregularized case.

The new operators optimized for cavity problems
were also compared with operators optimized for unbounded geometries,
yet applied to cavities.
In the case of the reference rectangular geometry,
the gain in the iteration count brought by the novel rational operators
is clear when compared with the unbounded ones.
On the other hand, when considering geometries deviating from the canonical one,
this iteration count gain decreases.
Nonetheless, even if the gain can be modest for large deviations,
the new rational operators always led to a lower iteration count
(than their unbounded counterparts) in the considered numerical experiments.
Of course, this does not guarantee that the novel rational operators
will always perform better, but this suggests that they are sufficiently robust
with respect to geometrical changes.
The computational cost of the different transmission operators
was also briefly discussed,
as well as some possible improvements for reducing it.

To conclude this paper, let us draw some final remarks
regarding the selection of an operator.
While a detailed analysis is out of the scope of this work,
some general trends can be extracted from this study.
First of all, when considering the wall-clock time,
as the computational cost of the novel rational operators is high,
a simple \tc{oo0}{u} operator might be a good default choice.
Nonetheless, as already mentioned,
alterations of the Schwarz scheme could improve this aspect.
On the other hand, when focusing only on the iteration count,
the \tc{pade}{c} operator (with a small regularization)
can be recommended as a first guess.
Indeed, since evanescent modes in cavities and unbounded problems
converge quickly to each other,
both Pad\'e based approximations should behave similarly for
these modes if $N$ is sufficiently large.
Let us note that this is coherent with the spectra
discussed in section~\ref{sec:validation:spectrum}.
Thus, the difference between \tc{pade}{c} and \tc{pade}{u}
lies mainly in the non-evanescent part,
for which the \tc{pade}{u} is clearly a poor approximation.
Therefore, \emph{in the worst case scenario}, it is legitimate to expect that
both \tc{pade}{c} and \tc{pade}{u} will behave similarly,
at least with a GMRES without restart, as shown in our numerical experiments.
Let us note that this expectation still needs to be investigated more formally,
which is left for a future work.

\section*{Declaration of competing interest}
The authors declare that they have no known competing financial interests
or personal relationships
that could have appeared to influence the work reported in this paper.

\section*{Acknowledgments}
This research project has been funded by the Deutsche Forschungsgemeinschaft
(DFG, German Research Foundation) -- Project number 445906998.
The work of Nicolas Marsic is also supported by the Graduate School CE
within the Centre for Computational Engineering
at the Technische Universit\"at Darmstadt.
Computational resources have been provided by the
Consortium des \'Equipements de Calcul Intensif (C\'ECI),
funded by the Fonds de la Recherche Scientifique de Belgique (F.R.S.-FNRS)
under Grant No. 2.5020.11 and by the Walloon Region.
The authors would like to express their gratitude to Mr.~Anthony Royer
for his help with the \Code{GmshFEM} and \Code{GmshDDM} frameworks.
In addition, the authors are grateful to
Ms.~Heike Koch, Mr.~Achim Wagner, Mr.~Dragos Munteanu, Mr.~Christian Schmitt,
Dr.~Wolfgang F.O. M\"uller and Dr.~David Colignon
for the administrative and technical support.
Finally, the authors would like to thank the anonymous Reviewers,
whose comments improved significantly the quality of this work.

\appendix
\section{Calculation of the convergence radius -- non-overlapping case}
\label{sec:dtn}
In this appendix, let us briefly discuss the main steps of the calculation
leading to the Fourier symbol of the Dirichlet-to-Neumann map
shown in equation~\eqref{eq:close:lambda:dtn}.
To this end, we start with case $s^2 \neq k^2$,
and write the solutions
of~\eqref{eq:ddm:fourier:ode:0} and~\eqref{eq:ddm:fourier:ode:1},
together with the boundary conditions~\eqref{eq:ddm:fourier:bnd:0} and
\eqref{eq:ddm:fourier:bnd:1} and the definition~\eqref{eq:ddm:fourier:bnd:x0}.
Formally, we obtain:
\begin{subequations}
  \begin{align}[left = \empheqlbrace]
    \widehat{p}_0^{n+1}(x,s)
    & = P_0^{n+1}(s)\,
        \frac{\sinh\left[+\alpha\left(\frac{\ell}{2}+x\right)\right]}
             {\sinh\left[+\alpha\left(\frac{\ell}{2}+\gamma\right)\right]}
    & \forall x\in\left[-\frac{\ell}{2}, \gamma\right],
      \forall s\in\mathbb{S}, s^2\neq{}k^2,\\
    \widehat{p}_1^{n+1}(x,s)
    & = P_1^{n+1}(s)\,
        \frac{\sinh\left[-\alpha\left(\frac{\ell}{2}-x\right)\right]}
             {\sinh\left[-\alpha\left(\frac{\ell}{2}-\gamma\right)\right]}
    & \forall x\in\left[\gamma, +\frac{\ell}{2}\right], \
      \forall s\in\mathbb{S}, s^2\neq{}k^2,
  \end{align}
\end{subequations}
where
\begin{subequations}
  \begin{align}[left = {\alpha(s) = \empheqlbrace}]
    & -\imag\sqrt{k^2-s^2} & \text{if}~s^2<k^2,\\
    & 0                    & \text{if}~s^2=k^2,\\
    & \sqrt{s^2-k^2}       & \text{if}~s^2>k^2.
  \end{align}
\end{subequations}
This solution can then be derived with respect to $x$ and
the result evaluated at the interface between the subdomains,
\ie at $x = \gamma(t) = t\ell - \ell/2$ with $t\in[0, 1]$
as discussed in section~\ref{sec:problem}.
We thus have that
\begin{subequations}
  \begin{align}[left = \empheqlbrace]
    \left.\deriv{\widehat{p}_0^{n+1}}{x}(x, s)\right|_{x=\gamma(t)}
    & = \alpha\,P_0^{n+1}\,%
      \coth\mathopen{}\left[+\alpha{}t\ell\mathclose{}\right]
    & \forall s\in\mathbb{S}, s^2\neq{}k^2,
    \\
    \left.\deriv{\widehat{p}_1^{n+1}}{x}(x, s)\right|_{x=\gamma(t)}
    & = \alpha\,P_1^{n+1}\,%
      \coth\mathopen{}\left[-\alpha(1-t)\ell\mathclose{}\right]
    & \forall s\in\mathbb{S}, s^2\neq{}k^2.
  \end{align}
\end{subequations}
The convergence radius $\rho(s)$ is then obtained
by simplifying the transmission conditions~\eqref{eq:ddm:fourier:trsm:0}
and~\eqref{eq:ddm:fourier:trsm:1} with the above expressions.
In particular, we can write
\begin{subequations}
  \begin{align}[left = \empheqlbrace]
    P_0^{n+1}
    & = P_0^{n-1}
      \frac{\lambda_{10}-\alpha\coth\left[\alpha{}t\ell\right]}
           {\lambda_{10}+\alpha\coth\left[\alpha{}(1-t)\ell\right]}
      \frac{\lambda_{01}-\alpha\coth\left[\alpha{}(1-t)\ell\right]}
           {\lambda_{01}+\alpha\coth\left[\alpha{}t\ell\right]}
      = P_0^{n-1}\rho^2(s),
    \\
    P_1^{n+1}
    & = P_1^{n-1}
      \frac{\lambda_{01}-\alpha\coth\left[\alpha{}(1-t)\ell\right]}
           {\lambda_{01}+\alpha\coth\left[\alpha{}t\ell\right]}
      \frac{\lambda_{10}-\alpha\coth\left[\alpha{}t\ell\right]}
           {\lambda_{10}+\alpha\coth\left[\alpha{}(1-t)\ell\right]}
      = P_1^{n-1}\rho^2(s).
  \end{align}
\end{subequations}
Equation~\eqref{eq:close:lambda:dtn} is thus recovered,
for the case $s^2 \neq k^2$, by identifying the above terms
and by exploiting the definition of $\alpha(s)$, $\ell_{01}$ and $\ell_{10}$,
that is
\begin{subequations}
  \begin{align}[left = \empheqlbrace]
    \ell_{01} & = t\ell,\\
    \ell_{10} & = (1-t)\ell.
  \end{align}
\end{subequations}

Let us now treat the case $s^2 = k^2$.
In this case $\widehat{p}_i^{n+1}(x,s)$ writes
\begin{subequations}
  \begin{align}[left = \empheqlbrace]
    \widehat{p}_0^{n+1}(x,s)
    & = P_0^{n+1} \frac{x+\ell/2}{\gamma+\ell/2}
    & \forall x\in\left[-\frac{\ell}{2}, \gamma\right], s^2 = k^2,
    \\
    \widehat{p}_1^{n+1}(x,s)
    & = P_1^{n+1} \frac{x-\ell/2}{\gamma-\ell/2}
    & \forall x\in\left[\gamma, +\frac{\ell}{2}\right], s^2 = k^2,
  \end{align}
\end{subequations}
and its derivative evaluated at $x=\gamma$ is simply:
\begin{subequations}
  \begin{align}[left = \empheqlbrace]
    \left.\deriv{\widehat{p}_0^{n+1}}{x}(x, s)\right|_{x=\gamma(t)}
    & = +P_0^{n+1}\frac{1}{t\ell}
    & \forall s\in\mathbb{S}, s^2 = k^2,
    \\
    \left.\deriv{\widehat{p}_1^{n+1}}{x}(x, s)\right|_{x=\gamma(t)}
    & = -P_1^{n+1}\frac{1}{(1-t)\ell}
    & \forall s\in\mathbb{S}, s^2 = k^2.
  \end{align}
\end{subequations}
Therefore, the transmission conditions~\eqref{eq:ddm:fourier:trsm:0}
and~\eqref{eq:ddm:fourier:trsm:1} become
\begin{subequations}
  \begin{align}[left = \empheqlbrace]
    P_0^{n+1}
    & = P_0^{n-1}
      \frac{\lambda_{10}-\frac{1}{t\ell}}
           {\lambda_{10}+\frac{1}{(1-t)\ell}}
      \frac{\lambda_{01}-\frac{1}{(1-t)\ell}}
           {\lambda_{01}+\frac{1}{t\ell}}
      = P_0^{n-1}\rho^2(s=k),
    \\
    P_1^{n+1}
    & = P_1^{n-1}
      \frac{\lambda_{01}-\frac{1}{(1-t)\ell}}
           {\lambda_{01}+\frac{1}{t\ell}}
      \frac{\lambda_{10}-\frac{1}{t\ell}}
           {\lambda_{10}+\frac{1}{(1-t)\ell}}
      = P_1^{n-1}\rho^2(s=k),
  \end{align}
\end{subequations}
revealing thus the convergence radius $\rho(s=k)$.
Again, equation~\eqref{eq:close:lambda:dtn} is recovered,
for the case $s^2 = k^2$, by identifying the above terms
and by exploiting the definition of $\ell_{01}$ and $\ell_{10}$.

\section{Calculation of the convergence radius -- overlapping case}
\label{sec:dtn:overlap}
Let us now assume an overlap of $2\delta$ such that
$\Sigma_{01}$ is located at $\gamma_0 = \gamma+\delta$ and
$\Sigma_{10}$ is located at $\gamma_1 = \gamma-\delta$,
where $\gamma(t) = t\ell-\ell/2$ with $t\in[0,1]$ as in the previous appendix.
In the case $s^2 \neq k^2$, the overlapping variant
of~\eqref{eq:ddm:fourier:ode:0} and~\eqref{eq:ddm:fourier:ode:1},
together with the boundary conditions~\eqref{eq:ddm:fourier:bnd:0} and
\eqref{eq:ddm:fourier:bnd:1} and the definition~\eqref{eq:ddm:fourier:bnd:x0}
admits the following solutions:
\begin{subequations}
  \begin{align}[left = \empheqlbrace]
    \widehat{p}_0^{n+1}(x,s)
    & = P_0^{n+1}(s)\,
        \frac{\sinh\left[+\alpha\left(\frac{\ell}{2}+x\right)\right]}
             {\sinh\left[+\alpha\left(\frac{\ell}{2}+\gamma_0\right)\right]}
    & \forall x\in\left[-\frac{\ell}{2}, \gamma_0\right],
      \forall s\in\mathbb{S}, s^2\neq{}k^2,\\
    \widehat{p}_1^{n+1}(x,s)
    & = P_1^{n+1}(s)\,
        \frac{\sinh\left[-\alpha\left(\frac{\ell}{2}-x\right)\right]}
             {\sinh\left[-\alpha\left(\frac{\ell}{2}-\gamma_1\right)\right]}
    & \forall x\in\left[\gamma_1, +\frac{\ell}{2}\right], \
      \forall s\in\mathbb{S}, s^2\neq{}k^2,
  \end{align}
\end{subequations}
which is similar as for the non-overlapping case
but with interfaces located at different positions,
$x=\gamma_0$ and $x=\gamma_1$, instead of a unique one, $x=\gamma$.
After evaluating those two solutions at $x=\gamma_0$ and $x=\gamma_1$,
we obtain:
\begin{subequations}
  \begin{align}[left = \empheqlbrace]
    \widehat{p}_0^{n+1}(\gamma_0,s)
    & = P_0^{n+1}(s)
    & \forall s\in\mathbb{S}, s^2\neq{}k^2,\\
    \widehat{p}_1^{n+1}(\gamma_0,s)
    & = P_1^{n+1}(s)\,
        \frac{\sinh\left[\alpha\left(\frac{\ell}{2}-\gamma-\delta\right)\right]}
             {\sinh\left[\alpha\left(\frac{\ell}{2}-\gamma+\delta\right)\right]}
    & \forall s\in\mathbb{S}, s^2\neq{}k^2,
  \end{align}
\end{subequations}
and
\begin{subequations}
  \begin{align}[left = \empheqlbrace]
    \widehat{p}_0^{n+1}(\gamma_1,s)
    & = P_0^{n+1}(s)
        \frac{\sinh\left[\alpha\left(\frac{\ell}{2}+\gamma-\delta\right)\right]}
             {\sinh\left[\alpha\left(\frac{\ell}{2}+\gamma+\delta\right)\right]}
    & \forall s\in\mathbb{S}, s^2\neq{}k^2,\\
    \widehat{p}_1^{n+1}(\gamma_1,s)
    & = P_1^{n+1}(s)
    & \forall s\in\mathbb{S}, s^2\neq{}k^2.
  \end{align}
\end{subequations}
Regarding the derivatives with respect to $x$ of
$\widehat{p}_0^{n+1}(x,s)$ and $\widehat{p}_0^{n+1}(x,s)$
at $x=\gamma_0$ and $x=\gamma_1$, we have:
\begin{subequations}
  \begin{align}[left = \empheqlbrace]
    \left.\deriv{\widehat{p}_0^{n+1}}{x}(x, s)\right|_{x=\gamma_0}
    & = +\alpha\,P_0^{n+1}\,%
        \frac{\cosh\left[\alpha\left(\frac{\ell}{2}+\gamma+\delta\right)\right]}
             {\sinh\left[\alpha\left(\frac{\ell}{2}+\gamma+\delta\right)\right]}
    & \forall s\in\mathbb{S}, s^2\neq{}k^2,
    \\
    \left.\deriv{\widehat{p}_1^{n+1}}{x}(x, s)\right|_{x=\gamma_0}
    & = -\alpha\,P_1^{n+1}\,%
        \frac{\cosh\left[\alpha\left(\frac{\ell}{2}-\gamma-\delta\right)\right]}
             {\sinh\left[\alpha\left(\frac{\ell}{2}-\gamma+\delta\right)\right]}
    & \forall s\in\mathbb{S}, s^2\neq{}k^2,
  \end{align}
\end{subequations}
and
\begin{subequations}
  \begin{align}[left = \empheqlbrace]
    \left.\deriv{\widehat{p}_0^{n+1}}{x}(x, s)\right|_{x=\gamma_1}
    & = +\alpha\,P_0^{n+1}\,%
        \frac{\cosh\left[\alpha\left(\frac{\ell}{2}+\gamma-\delta\right)\right]}
             {\sinh\left[\alpha\left(\frac{\ell}{2}+\gamma+\delta\right)\right]}
    & \forall s\in\mathbb{S}, s^2\neq{}k^2,
    \\
    \left.\deriv{\widehat{p}_1^{n+1}}{x}(x, s)\right|_{x=\gamma_1}
    & = -\alpha\,P_1^{n+1}\,%
        \frac{\cosh\left[\alpha\left(\frac{\ell}{2}-\gamma+\delta\right)\right]}
             {\sinh\left[\alpha\left(\frac{\ell}{2}-\gamma+\delta\right)\right]}
    & \forall s\in\mathbb{S}, s^2\neq{}k^2.
  \end{align}
\end{subequations}
Therefore, the overlapping variant of
the transmission conditions~\eqref{eq:ddm:fourier:trsm:0}
and~\eqref{eq:ddm:fourier:trsm:1} leads to
\begin{subequations}
  \begin{align}[left = \empheqlbrace]
    P_0^{n+1}
    & = P_0^{n-1}\,
      \frac{\lambda_{10}-\alpha\coth\left[\alpha\left(\frac{\ell}{2}+\gamma-\delta\right)\right]}
           {\lambda_{10}+\alpha\coth\left[\alpha\left(\frac{\ell}{2}-\gamma+\delta\right)\right]}
      \frac{\lambda_{01}-\alpha\coth\left[\alpha\left(\frac{\ell}{2}-\gamma-\delta\right)\right]}
           {\lambda_{01}+\alpha\coth\left[\alpha\left(\frac{\ell}{2}+\gamma+\delta\right)\right]}\nonumber\\
    & \hspace{4.5cm}
      \frac{\sinh\left[\alpha\left(\frac{\ell}{2}+\gamma-\delta\right)\right]}
           {\sinh\left[\alpha\left(\frac{\ell}{2}+\gamma+\delta\right)\right]}
      \frac{\sinh\left[\alpha\left(\frac{\ell}{2}-\gamma-\delta\right)\right]}
           {\sinh\left[\alpha\left(\frac{\ell}{2}-\gamma+\delta\right)\right]}
      = P_0^{n-1}\rho^2,\\
    P_1^{n+1}
    & = P_1^{n-1}\,
      \frac{\lambda_{01}-\alpha\coth\left[\alpha\left(\frac{\ell}{2}-\gamma-\delta\right)\right]}
           {\lambda_{01}+\alpha\coth\left[\alpha\left(\frac{\ell}{2}+\gamma+\delta\right)\right]}
      \frac{\lambda_{10}-\alpha\coth\left[\alpha\left(\frac{\ell}{2}+\gamma-\delta\right)\right]}
           {\lambda_{10}+\alpha\coth\left[\alpha\left(\frac{\ell}{2}-\gamma+\delta\right)\right]}\nonumber\\
    & \hspace{4.5cm}
      \frac{\sinh\left[\alpha\left(\frac{\ell}{2}-\gamma-\delta\right)\right]}
           {\sinh\left[\alpha\left(\frac{\ell}{2}-\gamma+\delta\right)\right]}
      \frac{\sinh\left[\alpha\left(\frac{\ell}{2}+\gamma-\delta\right)\right]}
           {\sinh\left[\alpha\left(\frac{\ell}{2}+\gamma+\delta\right)\right]}
      = P_1^{n-1}\rho^2.
  \end{align}
\end{subequations}
Equation~\eqref{eq:fourier:rho:overlap} is thus recovered,
for the case $s^2 \neq k^2$,
by identifying the above terms and by exploiting the definitions of
$\alpha(s)$, $\ell_{01}$, $\ell_{10}$, $\ell_{01}^\prime$ and
$\ell_{10}^\prime$, that are:
\begin{subequations}
  \begin{align}[left = \empheqlbrace]
    \ell_{01}        & = \frac{\ell}{2} + \gamma + \delta,\\
    \ell_{10}        & = \frac{\ell}{2} - \gamma + \delta,\\
    \ell_{01}^\prime & = \frac{\ell}{2} + \gamma - \delta,\\
    \ell_{10}^\prime & = \frac{\ell}{2} - \gamma - \delta.
  \end{align}
\end{subequations}

Finally, the case $s^2 = k^2$ is obtained in a similar way,
which leads to following recursion
\begin{subequations}
  \begin{align}[left = \empheqlbrace]
    P_0^{n+1}
    & = P_0^{n-1}\,
        \frac{\lambda_{10}-\left(\frac{\ell}{2}+\gamma-\delta\right)^{-1}}
             {\lambda_{10}+\left(\frac{\ell}{2}-\gamma+\delta\right)^{-1}}
        \frac{\lambda_{01}-\left(\frac{\ell}{2}-\gamma-\delta\right)^{-1}}
             {\lambda_{01}+\left(\frac{\ell}{2}+\gamma+\delta\right)^{-1}}
        \frac{\frac{\ell}{2}+\gamma-\delta}
             {\frac{\ell}{2}+\gamma+\delta}
        \frac{\frac{\ell}{2}-\gamma-\delta}
             {\frac{\ell}{2}-\gamma+\delta}
      = P_0^{n-1}\rho^2,\\
    P_1^{n+1}
    & = P_1^{n-1}\,
        \frac{\lambda_{01}-\left(\frac{\ell}{2}-\gamma-\delta\right)^{-1}}
             {\lambda_{01}+\left(\frac{\ell}{2}+\gamma+\delta\right)^{-1}}
        \frac{\lambda_{10}-\left(\frac{\ell}{2}+\gamma-\delta\right)^{-1}}
             {\lambda_{10}+\left(\frac{\ell}{2}-\gamma+\delta\right)^{-1}}
        \frac{{\ell}{2}-\gamma-\delta}
             {{\ell}{2}-\gamma+\delta}
        \frac{{\ell}{2}+\gamma-\delta}
             {{\ell}{2}+\gamma+\delta}
      = P_1^{n-1}\rho^2,
  \end{align}
\end{subequations}
revealing thus the convergence radius.
Equation~\eqref{eq:fourier:rho:overlap} is again recovered,
for the case $s^2 = k^2$,
by identifying the above terms and by exploiting the definitions of
$\alpha(s)$, $\ell_{01}$, $\ell_{10}$, $\ell_{01}^\prime$ and
$\ell_{10}^\prime$.

\newpage
\section{Tables with raw numbers}
\label{sec:raw}
In this section, we compiled the raw numbers of the following experiments.
\begin{itemize}
\item Increase in the number of rectangular subdomains
  (section~\ref{sec:val:dmany}, \Fig{fig:val:dmany}).
  \begin{table}[ht]
    \centering
    \begin{tabular}{cccccccc}
      \toprule
      $D$ &
      $\Sij{oo0}{u}$ & $\Sij{oo2}{u}$      & $\Sij{pade($64$)}{u}$ &
      $\Sij{oo0}{c}$ & $\Sij{ml($64$)}{c}$ &  $\Sij{ml($128$)}{c}$ &
                                              $\Sij{pade($64$)}{c}$ \\
      \midrule
       2 &  186 &   71 &   63 &  138 &  32 &  25 &  7 \\
       4 &  367 &  181 &  163 & 1063 &  46 &  35 &  8 \\
       6 &  527 &  278 &  262 & 1316 &  67 &  51 & 10 \\
       8 &  679 &  389 &  362 & 1584 &  79 &  61 & 14 \\
      10 &  825 &  488 &  462 & 2580 & 103 &  84 & 18 \\
      12 &  980 &  582 &  562 & 2241 & 161 & 122 & 20 \\
      14 & 1098 &  674 &  663 & 2586 & 216 & 161 & 26 \\
      16 & 1251 &  787 &  763 & 2950 & 221 & 155 & 23 \\
      18 & 1393 &  871 &  863 & 3322 & 342 & 264 & 34 \\
      20 & 1526 &  974 &  964 & 7160 & 293 & 222 & 33 \\
      22 & 1650 & 1079 & 1064 & 3985 & 369 & 343 & 42 \\
      24 & 1780 & 1180 & 1165 & 4447 & 309 & 229 & 27 \\
      26 & 1919 & 1260 & 1266 & 4825 & 459 & 416 & 50 \\
      28 & 2029 & 1339 & 1366 & 5283 & 492 & 386 & 54 \\
      30 & 2157 & 1428 & 1467 & 9376 & 479 & 347 & 42 \\
      32 & 2271 & 1527 & 1568 & 6180 & 686 & 508 & 62 \\
      \bottomrule
    \end{tabular}
    \caption{Number of GMRES iterations
      as a function of the number of subdomains.}
    \label{tab:val:dmany}
  \end{table}

\item Impact of the length-to-wavelength ratio
  (section~\ref{sec:val:k}, \Fig{fig:val:k}).
  \begin{table}[ht]
    \centering
    \begin{tabular}{ccccccc}
      \toprule
      $\ell/\lambda_w$ &
      $\Sij{oo0}{u}$ & $\Sij{oo2}{u}$      & $\Sij{pade($64$)}{u}$ &
      $\Sij{oo0}{c}$ & $\Sij{ml($64$)}{c}$ & $\Sij{pade($64$)}{c}$ \\
      \midrule
       5.0009 & 196 &  94 &  81 &  357 & 29 & 12 \\
       7.0009 & 181 & 109 & 102 &  423 & 31 & 12 \\
       9.0009 & 239 & 142 & 132 &  550 & 35 & 11 \\
      11.0009 & 267 & 166 & 158 &  662 & 44 & 12 \\
      13.0009 & 379 & 206 & 191 &  831 & 44 & 14 \\
      15.0009 & 413 & 233 & 218 &  919 & 48 & 14 \\
      17.0009 & 473 & 268 & 249 & 1078 & 79 & 14 \\
      19.0009 & 435 & 279 & 269 & 1136 & 64 & 14 \\
      21.0009 & 491 & 310 & 298 & 1271 & 65 & 12 \\
      23.0009 & 523 & 340 & 326 & 1380 & 95 & 14 \\
      25.0009 & 677 & 389 & 362 & 1583 & 78 & 14 \\
      \bottomrule
    \end{tabular}
    \caption{Number of GMRES iterations as a function
      of the length-to-wavelength ratio for a cavity with $D=8$ subdomains.}
    \label{tab:val:k}
  \end{table}

\newpage
\item Trapezoidal geometry -- modified Gram-Schmidt variant
  (section~\ref{sec:sensitivity:trapeze:mgs}, \Fig{fig:trapeze:data}).
  \begin{table}[ht]
    \centering
    \begin{tabular}{cccccccc}
      \toprule
      &
      \multicolumn{1}{c}{\tc{pade($64$)}{u}} &
      \multicolumn{3}{c}{  \tc{ml($64$)}{c}} &
      \multicolumn{3}{c}{\tc{pade($64$)}{c}} \\
      \cmidrule(r){2-2}\cmidrule(r){3-5}\cmidrule(r){6-8}
      $\delta$ &
      Unreg. &
      Unreg. & r($0.1$) & m($0.9, 64$) &
      Unreg. & r($0.1$) & m($0.9, 64$) \\
      \midrule
      0.00 &  762 & 220 & 330 & 286 &  18 & 241 & 266 \\
      0.01 &  759 & 190 & 241 & 227 & 146 & 218 & 216 \\
      0.02 &  760 & 234 & 297 & 294 & 198 & 275 & 286 \\
      0.04 &  759 & 321 & 375 & 373 & 241 & 351 & 357 \\
      0.08 &  774 & 342 & 378 & 378 & 305 & 370 & 383 \\
      0.16 &  765 & 414 & 401 & 425 & 377 & 406 & 434 \\
      0.32 &  870 & 498 & 536 & 535 & 435 & 521 & 543 \\
      0.64 & 1008 & 657 & 650 & 679 & 540 & 638 & 679 \\
      1.28 & 1243 & 808 & 743 & 786 & 611 & 733 & 788 \\
      \bottomrule
    \end{tabular}
    \caption{Iteration count of the GMRES solver
      when used on a trapezoidal geometry partitioned into $D=16$ subdomains.}
    \label{tab:trapeze:data}
  \end{table}

\item Rectangular cavity involving obstacles
  -- impact of the number of obstacles
  (section~\ref{sec:sensitivity:obstacles:N}, \Fig{fig:obstacles:N}).
  \begin{table}[ht]
    \centering
    \begin{tabular}{cccccccc}
      \toprule
      &
      \multicolumn{1}{c}{\tc{pade($64$)}{u}} &
      \multicolumn{3}{c}{  \tc{ml($64$)}{c}} &
      \multicolumn{3}{c}{\tc{pade($64$)}{c}} \\
      \cmidrule(r){2-2}\cmidrule(r){3-5}\cmidrule(r){6-8}
      $O$ &
      Unreg. &
      Unreg. & r($0.1$) & m($0.5, 4$) &
      Unreg. & r($0.1$) & m($0.5, 4$) \\
      \midrule
       1 & 420 & 183 &  21 & 387 & 399 & 230 & 218 \\
       3 & 462 & 202 &  76 & 415 & 427 & 264 & 278 \\
       5 & 472 & 237 & 130 & 405 & 423 & 284 & 329 \\
       7 & 562 & 326 & 220 & 465 & 476 & 366 & 413 \\
       9 & 578 & 323 & 300 & 476 & 489 & 393 & 459 \\
      11 & 498 & 323 & 388 & 439 & 457 & 381 & 452 \\
      13 & 567 & 393 & 480 & 472 & 488 & 457 & 525 \\
      15 & 647 & 471 & 584 & 520 & 538 & 535 & 595 \\
      17 & 662 & 542 & 606 & 578 & 568 & 589 & 609 \\
      \bottomrule
    \end{tabular}
    \caption{Iteration count of the GMRES solver
      when used on a rectangular cavity involving $O$ obstacles
      of radius $R=0.5\lambda_w$ and $D=17$ subdomains.}
    \label{tab:obstacles:N}
  \end{table}

\newpage
\item Rectangular cavity involving obstacles
  -- impact of the size of the obstacles
  (section~\ref{sec:sensitivity:obstacles:R}, \Fig{fig:obstacles:R}).
  \begin{table}[ht]
    \centering
    \begin{tabular}{cccccccc}
      \toprule
      &
      \multicolumn{1}{c}{\tc{pade($64$)}{u}} &
      \multicolumn{3}{c}{  \tc{ml($64$)}{c}} &
      \multicolumn{3}{c}{\tc{pade($64$)}{c}} \\
      \cmidrule(r){2-2}\cmidrule(r){3-5}\cmidrule(r){6-8}
      $R$ &
      Unreg. &
      Unreg. & r($0.1$) & m($0.5, 4$) &
      Unreg. & r($0.1$) & m($0.5, 4$) \\
      \midrule
      $0.025\lambda_w$ & 435 & 175 &  84 & 374 & 387 & 232 & 239 \\
      $0.050\lambda_w$ & 454 & 206 &  88 & 390 & 412 & 256 & 268 \\
      $0.100\lambda_w$ & 412 & 199 &  96 & 352 & 374 & 230 & 250 \\
      $0.200\lambda_w$ & 452 & 200 & 126 & 377 & 404 & 270 & 303 \\
      $0.400\lambda_w$ & 536 & 247 & 184 & 438 & 456 & 318 & 342 \\
      \bottomrule
    \end{tabular}
    \caption{Iteration count of the GMRES solver
      when used on a rectangular cavity involving $O=7$ obstacles
      and $D=17$ subdomains.}
    \label{tab:obstacles:R}
  \end{table}
\end{itemize}

\bibliographystyle{elsarticle-num}
\bibliography{ddm.bib}

\begin{thebibliography}{10}
\expandafter\ifx\csname url\endcsname\relax
  \def\url#1{\texttt{#1}}\fi
\expandafter\ifx\csname urlprefix\endcsname\relax\def\urlprefix{URL }\fi
\expandafter\ifx\csname href\endcsname\relax
  \def\href#1#2{#2} \def\path#1{#1}\fi

\bibitem{Ihlenburg1995}
F.~Ihlenburg, I.~Babu\v{s}ka, Finite element solution of the {H}elmholtz
  equation with high wave number part {I}: The h-version of the {F}{E}{M},
  Computers \& Mathematics with Applications 30~(9) (1995) 9--37.
\newblock \href {https://doi.org/10.1016/0898-1221(95)00144-N}
  {\path{doi:10.1016/0898-1221(95)00144-N}}.

\bibitem{Ernst2012}
O.~G. Ernst, M.~J. Gander, Why it is difficult to solve {H}elmholtz problems
  with classical iterative methods, in: I.~G. Graham, T.~Y. Hou, O.~Lakkis,
  R.~Scheichl (Eds.), Numerical Analysis of Multiscale Problems, Vol.~83 of
  Lecture Notes in Computational Science and Engineering, 2012, pp. 325--363.
\newblock \href {https://doi.org/10.1007/978-3-642-22061-6_10}
  {\path{doi:10.1007/978-3-642-22061-6_10}}.

\bibitem{Ihlenburg1997}
F.~Ihlenburg, I.~Babu\v{s}ka, Finite element solution of the {H}elmholtz
  equation with high wave number part {I}{I}: The h-p version of the {F}{E}{M},
  SIAM Journal on Numerical Analysis 34~(1) (1997) 315--358.
\newblock \href {https://doi.org/10.1137/S0036142994272337}
  {\path{doi:10.1137/S0036142994272337}}.

\bibitem{Moiola2014}
A.~Moiola, E.~A. Spence, Is the {H}elmholtz equation really sign-indefinite?,
  SIAM Review 56~(2) (2014) 274--312.
\newblock \href {https://doi.org/10.1137/120901301}
  {\path{doi:10.1137/120901301}}.

\bibitem{Diwan2019}
G.~C. Diwan, A.~Moiola, E.~A. Spence, Can coercive formulations lead to fast
  and accurate solution of the {H}elmholtz equation?, Journal of Computational
  and Applied Mathematics 352 (2019) 110--131.
\newblock \href {https://doi.org/10.1016/j.cam.2018.11.035}
  {\path{doi:10.1016/j.cam.2018.11.035}}.

\bibitem{Yannakakis1981}
M.~Yannakakis, Computing the minimum fill-in is {N}{P}-complete, SIAM Journal
  on Algebraic Discrete Methods 2~(1) (1981) 77--79.
\newblock \href {https://doi.org/10.1137/0602010} {\path{doi:10.1137/0602010}}.

\bibitem{Marsic2018b}
N.~Marsic, H.~De~Gersem, G.~Dem\'esy, A.~Nicolet, C.~Geuzaine, Modal analysis
  of the ultrahigh finesse {H}aroche {Q}{E}{D} cavity, New Journal of Physics
  20~(4) (2018) 043058.
\newblock \href {https://doi.org/10.1088/1367-2630/aab6fd}
  {\path{doi:10.1088/1367-2630/aab6fd}}.

\bibitem{Despres1990}
B.~Despr\'es,
  \href{https://gallica.bnf.fr/ark:/12148/bpt6k57815213}{D\'ecomposition de
  domaine et probl\`eme de {H}elmholtz}, Comptes Rendus de l'Acad\'emie des
  Sciences 311 (1990) 313--316.
\newline\urlprefix\url{https://gallica.bnf.fr/ark:/12148/bpt6k57815213}

\bibitem{Boubendir2007}
Y.~Boubendir, An analysis of the {BEM}-{FEM} non-overlapping domain
  decomposition method for a scattering problem, Journal of Computational and
  Applied Mathematics 204~(2) (2007) 282--291.
\newblock \href {https://doi.org/10.1016/j.cam.2006.02.044}
  {\path{doi:10.1016/j.cam.2006.02.044}}.

\bibitem{Gander2002}
M.~J. Gander, F.~Magoul\`es, F.~Nataf, Optimized {S}chwarz methods without
  overlap for the {H}elmholtz equation, SIAM Journal on Scientific Computing
  24~(1) (2002) 38--60.
\newblock \href {https://doi.org/10.1137/S1064827501387012}
  {\path{doi:10.1137/S1064827501387012}}.

\bibitem{Boubendir2012}
Y.~Boubendir, X.~Antoine, C.~Geuzaine, A quasi-optimal non-overlapping domain
  decomposition algorithm for the {H}elmholtz equation, Journal of
  Computational Physics 231~(2) (2012) 262--280.
\newblock \href {https://doi.org/10.1016/j.jcp.2011.08.007}
  {\path{doi:10.1016/j.jcp.2011.08.007}}.

\bibitem{Vion2014a}
A.~Vion, C.~Geuzaine, Double sweep preconditioner for optimized {S}chwarz
  methods applied to the {H}elmholtz problem, Journal of Computational Physics
  266 (2014) 171--190.
\newblock \href {https://doi.org/10.1016/j.jcp.2014.02.015}
  {\path{doi:10.1016/j.jcp.2014.02.015}}.

\bibitem{Gander2019}
M.~J. Gander, H.~Zhang, A class of iterative solvers for the {H}elmholtz
  equation: Factorizations, sweeping preconditioners, source transfer, single
  layer potentials, polarized traces, and optimized {S}chwarz methods, SIAM
  Review 61~(1) (2019) 3--76.
\newblock \href {https://doi.org/10.1137/16m109781x}
  {\path{doi:10.1137/16m109781x}}.

\bibitem{Dolean2015a}
V.~Dolean, P.~Jolivet, F.~Nataf, An introduction to domain decomposition
  methods: algorithms, theory and parallel implementation, Society for
  Industrial and Applied Mathematics, 2015.
\newblock \href {https://doi.org/10.1137/1.9781611974065}
  {\path{doi:10.1137/1.9781611974065}}.

\bibitem{Peng2011}
Z.~Peng, J.-F. Lee, Non-conformal domain decomposition method with mixed true
  second order transmission condition for solving large finite antenna arrays,
  IEEE Transactions on Antennas and Propagation 59~(5) (2011) 1638--1651.
\newblock \href {https://doi.org/10.1109/TAP.2011.2123067}
  {\path{doi:10.1109/TAP.2011.2123067}}.

\bibitem{Tournier2017}
P.-H. Tournier, M.~Bonazzoli, V.~Dolean, F.~Rapetti, F.~Hecht, F.~Nataf,
  I.~Aliferis, I.~El~Kanfoud, C.~Migliaccio, M.~de~Buhan, M.~Darbas,
  S.~Semenov, C.~Pichot, Numerical modeling and high-speed parallel computing:
  New perspectives on tomographic microwave imaging for brain stroke detection
  and monitoring., IEEE Antennas and Propagation Magazine 59~(5) (2017)
  98--110.
\newblock \href {https://doi.org/10.1109/map.2017.2731199}
  {\path{doi:10.1109/map.2017.2731199}}.

\bibitem{Marsic2016a}
N.~Marsic, C.~Waltz, J.-F. Lee, C.~Geuzaine, Domain decomposition methods for
  time-harmonic electromagnetic waves with high order whitney forms, IEEE
  Transactions on Magnetics 52~(3) (2016) 1--4.
\newblock \href {https://doi.org/10.1109/TMAG.2015.2476510}
  {\path{doi:10.1109/TMAG.2015.2476510}}.

\bibitem{Saleh2007}
B.~E.~A. Saleh, M.~C. Teich, Fundamentals of Photonics, 2nd Edition,
  Wiley-Interscience, 2007.
\newblock \href {https://doi.org/10.1002/0471213748}
  {\path{doi:10.1002/0471213748}}.

\bibitem{Ko2006}
K.~Ko, N.~Folwell, L.~Ge, A.~Guetz, L.~Lee, Z.~Li, C.~Ng, E.~Prudencio,
  G.~Schussman, R.~Uplenchwar, L.~Xiao, Advances in electromagnetic modelling
  through high performance computing, Physica C: Superconductivity and its
  Applications 441~(1-2) (2006) 258--262.
\newblock \href {https://doi.org/10.1016/j.physc.2006.03.139}
  {\path{doi:10.1016/j.physc.2006.03.139}}.

\bibitem{Ahlfors1979}
L.~V. Ahlfors, Complex analysis, 3rd Edition, McGraw-Hill, New York, NY, 1979.

\bibitem{Baker1996}
G.~A. Baker, P.~Graves-Morris, Padé Approximants, 2nd Edition, Cambridge
  University Press, 1996.
\newblock \href {https://doi.org/10.1017/cbo9780511530074}
  {\path{doi:10.1017/cbo9780511530074}}.

\bibitem{Oldham2009}
K.~Oldham, J.~Myland, J.~Spanier, An Atlas of Functions, 2nd Edition,
  Springer-Verlag New York, 2009.
\newblock \href {https://doi.org/10.1007/978-0-387-48807-3}
  {\path{doi:10.1007/978-0-387-48807-3}}.

\bibitem{Bini2014}
D.~A. Bini, L.~Robol, Solving secular and polynomial equations: A
  multiprecision algorithm, Journal of Computational and Applied Mathematics
  272 (2014) 276--292.
\newblock \href {https://doi.org/10.1016/j.cam.2013.04.037}
  {\path{doi:10.1016/j.cam.2013.04.037}}.

\bibitem{Milinazzo1997}
F.~A. Milinazzo, C.~A. Zala, G.~H. Brooke, Rational square-root approximations
  for parabolic equation algorithms, The Journal of the Acoustical Society of
  America 101~(2) (1997) 760--766.
\newblock \href {https://doi.org/10.1121/1.418038}
  {\path{doi:10.1121/1.418038}}.

\bibitem{Titchmarsh1976}
E.~C. Titchmarsh, The theory of functions, 2nd Edition, Oxford University
  Press, London, England, 1976.

\bibitem{Royer2021}
A.~Royer, E.~B\'echet, C.~Geuzaine, Gmsh-{F}em: An efficient finite element
  library based on {G}msh, in: 14th WCCM-ECCOMAS Congress, 2021.
\newblock \href {https://doi.org/10.23967/wccm-eccomas.2020.161}
  {\path{doi:10.23967/wccm-eccomas.2020.161}}.

\bibitem{Saad2003}
Y.~Saad, Iterative methods for sparse linear systems, 2nd Edition, Society for
  Industrial and Applied Mathematics, 2003.
\newblock \href {https://doi.org/10.1137/1.9780898718003}
  {\path{doi:10.1137/1.9780898718003}}.

\bibitem{Saad1986}
Y.~Saad, M.~H. Schultz, {GMRES}: a generalized minimal residual algorithm for
  solving nonsymmetric linear systems, SIAM Journal on Scientific and
  Statistical Computing 7~(3) (1986) 856--869.
\newblock \href {https://doi.org/10.1137/0907058} {\path{doi:10.1137/0907058}}.

\bibitem{Balay1997}
S.~Balay, W.~D. Gropp, L.~C. McInnes, B.~F. Smith, Efficient management of
  parallelism in object oriented numerical software libraries, in: E.~Arge,
  A.~M. Bruaset, H.~P. Langtangen (Eds.), Modern Software Tools in Scientific
  Computing, Birkh\"auser Press, 1997, pp. 163--202.
\newblock \href {https://doi.org/10.1007/978-1-4612-1986-6_8}
  {\path{doi:10.1007/978-1-4612-1986-6_8}}.

\bibitem{Amestoy2001}
P.~R. Amestoy, I.~S. Duff, J.~Koster, J.-Y. L'Excellent, A fully asynchronous
  multifrontal solver using distributed dynamic scheduling, SIAM Journal on
  Matrix Analysis and Applications 23~(1) (2001) 15--41.
\newblock \href {https://doi.org/10.1137/S0895479899358194}
  {\path{doi:10.1137/S0895479899358194}}.

\bibitem{Greenbaum1996}
A.~Greenbaum, V.~Pták, Z.~Strakoš, Any nonincreasing convergence curve is
  possible for {GMRES}, SIAM Journal on Matrix Analysis and Applications 17~(3)
  (1996) 465--469.
\newblock \href {https://doi.org/10.1137/s0895479894275030}
  {\path{doi:10.1137/s0895479894275030}}.

\bibitem{Haider2019}
D.~M. Haider, H.~De~Gersem, J.~Golm, T.~Koettig, F.~Kurian, N.~Marsic, W.~F.~O.
  M\"uller, M.~Schmelz, M.~Schwickert, T.~Sieber, R.~Stolz, T.~St\"ohlker,
  V.~Tympel, F.~Ucar, V.~Zakosarenko, Versatile beamline cryostat for the
  cryogenic current comparator ({CCC}) for {FAIR}, in: Proceedings of the 8th
  International Beam Instrumentation Conference ({IBIC'19}), no.~8 in
  International Beam Instrumentation Conference, JACoW Publishing, Geneva,
  Switzerland, 2019, pp. 78--81.
\newblock \href {https://doi.org/10.18429/JACoW-IBIC2019-MOPP007}
  {\path{doi:10.18429/JACoW-IBIC2019-MOPP007}}.

\bibitem{Seidel2018}
P.~Seidel, V.~Tympel, R.~Neubert, J.~Golm, M.~Schmelz, R.~Stolz,
  V.~Zakosarenko, T.~Sieber, M.~Schwickert, F.~Kurian, F.~Schmidl,
  T.~St\"ohlker, Cryogenic current comparators for larger beamlines, IEEE
  Transactions on Applied Superconductivity 28~(4) (2018) 1--5.
\newblock \href {https://doi.org/10.1109/tasc.2018.2815647}
  {\path{doi:10.1109/tasc.2018.2815647}}.

\bibitem{Boubendir2018}
Y.~Boubendir, D.~Midura, Non-overlapping domain decomposition algorithm based
  on modified transmission conditions for the {H}elmholtz equation, Computers
  \& Mathematics with Applications 75~(6) (2018) 1900--1911.
\newblock \href {https://doi.org/10.1016/j.camwa.2017.07.027}
  {\path{doi:10.1016/j.camwa.2017.07.027}}.

\end{thebibliography}
\end{document}